# Issues in the Estimation of Mis-Specified Models of Fractionally Integrated Processes[†]

Gael M. Martin,[‡] K. Nadarajah & D. S. Poskitt
*Department of Econometrics & Business Statistics, Monash University*

September 18, 2018

**Abstract**

We provide a comprehensive set of new results on the impact of mis-specifying the short run dynamics in fractionally integrated processes. We show that four alternative parametric estimators – frequency domain maximum likelihood, Whittle, time domain maximum likelihood and conditional sum of squares – converge to the same pseudo-true value under common mis-specification, and that they possess a common asymptotic distribution. The results are derived assuming a completely general parametric specification for the short run dynamics of the estimated (mis-specified) fractional model, and with long memory, short memory and antipersistence in both the model and the true data generating process accommodated. As well as providing new theoretical insights, we undertake an extensive set of numerical explorations, beginning with the numerical evaluation, and implementation, of the (common) asymptotic distribution that holds under the most extreme form of mis-specification. Simulation experiments are then conducted to assess the relative finite sample performance of all four mis-specified estimators, initially under the assumption of a known mean, as accords with the theoretical derivations. The importance of the known mean assumption is illustrated via the production of an alternative set of bias and mean squared error results, in which the estimators are applied to demeaned data. The paper concludes with a discussion of open problems.

*Keywords*: Long memory models, pseudo-true parameter, frequency domain estimators, time domain estimators, Whittle, conditional sum of squares.

*MSC2010 subject classifications*: Primary 62M10, 62M15; Secondary 62G09
*JEL classifications*: C18, C22, C52

## 1 Introduction

Let $\{y_t\}$, $t \in \mathbb{Z}$, be a (strictly) stationary process with mean $\mu_0$ and spectral density $f_0(\lambda)$, $\lambda \in [-\pi, \pi]$, that is such that

$$f_0(\lambda) \sim |\lambda|^{-2d_0} L_0(\lambda) \quad \text{as} \quad \lambda \to 0,$$

where $0 \leq |d_0| < 0.5$ and $L_0(\lambda)$ is a positive function that is slowly varying at 0. Prototypical examples of processes of this type are fractional noise, obtained as the increments of self-similar processes, and fractional autoregressive moving average processes. The process $\{y_t\}$ is said to exhibit long memory (or long-range dependence) when $0 < d_0 < 0.5$, short memory (or short-range dependence) when $d_0 = 0$, and antipersistence when $-0.5 < d_0 < 0$, and in this paper we undertake an extensive examination of the consequences for estimation of such processes of mis-specifying the short run dynamics. In so doing we provide a significant extension of earlier work on this particular form of mis-specification in Yajima (1992) and Chen and Deo (2006), as well as complementing work that focuses on other types of mis-specification in fractional settings, as in Hassler (1994) and Crato and Taylor (1996). Our work also complements that of Robinson (2014), where mis-specification of the local-to zero characterisation of long memory is examined, and that in Cavaliere, Nielsen and Taylor (2017), where a comprehensive treatment of inference in fractional models under very general forms of heteroscedasticity is provided. Whilst mis-specification *per se* is not the focus of the latter paper, the proof of convergence to a pseudo-true parameter of the conditional sum of squares (CSS) estimator under the imposition of incorrect linear restrictions bears some relationship with our more general results on mis-specified estimators in the fractional setting. Our results also generalise the existing literature on the properties of various parametric estimators - including their asymptotic equivalence - in *correctly specified* long memory models; see Fox and Taqqu (1986), Dahlhaus (1989), Giraitis and Surgailis (1990), Sowell (1992), Beran (1995), Robinson (2006) and Hualde and Robinson (2011), among others.

---

[†]This research has been supported by Australian Research Council (ARC) Discovery Grants DP150101728 and DP170100729.

[‡]Corresponding author: Gael Martin, Department of Econometrics and Business Statistics, Monash University, Clayton, Victoria 3800, Australia. Tel.: +61-3-9905-1189; fax: +61-3-9905-5474; email: gael.martin@monash.edu.

We begin by showing that four alternative parametric techniques – frequency domain maximum likelihood (FML), Whittle, time domain maximum likelihood (TML) and CSS – converge to a common pseudo-true parameter value when the short memory component is mis-specified.[1] Convergence is established for all three forms of dependence in the true data generating process (TDGP) - long memory, short memory and antipersistence. We establish convergence by demonstrating that when the mis-specified model is evaluated at points in the parameter space where the fractional index $d$ exceeds $d_0 - 0.5$ the FML criterion function has a deterministic limit, but that the FML criterion function is *divergent* otherwise. The difference in the behaviour of the FML criterion function on subsets of the parameter space implies that the objective function does not behave uniformly. (See Cavaliere, Nielsen and Taylor, 2017; Hualde and Robinson, 2011; Robinson, 1995, for related discussion.) This lack of uniformity makes proofs of convergence across the whole parameter space more complex than usual, but solutions presented in the previously cited references can be tailored to the current situation. We then show that under common mis-specification the criterion functions that define all three alternative estimators behave in a manner similar to that of the FML criterion. All four estimators are, accordingly, shown to converge to the same pseudo-true parameter value – by definition the common value that optimizes all four limiting objective functions.

Secondly, we derive closed-form representations for the first-order conditions that define the pseudo-true parameters for completely general autoregressive fractionally integrated moving average (ARFIMA) model structures – both true and mis-specified. This represents a substantial extension of the analysis in Chen and Deo (2006), in which the FML estimator under mis-specification was first investigated, but with expressions for the relevant first-order conditions provided for certain special specifications only, and with convergence established solely for *long* memory Gaussian processes.

Thirdly, we extend the asymptotic theory established by Chen and Deo (2006) for the FML estimator in the long memory Gaussian process case to the other three estimators, under long memory, short memory and antipersistence for both the TDGP and the estimated model, and without the imposition of Gaussianity. We show that all four methods are asymptotically equivalent in that they converge in distribution under common mis-specification. The convergence rate and nature of the asymptotic distribution is determined by the deviation of the pseudo-true value of the fractional index, $d_1$ say, from the true value, $d_0$, with three critical ranges for $d^* = d_0 - d_1 < 0.5$ given by $d^* > 0.25$, $d^* = 0.25$ and $d^* < 0.25$. This nonstandard distributional behaviour for all four parametric estimators introduces a further degree of complexity into the analysis, and contrasts sharply with earlier results established in correctly specified models, where separate modeling of the short-run and fractional dynamics results in the asymptotic distribution of the parameter estimates being normal and free of the fractional indices.

---

[1] Given that each of these estimators can be derived from a Gaussian likelihood, but we do not presuppose Gaussianity, each could be designated as a 'quasi' maximum likelihood estimator in the usual way; however for the sake of notational simplicity we avoid this qualifying term.

A fourth contribution of our paper is the provision of a closed-form representation of the (common) asymptotic distribution that obtains under the most extreme type of mis-specification – whereby both a $\sqrt{n}$ rate of convergence and limiting Gaussianity is lost – together with a demonstration of how to implement the distribution numerically using appropriate truncation of the series expansion that characterizes the distribution. This then enables us to illustrate graphically the differences in the rates at which the finite sample distributions of the four different estimators approach the (common) asymptotic distribution. Notably, for $d^* \geq 0.25$, there is a distinct grouping into frequency domain and time domain techniques; with the latter tending to replicate the asymptotic distribution more closely than the former in small samples.

Finally, we perform an extensive set of simulation experiments in which the relative finite sample performance of all four mis-specified estimators is assessed. The experiments are first conducted assuming a known (zero) mean, in line with the theoretical derivations in the paper, and then re-run with the mean estimated. The ranking of the estimators, in terms of bias and MSE, is shown to depend heavily on whether the mean is specified or estimated, a conclusion that parallels results documented previously for correctly specified ARFIMA models (see, for example, Cheung and Diebold, 1994; Nielsen and Frederiksen, 2005; Sowell, 1992).

We note that in defining the Whittle estimator we focus on a particular approximation to the frequency domain Gaussian (negative) log likelihood, in which *sums* over Fourier frequencies are used to approximate the relevant integrals. Despite the analytical equivalence of this estimator with the FML estimator for large $n$, the small sample performances of the two procedures will be seen to differ systematically. In common with the FML approach, this form of Whittle estimator is invariant to the mean of the process. For interest, we also present selected numerical results pertaining to the integral-based form of the Whittle estimator (referred to hereafter as 'exact Whittle'), both for the known mean case and when the mean is unknown, the lack of invariance of this estimator to the mean rendering this latter exercise of particular interest.

The paper is organized as follows. In Section 2 we define the estimation problem, namely producing an estimate of the parameters of a fractionally integrated model when the component of the model that characterizes the short term dynamics is mis-specified. The criterion functions that define the Whittle, TML and CSS estimators, as well as the FML estimator, are specified, and we demonstrate that all four estimators possess a common probability limit under mis-specification. The limiting form of the criterion function for a mis-specified ARFIMA model is presented in Section 3, under complete generality for the short memory dynamics in the true process and estimated model, and closed-form expressions for the first-order conditions that define the pseudo-true values of the parameters are then given. The asymptotic equivalence of all four estimation methods is proved in Section 4. The finite sample performance of the alternative estimators of $d$ in the mis-specified model – with reference to estimating the pseudo-true value of $d$ – is documented in Section 5. The form of the sampling distribution is recorded, as is the bias and mean squared error (MSE), under different degrees of mis-specification, for all four estimators. Bias and MSE results are also documented for the exact Whittle estimator. The experiments are first conducted assuming a known (zero) mean, in line with the theoretical derivations in the paper. In this case, the CSS estimator exhibits superior performance, in terms of bias and mean squared error, across a range of mis-specification settings, whilst the performance of the FML estimator is notably inferior. We then re-run the simulations using demeaned data. Only the time domain estimators, plus the exact Whittle estimator, are affected by this change, and with the rate of convergence of the sample mean being slow under long memory we find that the superiority of the time domain estimators is diminished - the (sums-based) Whittle estimator now being the best performer overall. The paper concludes in Section 6 with a brief summary and some discussion of several issues that arise from the work. The proofs of the results presented in the paper are assembled in Appendix A. Appendix B contains certain technical derivations referenced in the text.

## 2 Estimation Under Mis-specification of the Short Run Dynamics

Assume that $\{y_t\}$ is generated from a TDGP that is a purely-nondeterministic stationary and ergodic process with spectral density given by

$$\frac{\sigma_0^2}{2\pi} f_0(\lambda) = \frac{\sigma_0^2}{2\pi} g_0(\lambda) \left(2\sin(\lambda/2)\right)^{-2d_0}, \tag{1}$$

where $\sigma_0^2$ is the innovation variance, $g_0(\lambda)$ is a real valued symmetric function of $\lambda$ defined on $[-\pi, \pi]$ that is bounded above and bounded away from zero, and $-0.5 < d_0 < 0.5$. Then there exists a zero mean process $\{\varepsilon_t\}$ of uncorrelated random variables with variance $\sigma_0^2$ such that $\{y_t\}$ has the moving average representation

$$y_t = \mu_0 + \sum_{j=0}^{\infty} b_{0j} \varepsilon_{t-j}, \quad t \in \mathbb{Z} = 0, \pm 1, \ldots, \tag{2}$$

where $\{b_{0j}\}$ is a sequence of constants satisfying $b_{00} = 1$ and $\sum_{j=0}^{\infty} b_{0j}^2 < \infty$, and $f_0(\lambda) = |b_0(\exp(\imath\lambda))|^2$, $\lambda \in [-\pi, \pi]$, with $(1-z)^{d_0} b_0(z) = c_0(z) = \sum_{j=0}^{\infty} c_{0j} z^j$ and $0 < |c_0(z)|$, $|z| \leq 1$. We will suppose that $c(\exp(\imath\lambda))$ is differentiable in $\lambda$ for all $\lambda \neq 0$ with a derivative that is of order $O(|\lambda|^{-1})$ as $\lambda \to 0$, and that

(A.1) For all $t \in \mathbb{Z}$ we have $E_0[\varepsilon_t|\mathbb{F}_{t-1}] = 0$ and $E_0[\varepsilon_t^2|\mathbb{F}_{t-1}] = \sigma_0^2$, a.s. where $\mathbb{F}_{t-1}$ in the conditional expectations is the sigma-field of events generated by $\varepsilon_s$, $s \leq t-1$. Here, and in what follows, the zero subscript denotes that the moments are defined with respect to the TDGP.

The conditions imposed on $c_0(z)$ imply that $g_0(\lambda)$ corresponds to the spectrum of an invertible short-memory process that is bounded and bounded away from zero for all $\lambda \in [-\pi, \pi]$ and the TDGP satisfies *Conditions A* of Hannan (1973, page 131). Assumption A.1 was introduced into time series analysis by Hannan (1973) and has been employed by several authors in investigations of both short memory and fractional linear processes since. The assumption that $\{\varepsilon_t\}$ is a conditionally homoscedastic martingale difference process circumvents the need to assume independence or identical distributions for the innovations, but rules out heteroscedasticity (see Cavaliere, Nielsen and Taylor, 2017, pages 5-6).

The model to be estimated is a parametric specification for the spectral density of $\{y_t\}$ of the form

$$\frac{\sigma^2}{2\pi} f_1(\boldsymbol{\eta}, \lambda) = \frac{\sigma^2}{2\pi} g_1(\boldsymbol{\beta}, \lambda) \left(2\sin(\lambda/2)\right)^{-2d}, \tag{3}$$

where $g_1(\boldsymbol{\beta}, \lambda)$ is a real valued symmetric function of $\lambda$ defined on $[-\pi, \pi]$. The parameter of interest will be taken as $\boldsymbol{\eta} = (d, \boldsymbol{\beta}^\top)^\top$, where $d \in (-0.5, 0.5)$ and $\boldsymbol{\beta} \in \mathbb{B}$, where $\mathbb{B}$ is an $l$-dimensional compact convex set in $\mathbb{R}^l$. The variance $\sigma^2$ will be viewed as a supplementary or nuisance parameter. The model is to be estimated from a realization $y_t$, $t = 1, \ldots, n$, of $\{y_t\}$ and, in order that the structure of the model should parallel the assumed properties of the TDGP, it will be assumed that the model is specified in such a way that:

(A.2) For all $\boldsymbol{\beta} \in \mathbb{B}$, $\int_{-\pi}^{\pi} \log g_1(\boldsymbol{\beta}, \lambda) d\lambda = 0$, and $\boldsymbol{\beta} \neq \boldsymbol{\beta}'$ implies that $g_1(\boldsymbol{\beta}, \lambda) \neq g_1(\boldsymbol{\beta}', \lambda)$ on a set of positive Lebesgue measure.

(A.3) The function $g_1(\boldsymbol{\beta}, \lambda)$ is differentiable with respect to $\lambda$, with derivative $\partial g_1(\boldsymbol{\beta}, \lambda)/\partial \lambda$ continuous at all $(\boldsymbol{\beta}, \lambda)$, $\lambda \neq 0$, and $|\partial g_1(\boldsymbol{\beta}, \lambda)/\partial \lambda| = O(|\lambda|^{-1})$ as $\lambda \to 0$. Furthermore, $\inf_{\boldsymbol{\beta}} \inf_\lambda g_1(\boldsymbol{\beta}, \lambda) > 0$ and $\sup_{\boldsymbol{\beta}} \sup_\lambda g_1(\boldsymbol{\beta}, \lambda) < \infty$.

If there exists a subset of $[-\pi, \pi]$ with non-zero Lebesgue measure in which $g_1(\boldsymbol{\beta}, \lambda) \neq g_0(\lambda)$ for all $\boldsymbol{\beta} \in \mathbb{B}$ then the model will be referred to as a mis-specified model (MisM).

The above TDGP and modelling assumptions encompass the standard parametric models, such as fractional noise, and fractional exponential and ARFIMA processes. (A detailed outline of the properties of such processes is provided in Beran, 1994) We will return to a discussion of these regularity conditions later, where a strengthening of these conditions – detailed below – will be required in order to derive our asymptotic distribution theory. Meanwhile we note (for future reference) that an ARFIMA model for a time series $\{y_t\}$ may be defined as follows,

$$\phi(L)(1-L)^d\{y_t - \mu\} = \theta(L)\varepsilon_t, \qquad (4)$$

where $\mu = E(y_t)$, $L$ is the lag operator such that $L^k y_t = y_{t-k}$, and $\phi(z) = 1 + \phi_1 z + ... + \phi_p z^p$ and $\theta(z) = 1 + \theta_1 z + ... + \theta_q z^q$ are the autoregressive and moving average operators respectively, where it is assumed that $\phi(z)$ and $\theta(z)$ have no common roots and that the roots lie outside the unit circle. The errors $\{\varepsilon_t\}$ are assumed to be a white noise sequence with finite variance $\sigma_\varepsilon^2 > 0$. For $|d| < 0.5$, $\{y_t\}$ can be represented as an infinite-order moving average of $\{\varepsilon_t\}$ with square-summable coefficients and, hence, on the assumption that the specification in (4) is correct, $\{y_t\}$ is defined as the limit in mean square of a covariance-stationary process. When $0 < d < 0.5$ neither the moving average coefficients nor the autocovariances of the process are absolutely summable, declining at a hyperbolic rate rather than the exponential rate typical of an ARMA process, with the term 'long memory' invoked accordingly. Thus, for an ARFIMA model we have $g_1(\boldsymbol{\beta}, \lambda) = |\theta(e^{i\lambda})|^2/|\phi(e^{i\lambda})|^2$ where $\boldsymbol{\beta} = (\phi_1, \phi_2, ..., \phi_p, \theta_1, \theta_2, ..., \theta_q)^T$ and an $ARFIMA(p, d, q)$ model will be mis-specified if the realizations are generated from a true $ARFIMA(p_0, d_0, q_0)$ process and any of $\{p \neq p_0 \cup q \neq q_0\} \setminus \{p_0 \leq p \cap q_0 \leq q\}$ obtain.

We consider estimators of the parameter of interest, $\boldsymbol{\eta} = (d, \boldsymbol{\beta}^T)^T$, that are obtained by minimizing a criterion function $Q_n(\boldsymbol{\eta})$ over a user-assigned compact subset of the parameter space $(-0.5, 0.5) \times \mathbb{B}$,

$$\mathbb{E}_\delta = \mathbb{D}_\delta \times \mathbb{B} \quad \text{where} \quad \mathbb{D}_\delta = \{d : |d| \leq 0.5 - \delta\}, \text{ for some } 0 < \delta \ll 0.5. \qquad (5)$$

The bound on $|d|$ must be set by the practitioner via some criterion that reflects numerical precision. Under mis-specification the generic estimator, denoted by $\widehat{\boldsymbol{\eta}}_1$ for the time being, is obtained by minimizing $Q_n(\boldsymbol{\eta})$ assuming that $\{y_t\}$ follows the MisM.[2] In Section 2.1 we specify the form of $Q_n(\boldsymbol{\eta})$ associated with the FML estimator considered in Chen and Deo (2006) and outline its relationship with the criterion functions underlying two alternative versions of the frequency domain estimator introduced by Whittle, making it clear which form of Whittle estimator is the focus of our theoretical investigations. In Section 2.2 we define the two time domain estimators that we consider here, TML and CSS, and their associated criterion functions.

---

[2] We follow the usual convention by denoting the estimator obtained under mis-specification as $\widehat{\boldsymbol{\eta}}_1$ rather than simply by $\widehat{\boldsymbol{\eta}}$, say. This is to make it explicit that the estimator is obtained under mis-specification and does not correspond to the estimator produced under the correct specification of the model, which could be denoted by $\widehat{\boldsymbol{\eta}}_0$.

Anticipating the convergence results that follow later in this section, for any given $Q_n(\boldsymbol{\eta})$ a law of large numbers can be combined with standard arguments to establish that on compact subsets of $\{\mathbb{D}_\delta \cap \{d : (d_0 - d) < 0.5\}\} \times \mathbb{B}$, i.e. subsets of $\mathbb{E}_\delta$ where $\mathbb{D}_\delta$ intersects with $\{d : (d_0 - d) < 0.5\}$, the criterion $Q_n(\boldsymbol{\eta})$ will converge uniformly to the non-stochastic limiting objective function

$$Q(\boldsymbol{\eta}) = \lim_{n \to \infty} E_0 \left[ Q_n(\boldsymbol{\eta}) \right] = \frac{\sigma_0^2}{2\pi} \int_0^\pi \frac{f_0(\lambda)}{f_1(\boldsymbol{\eta},\lambda)} d\lambda \,. \tag{6}$$

If, on the other hand, $Q_n(\boldsymbol{\eta})$ is evaluated on a subset of $\mathbb{E}_\delta$ where $\mathbb{D}_\delta$ intersects with $\{d : (d_0 - d) \geq 0.5\}$, then the criterion function is divergent. The latter corresponds to the integral on the right hand side in (6) being assigned the value $\infty$ if $(d_0 - d) \geq 0.5$ (see the comment by Hannan on his Lemma 2 in Hannan, 1973, page 134). This difference in behaviour of the criterion function about the point $d_0 - d = 0.5$ implies that $Q_n(\boldsymbol{\eta})$ does not converge uniformly on subsets of the parameter space that include this point. Nevertheless, as will be demonstrated below, provided that $\boldsymbol{\eta}_1 \in \mathbb{E}_\delta$, where $\boldsymbol{\eta}_1$ is the minimizer of $Q(\boldsymbol{\eta})$, $Q_n(\widehat{\boldsymbol{\eta}}_1)$ will converge to $Q(\boldsymbol{\eta}_1)$ and $\widehat{\boldsymbol{\eta}}_1$ will converge to $\boldsymbol{\eta}_1$ as a consequence.

In Section 2.3 we derive our asymptotic results pertaining to the convergence of $Q_n(\boldsymbol{\eta})$ and demonstrate the relationships between the limiting criterion functions of the Whittle, TML and CSS estimators to the limiting criterion function of the FML estimator. The value that minimizes the limiting criterion function of all four estimators is shown to be identical, and the asymptotic convergence of all four estimators to the common pseudo-true parameter, $\boldsymbol{\eta}_1$, is thereby established. In the theoretical derivations we adopt the assumption of a known mean for both the true and estimated models, with a zero value specified without loss of generality.

## 2.1 Frequency domain estimators

In their paper Chen and Deo (2006) focus on the estimator of $\boldsymbol{\eta} = (d, \boldsymbol{\beta}^T)^T$ defined as the value of $\boldsymbol{\eta}$ that minimizes the objective function

$$Q_n^{(1)}(\boldsymbol{\eta}) = \frac{2\pi}{n} \sum_{j=1}^{\lfloor n/2 \rfloor} \frac{I(\lambda_j)}{f_1(\boldsymbol{\eta},\lambda_j)} \,, \tag{7}$$

where $I(\lambda_j)$ is the periodogram, defined as $I(\lambda) = \frac{1}{2\pi n} |\sum_{t=1}^n y_t \exp(-i\lambda t)|^2$ evaluated at the Fourier frequencies $\lambda_j = 2\pi j/n$; $(j = 1, ..., \lfloor n/2 \rfloor)$, $\lfloor x \rfloor$ is the largest integer not greater than $x$. We have labeled this the FML estimator. The objective function in (7) is an approximation to the frequency domain Gaussian (negative) log-likelihood introduced initially by Whittle (1953) for short-range dependent processes, namely

$$W_n(\sigma^2, \boldsymbol{\eta}) = \int_{-\pi}^{\pi} \left\{ \log \frac{\sigma^2}{2\pi} f_1(\boldsymbol{\eta}, \lambda) + \frac{2\pi I(\lambda)}{\sigma^2 f_1(\boldsymbol{\eta}, \lambda)} \right\} d\lambda, \tag{8}$$

and it coincides with the frequency domain objective function considered in Hannan (1973). Concentrating out $\sigma^2$ in (8) and minimizing the associated profile function with respect to $\boldsymbol{\eta}$ produces what we refer to as the exact Whittle estimator.

An alternative approximation to the Whittle criterion function in (8), considered for example in Beran (1994), is

$$Q_n^{(2)}(\sigma^2, \boldsymbol{\eta}) = \frac{2\pi}{n} \sum_{j=1}^{\lfloor n/2 \rfloor} \log \left[ \frac{\sigma^2}{2\pi} f_1(\boldsymbol{\eta},\lambda_j) \right] + \frac{(2\pi)^2}{\sigma^2 n} \sum_{j=1}^{\lfloor n/2 \rfloor} \frac{I(\lambda_j)}{f_1(\boldsymbol{\eta},\lambda_j)}. \tag{9}$$

Taking $\boldsymbol{\eta}$ as the parameter of interest and concentrating $Q_n^{(2)}(\sigma^2, \boldsymbol{\eta})$ with respect to $\sigma^2$ indicates that the value of $\sigma^2$ that minimises (9) is given by $\widehat{\sigma}^2(\boldsymbol{\eta}) = 2Q_n^{(1)}(\boldsymbol{\eta})$. Substituting back in to (9) yields the (negative) profile likelihood,

$$Q_n^{(2)}(\boldsymbol{\eta}) = \frac{2\pi}{2} \log\left(\frac{\widehat{\sigma}^2(\boldsymbol{\eta})}{2\pi}\right) + \frac{2\pi}{n} \sum_{j=1}^{\lfloor n/2 \rfloor} \log f_1(\boldsymbol{\eta}, \lambda_j) + \pi.$$

Minimization of $Q_n^{(2)}(\boldsymbol{\eta})$ with respect to $\boldsymbol{\eta}$ yields what we call (simply) the Whittle estimator, and which is the form of Whittle procedure that features in our theoretical derivations. Since $\lim_{n\to\infty} \frac{2\pi}{n} \sum_{j=1}^{\lfloor n/2 \rfloor} \log f_1(\boldsymbol{\eta}, \lambda_j) = 0$ (see Appendix A.3.1) it follows that this estimator is equivalent to the FML estimator for large $n$. However, as indicated in Boes, Davis and Gupta (1989), and as will be seen in the simulation results documented in Section 5, the finite sample performance of these two estimators differs. For interest (and as prompted by a referee) we also report selected numerical results on the finite sample performance of the exact Whittle estimator described above.

## 2.2 Time domain estimators

The criterion functions of the two alternative time domain estimators are defined as follows:

- Let $\mathbf{Y}^T = (y_1, y_2, ..., y_n)$ and denote the variance covariance matrix of $\mathbf{Y}$ derived from the mis-specified model by $\sigma^2 \boldsymbol{\Sigma}_\eta = [\gamma_1(i-j)]$, $i, j = 1, 2, ..., n$, where

$$\gamma_1(\tau) = \gamma_1(-\tau) = \frac{\sigma^2}{2\pi} \int_{-\pi}^{\pi} f_1(\boldsymbol{\eta}, \lambda) e^{i\lambda\tau} d\lambda.$$

The Gaussian log-likelihood function for the TML estimator is

$$-\frac{1}{2}\left(n\log(2\pi\sigma^2) + \log|\boldsymbol{\Sigma}_\eta| + \frac{1}{\sigma^2}(\mathbf{Y} - \mu\mathbf{l})^T \boldsymbol{\Sigma}_\eta^{-1} (\mathbf{Y} - \mu\mathbf{l})\right), \tag{10}$$

where $\mathbf{l}^T = (1, 1, ..., 1)$, and maximizing (10) is equivalent to minimizing the criterion function

$$Q_n^{(3)}(\sigma^2, \boldsymbol{\eta}) = \log\sigma^2 + \frac{1}{n}\log|\boldsymbol{\Sigma}_\eta| + \frac{1}{n\sigma^2}(\mathbf{Y} - \mu\mathbf{l})^T \boldsymbol{\Sigma}_\eta^{-1} (\mathbf{Y} - \mu\mathbf{l}). \tag{11}$$

- To construct the CSS estimator note that we can expand $(1-z)^d$ in a binomial expansion as

$$(1-z)^d = \sum_{j=0}^{\infty} \frac{\Gamma(j-d)}{\Gamma(j+1)\Gamma(-d)} z^j, \tag{12}$$

where $\Gamma(\cdot)$ is the gamma function. Furthermore, since $g_1(\boldsymbol{\beta}, \lambda)$ is bounded, by Assumption (A.3), we can employ the method of Whittle (Whittle, 1984, §2.8) to construct an autoregressive operator $\alpha(\boldsymbol{\beta}, z) = \sum_{i=0}^{\infty} \alpha_i(\boldsymbol{\beta}) z^i$ such that $g_1(\boldsymbol{\beta}, \lambda) = |\alpha(\boldsymbol{\beta}, e^{i\lambda})|^{-2}$. The objective function of the CSS estimation method then becomes

$$Q_n^{(4)}(\boldsymbol{\eta}) = \frac{1}{n} \sum_{t=1}^{n} e_t^2, \tag{13}$$

where

$$e_t = \sum_{i=0}^{t-1} \tau_i(\boldsymbol{\eta})(y_{t-i} - \mu) \tag{14}$$

and the coefficients $\tau_j(\boldsymbol{\eta})$, $j = 0, 1, 2, \ldots$, are given by $\tau_0(\boldsymbol{\eta}) = 1$ and

$$\tau_j(\boldsymbol{\eta}) = \sum_{s=0}^{j} \frac{\alpha_{j-s}(\boldsymbol{\beta})\Gamma(j-d)}{\Gamma(j+1)\Gamma(-d)}, \quad j = 1, 2, \ldots. \quad (15)$$

As with the FML estimator, the CSS estimate of $\sigma^2$ is given implicitly by the minimum value of the criterion function.

We can think of the CSS estimator as providing an approximation to the TML estimator that parallels the approximation of the FML and (sums-based) Whittle estimators to the exact Whittle estimator.

## 2.3 Convergence Properties

In Chen and Deo (2006) it is shown that if $\{y_t\}$ is a long-range dependent Gaussian process, then on subsets of the parameter space of the form $(\delta, 0.5 - \delta) \times \Phi$, where $0 < \delta < 0.25$ and $\Phi$ is a compact convex set, we have (for $Q_n^{(1)}(\boldsymbol{\eta})$ defined in (7)) $\text{plim}_{n \to \infty} |Q_n^{(1)}(\boldsymbol{\eta}) - Q(\boldsymbol{\eta})| = 0$ (Chen and Deo, 2006, Lemma 2). The minimum of the limiting objective function $Q(\boldsymbol{\eta})$ then defines a pseudo-true parameter value to which the FML estimator will converge, since with the addition of the assumption that there exists a unique vector $\boldsymbol{\eta}_1 = (d_1, \boldsymbol{\beta}_1^T)^T \in (\delta, 0.5 - \delta) \times \Phi$ that minimizes $Q(\boldsymbol{\eta})$, it follows that the FML estimator will converge to $\boldsymbol{\eta}_1$.

Because Chen and Deo assumed that the TDGP was a long memory process and that in the MisM the fractional index was similarly confined to the long memory region, they did not explicitly consider the case where $(d_0 - d) \geq 0.5$. In contrast, as noted with reference to the TDGP in (1), our work allows for $0 \leq |d_0| < 0.5$, and involves the specification of the appropriate user-assigned compact subset for $\boldsymbol{\eta} = (d, \boldsymbol{\beta}^T)^T$ in (5). This implies a wider range of values for $(d_0 - d)$ and, hence, the need for our analysis to deal with the differing behaviour of $Q_n^{(1)}(\boldsymbol{\eta})$ about the point $d_0 - d = 0.5$ alluded to above. To achieve this, we divide the parameter space $\mathbb{E}_\delta$ into three disjoint sub-sets:

1. $\mathbb{E}_\delta^0 = \mathbb{D}_\delta^0 \times \mathbb{B}$ where $\mathbb{D}_\delta^0 = \mathbb{D}_\delta \cap \{d : -(1-2\delta) \leq (d_0 - d) \leq 0.5 - \delta\}$,

2. $\overline{\mathbb{E}}_{\delta 1}^0 = \overline{\mathbb{D}}_{\delta 1}^0 \times \mathbb{B}$ where $\overline{\mathbb{D}}_{\delta 1}^0 = \mathbb{D}_\delta \cap \{d : 0.5 - \delta < (d_0 - d) < 0.5\}$ and

3. $\overline{\mathbb{E}}_{\delta 2}^0 = \overline{\mathbb{D}}_{\delta 2}^0 \times \mathbb{B}$ where $\overline{\mathbb{D}}_{\delta 2}^0 = \mathbb{D}_\delta \cap \{d : 0.5 \leq (d_0 - d) \leq 1 - 2\delta\}.$

The superscript '0' is used to indicate that the relevant subspaces relate to the deviation $(d_0 - d)$ assuming that $d_0 \in \mathbb{D}_\delta$. The notation in 2. and 3. is used to denote the breakdown of the complement of the set in 1, $\mathbb{E}_\delta^0$, into two disjoint subsets, $\overline{\mathbb{E}}_{\delta 1}^0$ and $\overline{\mathbb{E}}_{\delta 2}^0$. This division of the parameter space of $(d_0 - d)$ is depicted graphically in Figure 1.

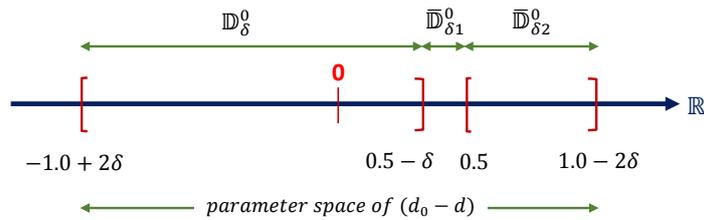

Figure 1: Graphical illustration of the division of the parameter space of $(d_0 - d)$

We will establish that on the subset $\mathbb{E}_\delta^0$ we have $\lim_{n \to \infty} Q_n^{(1)}(\boldsymbol{\eta}) = Q(\boldsymbol{\eta})$ almost surely and uniformly in $\boldsymbol{\eta}$, where $Q(\boldsymbol{\eta})$ is defined as in (6), whereas $Q_n^{(1)}(\boldsymbol{\eta})$ is of order $O(\delta^{-1})$ on $\overline{\mathbb{E}}_{\delta 1}^0$ and is divergent as $n \to \infty$ on $\overline{\mathbb{E}}_{\delta 2}^0$. This is the content of Lemmas 1, 2, 3 and 4 below.

Proposition 1 then establishes that the FML estimator converges to $\boldsymbol{\eta}_1 = \arg\min_{\boldsymbol{\eta}} Q(\boldsymbol{\eta})$. We will also establish that the convergence/divergence properties of each of the three alternative estimators, Whittle, TML and CSS, is the same as that of the FML estimator. The upshot of this is summarised in Theorem 1.

**Lemma 1** *Suppose that the TDGP of $\{y_t\}$ is as prescribed in equations (1) and (2) and that the MisM is specified as in (3). Assume also that Assumptions $A.1 - A.3$ are satisfied. Then for any constant $\nu_f > 0$*

$$\left| \frac{2\pi}{n} \sum_{j=1}^{\lfloor n/2 \rfloor} \frac{I(\lambda_j)}{f_1(\boldsymbol{\eta},\lambda_j) + \nu_f} - \frac{\sigma_0^2}{2\pi} \int_0^\pi \frac{f_0(\lambda)}{f_1(\boldsymbol{\eta},\lambda) + \nu_f} d\lambda \right|$$

*converges to zero almost surely and uniformly in $\boldsymbol{\eta}$ on $\mathbb{E}_\delta^0$.*

Since, obviously, $f_1(\boldsymbol{\eta},\lambda) < f_1(\boldsymbol{\eta},\lambda) + \nu_f$ it follows from Lemma 1 that,

$$\liminf_{n\to\infty} Q_n^{(1)}(\boldsymbol{\eta}) \geq \lim_{n\to\infty} \frac{2\pi}{n} \sum_{j=1}^{\lfloor n/2 \rfloor} \frac{I(\lambda_j)}{f_1(\boldsymbol{\eta},\lambda_j) + \nu_f}$$

$$= \frac{\sigma_0^2}{2\pi} \int_0^\pi \frac{f_0(\lambda)}{f_1(\boldsymbol{\eta},\lambda) + \nu_f} \quad a.s.$$

uniformly in $\boldsymbol{\eta}$ on $\mathbb{E}_\delta^0$. Letting $\delta_f \to 0$ and applying Lebegue's monotone convergence theorem gives

$$\liminf_{n\to\infty} Q_n^{(1)}(\boldsymbol{\eta}) \geq Q(\boldsymbol{\eta}) = \frac{\sigma_0^2}{2\pi} \int_0^\pi \frac{f_0(\lambda)}{f_1(\boldsymbol{\eta},\lambda)} d\lambda \quad a.s.$$

To establish that $Q(\boldsymbol{\eta})$ also provides a limit superior for $Q_n^{(1)}(\boldsymbol{\eta})$ when $\boldsymbol{\eta} \in \mathbb{E}_\delta^0$ we will use the following lemma.

**Lemma 2** *Suppose that the conditions of Lemma 1 hold. Set*

$$h_1(\boldsymbol{\eta},\lambda) = \begin{cases} f_1(\boldsymbol{\eta},\lambda), & f_1(\boldsymbol{\eta},\lambda) \geq \nu_f \\ \nu_f, & f_1(\boldsymbol{\eta},\lambda) < \nu_f, \end{cases}$$

*where $\nu_f > 0$. Then for all $\nu_f > 0$,*

$$\left| \frac{2\pi}{n} \sum_{j=1}^{\lfloor n/2 \rfloor} \frac{I(\lambda_j)}{h_1(\boldsymbol{\eta},\lambda_j)} - \frac{\sigma_0^2}{2\pi} \int_0^\pi \frac{f_0(\lambda)}{h_1(\boldsymbol{\eta},\lambda)} d\lambda \right|$$

*converges to zero almost surely uniformly in $\boldsymbol{\eta}$ on $\mathbb{E}_\delta^0$.*

The following lemma shows that the limiting form of the FML criterion function presented by Chen and Deo (2006), for Gaussian processes (specifically) and only in the case where both $d$ and $d_0$ lie in the interval $(0, 0.5)$, holds more generally, and can incorporate all three forms of memory - long memory, short memory and antipersistence - in both the true and estimated models.

**Lemma 3** *Suppose that the conditions of Lemmas 1 and 2 hold. Then*

$$\lim_{n\to\infty} \sup_{\boldsymbol{\eta}\in\mathbb{E}_\delta^0} |Q_n^{(1)}(\boldsymbol{\eta}) - Q(\boldsymbol{\eta})| = 0\,.$$

Lemma 4 then indicates that for points in $\mathbb{E}_\delta$ where $(d_0 - d) > 0.5 - \delta$, $0 < \delta < 0.5$, uniform convergence of the criterion function $Q_n^{(1)}(\boldsymbol{\eta})$ fails.

**Lemma 4** *Suppose that the TDGP of $\{y_t\}$ is as prescribed in equations (1) and (2) and that the MisM is specified as in (3). Assume also that Assumptions $A.1 - A.3$ are satisfied. Then for all $\boldsymbol{\eta} \in \overline{\mathbb{E}}_{\delta 1}^0$ we have $\liminf_{n \to \infty} Q_n^{(1)}(\boldsymbol{\eta}) = O(\delta^{-1})$ and for $\boldsymbol{\eta} \in \overline{\mathbb{E}}_{\delta 2}^0$*

$$\liminf_{n \to \infty} Q_n^{(1)}(\boldsymbol{\eta}) \geq C > 0$$

*almost surely for all $C$, no matter how large.*

Note that Lemma 4 implies that as $n$ increases, and for all $\delta$ sufficiently small, $\widehat{\boldsymbol{\eta}}_1^{(1)} = \arg\min_{\boldsymbol{\eta}} Q_n^{(1)}(\boldsymbol{\eta})$ cannot lie in $\overline{\mathbb{E}}_{\delta 1}^0 \cup \overline{\mathbb{E}}_{\delta 2}^0$. Proposition 1 now follows as an almost immediate corollary of the previous developments if we suppose that the following additional assumption holds:

(A.4) There exists a unique pseudo-true parameter vector $\boldsymbol{\eta}_1 = (d_1, \boldsymbol{\beta}_1^T)^T$ belonging to the subset $\mathbb{E}_\delta^0$ that satisfies $\boldsymbol{\eta}_1 = \arg\min_{\boldsymbol{\eta}} Q(\boldsymbol{\eta})$.

Proposition 1 establishes the convergence of the FML estimator to $\boldsymbol{\eta}_1$ under the same generality for both the TDGP and MisM as highlighted above (*cf.* Chen and Deo, 2006, Corollary 1).

**Proposition 1** *Suppose that the TDGP of $\{y_t\}$ is as prescribed in equations (1) and (2) and that the MisM is specified as in (3). Assume also that Assumptions $A.1 - A.4$ are satisfied. Let $\widehat{\boldsymbol{\eta}}_1^{(1)}$ denote the FML estimator obtained by minimising the criterion function $Q_n^{(1)}(\boldsymbol{\eta})$ over $\mathbb{E}_\delta$. Then $\lim_{n \to \infty} Q_n^{(1)}(\widehat{\boldsymbol{\eta}}_1^{(1)}) = Q(\boldsymbol{\eta}_1)$ and $\widehat{\boldsymbol{\eta}}_1^{(1)} \to \boldsymbol{\eta}_1$ almost surely.*

Index now by $i = 2, 3$ and $4$ the estimators associated with the Whittle, TML and CSS criterion functions respectively; that is $\widehat{\boldsymbol{\eta}}_1^{(i)}$ minimises $Q_n^{(i)}(\cdot)$, $i = 2, 3, 4$, with each viewed as a function of $\boldsymbol{\eta}$. Given the relationships between $Q_n^{(1)}(\cdot)$ and $Q_n^{(i)}(\cdot)$, $i = 2, 3, 4$, as outlined in the appendix, it follows that $\widehat{\boldsymbol{\eta}}_1^{(i)}$, $i = 1, 2, 3, 4$, must share the same convergence properties. Thus we can state the following theorem:

**Theorem 1** *Suppose that the TDGP of $\{y_t\}$ is as prescribed in equations (1) and (2) and that the MisM is specified as in (3). Assume also that Assumptions $A.1 - A.4$ are satisfied. Let $\widehat{\boldsymbol{\eta}}_1^{(i)}$, $i = 1, 2, 3, 4$, denote, respectively, the FML, Whittle, TML and CSS estimators of the parameter vector $\boldsymbol{\eta} = (d, \boldsymbol{\beta}^T)^T$ of the MisM. Then $\lim_{n \to \infty} \|\widehat{\boldsymbol{\eta}}_1^{(i)} - \widehat{\boldsymbol{\eta}}_1^{(j)}\| = 0$ almost surely for all $i, j = 1, 2, 3, 4$, where the common limiting value of $\widehat{\boldsymbol{\eta}}_1^{(i)}$, $i = 1, 2, 3, 4$, is $\boldsymbol{\eta}_1 = \arg\min_{\boldsymbol{\eta}} Q(\boldsymbol{\eta})$.*

Having established that the four parametric estimators converge towards a common $\boldsymbol{\eta}_1$, we can as a consequence now broaden the applicability of the asymptotic distributional results derived by Chen and Deo (2006) for the FML estimator. This we do in Section 4 by establishing that all four alternative parametric estimators converge in distribution for all three forms of memory - long memory, short memory and antipersistence.

Prior to doing this, however, we note that Cavaliere, Nielsen and Taylor (2017) have shown that if $g_0(\lambda)$ has a parametric form that is known, but the parameter values that characterise it are not, then the CSS estimator will convergence to a pseudo-true value if it is evaluated whilst imposing incorrect linear parameter constraints (Cavaliere, Nielsen and Taylor, 2017, Theorem 5(ii)). Theorem 1 provides a generalisation of this result by extending it to the FML, Whittle and TML estimators, and by also allowing for the possibility that the parametric form of the model, $g_1(\lambda)$, may itself be mis-specified. In what follows we indicate the precise form of the limiting objective function $Q(\boldsymbol{\eta})$, and the associated first-order conditions that define the pseudo-true value $\boldsymbol{\eta}_1$ of the four estimation procedures, in the ARFIMA case.

## 3 Pseudo-True Parameters Under ARFIMA Mis-Specification

Under Assumptions $A.1 - A.4$, the value of $\boldsymbol{\eta}_1 = \arg\min_{\boldsymbol{\eta}} Q(\boldsymbol{\eta})$ can be determined as the solution of the first-order condition $\partial Q(\boldsymbol{\eta})/\partial \boldsymbol{\eta} = 0$, and Chen and Deo (2006) illustrate the relationship between $\partial \log Q(\boldsymbol{\eta})/\partial d$ and the deviation $d^* = d_0 - d_1$ for the simple special case in which the TDGP is an $ARFIMA(0, d_0, 1)$ and the MisM is an $ARFIMA(0, d, 0)$. They then cite (without providing detailed derivations) certain results that obtain when the MisM is an $ARFIMA(1, d, 0)$. Here we provide a significant generalization, by deriving expressions for both $Q(\boldsymbol{\eta})$ and the first-order conditions that define the pseudo-true parameters, under the full $ARFIMA(p_0, d_0, q_0)/ARFIMA(p, d, q)$ dichotomy for the true process and the estimated model. Representations of the associated expressions via polynomial and power series expansions suitable for the analytical investigation of $Q(\boldsymbol{\eta})$ are presented. It is normally not possible to solve the first order conditions $\partial Q(\boldsymbol{\eta})/\partial \boldsymbol{\eta} = 0$ exactly as they are both nonlinear and (in general) defined as infinite sums. Instead one would determine the estimate numerically, via a Newton iteration for example, with the series expansions replaced by finite sums. An evaluation of the magnitude of the approximation error produced by any power series truncation that might arise from such a numerical implementation is given. The results are then illustrated in the special case where $p_0 = q = 0$, in which case true MA short memory dynamics of an arbitrary order are mis-specified as AR dynamics of an arbitrary order. In this particular case, as will be seen, no truncation error arises in the computations.

To begin, denote the spectral density of the TDGP, a general $ARFIMA(p_0, d_0, q_0)$ process, by

$$\frac{\sigma_0^2}{2\pi} f_0(\lambda) = \frac{\sigma_0^2}{2\pi} \frac{\left|1 + \theta_{10} e^{i\lambda} + ... + \theta_{q_0 0} e^{iq_0 \lambda}\right|^2}{\left|1 + \phi_{10} e^{i\lambda} + ... + \phi_{p_0 0} e^{ip_0 \lambda}\right|^2} |2\sin(\lambda/2)|^{-2d_0},$$

and that of the MisM, an $ARFIMA(p, d, q)$ model, by

$$\frac{\sigma^2}{2\pi} f_1(\boldsymbol{\eta}, \lambda) = \frac{\sigma^2}{2\pi} \frac{\left|1 + \theta_1 e^{i\lambda} + ... + \theta_q e^{iq\lambda}\right|^2}{\left|1 + \phi_1 e^{i\lambda} + ... + \phi_p e^{ip\lambda}\right|^2} |2\sin(\lambda/2)|^{-2d}.$$

Substituting these expressions into the limiting objective function we obtain the representation

$$Q(\boldsymbol{\eta}) = \frac{\sigma_0^2}{2\pi} \int_0^\pi \frac{f_0(\lambda)}{f_1(\boldsymbol{\eta}, \lambda)} d\lambda = \frac{\sigma_0^2}{2\pi} \int_0^\pi \frac{|A_\beta(e^{i\lambda})|^2}{|B_\beta(e^{i\lambda})|^2} |2\sin(\lambda/2)|^{-2(d_0 - d)} d\lambda, \qquad (16)$$

where

$$A_\beta(z) = \sum_{j=0}^{\underline{q}} a_j z^j = \theta_0(z)\phi(z) = (1 + \theta_{10} z + ... + \theta_{q_0 0} z^{q_0})(1 + \phi_1 z + ... + \phi_p z^p), \qquad (17)$$

with $\underline{q} = q_0 + p$ and

$$B_\beta(z) = \sum_{j=0}^{\underline{p}} b_j z^j = \phi_0(z)\theta(z) = (1 + \phi_{10} z + ... + \phi_{p_0 0} z^{p_0})(1 + \theta_1 z + ... + \theta_q z^q), \qquad (18)$$

with $\underline{p} = p_0 + q$. The expression for $Q(\boldsymbol{\eta})$ in (16) takes the form of the variance of an ARFIMA process with MA operator $A_\beta(z)$, AR operator $B_\beta(z)$ and fractional index $d_0 - d$. It follows that $Q(\boldsymbol{\eta})$ could be evaluated using the procedures presented in Sowell (1992). Sowell's algorithms are based upon series expansions in gamma and hypergeometric functions however, and although they are suitable for numerical calculations, they do not readily lend themselves to the analytical investigation of $Q(\boldsymbol{\psi})$. We therefore seek an alternative formulation.

Let $C(z) = \sum_{j=0}^{\infty} c_j z^j = A_\beta(z)/B_\beta(z)$ where $A_\beta(z)$ and $B_\beta(z)$ are as defined in (17) and (18) respectively. Then (16) can be expanded to give

$$Q(\boldsymbol{\eta}) = 2^{1-2(d_0-d)}\frac{\sigma_0^2}{2\pi}\left[\sum_{j=0}^{\infty}\sum_{k=0}^{\infty} c_j c_k \int_0^{\pi/2} \cos\left(2(j-k)\lambda\right)\sin(\lambda)^{-2(d_0-d)}d\lambda\right].$$

Using standard results for the integral $\int_0^{\pi}(\sin x)^{v-1}\cos(ax)dx$ from Gradshteyn and Ryzhik (2007, p 397) yields, after some algebraic manipulation,

$$Q(\boldsymbol{\eta}) = \frac{\sigma_0^2}{2(1-2(d_0-d))}\left[\sum_{j=0}^{\infty}\sum_{k=0}^{\infty}\frac{c_j c_k \cos((j-k)\pi)}{\mathcal{B}(1-(d_0-d)+(j-k),1-(d_0-d)-(j-k))}\right],$$

where $\mathcal{B}(a,b)$ denotes the Beta function. This expression can in turn be simplified to

$$Q(\boldsymbol{\eta}) = \frac{\sigma_0^2 \Gamma(1-2(d_0-d))}{2\Gamma^2(1-(d_0-d))}\}K(\boldsymbol{\eta}), \tag{19}$$

where

$$K(\boldsymbol{\eta}) = \sum_{j=0}^{\infty} c_j^2 + 2\sum_{k=0}^{\infty}\sum_{j=k+1}^{\infty} c_j c_k \rho(j-k)$$

and

$$\rho(h) = \prod_{i=1}^{h}\left(\frac{(d_0-d)+i-1}{i-(d_0-d)}\right), \quad h=1,2,\ldots.$$

Using (19) we now derive the form of the first-order conditions that define $\boldsymbol{\eta}_1$, namely $\partial Q(\boldsymbol{\eta})/\partial\boldsymbol{\eta} = 0$. Differentiating $Q(\boldsymbol{\eta})$ first with respect to $\beta_r$, $r=1,\ldots,l$, and then $d$ gives:

$$\frac{\partial Q(\boldsymbol{\eta})}{\partial \beta_r} = \{\frac{\sigma_0^2 \Gamma(1-2(d_0-d))}{2\Gamma^2(1-(d_0-d))}\}\frac{\partial K(\boldsymbol{\eta})}{\partial \beta_r}, \quad r=1,2,...,l,$$

where

$$\frac{\partial K(\boldsymbol{\eta})}{\partial \beta_r} = \sum_{j=1}^{\infty} 2c_j \frac{\partial c_j}{\partial \beta_r} + 2\sum_{k=0}^{\infty}\sum_{j=k+1}^{\infty}(c_k\frac{\partial c_j}{\partial \beta_r} + \frac{\partial c_k}{\partial \beta_r}c_j)\rho(j-k),$$

and

$$\frac{\partial Q(\boldsymbol{\eta})}{\partial d} = \{\frac{\sigma_0^2 \Gamma(1-2(d_0-d))}{2\Gamma^2(1-(d_0-d))}\}\left\{2\left(\Psi[1-2(d_0-d)]-\Psi[1-(d_0-d)]\right)K(\boldsymbol{\eta}) + \frac{\partial K(\boldsymbol{\eta})}{\partial d}\right\},$$

where $\Psi(\cdot)$ denotes the digamma function and

$$\frac{\partial K(\boldsymbol{\eta})}{\partial d} = 2\sum_{k=0}^{\infty}\sum_{j=k+1}^{\infty} c_j c_k \rho(j-k)\left\{2\Psi[1-(d_0-d)]\right.$$
$$\left. -\Psi[1-(d_0-d)+(j-k)]-\Psi[1-(d_0-d)-(j-k)]\right\}.$$

Eliminating the common (non-zero) factor $\{\pi\frac{\sigma_{\varepsilon 0}^2}{\sigma^2}\frac{\Gamma(1-2(d_0-d))}{\Gamma^2(1-(d_0-d))}\}$ from both $\partial Q(\boldsymbol{\eta})/\partial\boldsymbol{\beta}$ and $\partial Q(\boldsymbol{\eta})/\partial d$, it follows that the pseudo-true parameter values of the $ARFIMA(p,d,q)$ MisM can be obtained by solving

$$\frac{\partial K(\boldsymbol{\eta})}{\partial \beta_r} = 0, \quad r=1,2,...,l, \tag{20}$$

and
$$2(\Psi[1-2(d_0-d)] - \Psi[1-(d_0-d)])K(\boldsymbol{\eta}) + \frac{\partial K(\boldsymbol{\eta})}{\partial d} = 0 \tag{21}$$

for $\beta_{r1}$, $r = 1, \ldots, l$, and $d_1$ using appropriate algebraic and numerical procedures. A corollary of the following theorem is that $\boldsymbol{\eta}_1$ can be calculated to any desired degree of numerical accuracy by truncating the series expansions in the expressions for $K(\boldsymbol{\eta})$, $\partial K(\boldsymbol{\eta})/\partial \boldsymbol{\beta}$ and $\partial K(\boldsymbol{\eta})/\partial d$ after a suitable number of $N$ terms before substituting into (20) and (21) and solving (numerically) for $\phi_{i1}$, $i = 1, 2, ..., p$, $\theta_{j1}$, $j = 1, 2, ..., q$, and $d_1$.

**Theorem 2** *Set $C_N(z) = \sum_{j=0}^{N} c_j z^j$ and let $Q_N(\boldsymbol{\eta}) = (\sigma_0^2/\sigma_\varepsilon^2) I_N$ where the integral $I_N = \int_0^\pi |C_N(\exp(-i\lambda))|^2 |2\sin(\lambda/2)|^{-2(d_0-d)} d\lambda$. Then*

$$Q(\boldsymbol{\eta}) = Q_N(\boldsymbol{\eta}) + R_N = \left\{ \frac{\sigma_0^2 \Gamma(1-2(d_0-d))}{2\Gamma^2(1-(d_0-d))} \right\} K_N(\boldsymbol{\eta}) + R_N$$

*where*

$$K_N(\boldsymbol{\eta}) = \sum_{j=0}^{N} c_j^2 + 2 \sum_{k=0}^{N-1} \sum_{j=k+1}^{N} c_j c_k \rho(j-k)$$

*and there exists a $\zeta$, $0 < \zeta < 1$, such that $R_N = O(\zeta^{(N+1)}) = o(N^{-1})$. Furthermore, $\partial Q_N(\boldsymbol{\eta})/\partial \boldsymbol{\eta} = \partial Q(\boldsymbol{\eta})/\partial \boldsymbol{\eta} + o(N^{-1})$.*

By way of illustration, consider the case of mis-specifying a true $ARFIMA(0, d_0, q_0)$ process by an $ARFIMA(p, d, 0)$ model. When $p_0 = q = 0$ we have $B_\beta(z) \equiv 1$ and $C(z)$ is polynomial, $C(z) = 1 + \sum_{j=1}^{q} c_j z^j$ where $c_j = \sum_{r=\max\{0,j-p\}}^{\min\{j,p\}} \theta_{(j-r)0} \phi_r$. Abbreviating the latter to $\sum_r \theta_{(j-r)0} \phi_r$, this then gives us:

$$K(d, \phi_1, \ldots, \phi_p) = \sum_{j=0}^{q} (\sum_r \theta_{(j-r)0} \phi_r)^2 +$$
$$2 \sum_{k=0}^{q-1} \sum_{j=k+1}^{q} (\sum_r \theta_{(j-r)0} \phi_r)(\sum_r \theta_{(k-r)0} \phi_r) \rho(j-k) \, ;$$

and setting $\theta_{s0} \equiv 0$, $s \ni [0, 1, \ldots, q_0]$,

$$\frac{\partial K(d, \phi_1, \ldots, \phi_p)}{\partial \phi_r} = \sum_{j=1}^{q} 2(\sum_r \theta_{(j-r)0} \phi_r) \theta_{(j-r)0} +$$
$$2 \sum_{k=0}^{q-1} \sum_{j=k+1}^{q} \left\{ (\sum_r \theta_{(j-r)0} \phi_r) \theta_{(k-r)0} + \theta_{(j-r)0} (\sum_r \theta_{(k-r)0} \phi_r) \right\} \rho(j-k) \, ,$$

$r = 1, \ldots, p$, and

$$\frac{\partial K(d, \phi_1, \ldots, \phi_p)}{\partial d} = 2 \sum_{k=0}^{q-1} \sum_{j=k+1}^{q} (\sum_r \theta_{(j-r)0} \phi_r)(\sum_r \theta_{(k-r)0} \phi_r) \rho(j-k) \times$$
$$(2\Psi[1-(d_0-d)] - \Psi[1-(d_0-d)+(j-k)] - \Psi[1-(d_0-d)-(j-k)])$$

for the required derivatives. The pseudo-true values $\phi_{r1}$, $r = 1, \ldots, p$, and $d_1$ can now be obtained by solving (20) and (21) having inserted these exact expressions for $K(d, \phi_1, \ldots, \phi_p)$, $\partial K(d, \phi_1, \ldots, \phi_p)/\partial \phi_r$, $r = 1, \ldots, p$, and $\partial K(d, \phi_1, \ldots, \phi_p)/\partial d$ into the equations.

Let us further highlight some features of this special case by focussing on the example where the TDGP is an $ARFIMA(0, d_0, 1)$ and the MisM an $ARFIMA(1, d, 0)$. In this example $q = 2$ and $C(z) = 1 + c_1 z + c_2 z^2$ where, neglecting the first order MA and AR coefficient subscripts, $c_1 = (\theta_0 + \phi)$ and $c_2 = \theta_0 \phi$. The second factor of the criterion function in (19) is now

$$K(d, \phi) = 1 + (\theta_0 + \phi)^2 + (\theta_0 \phi)^2 \\ + \frac{2\left[\theta_0\phi(d_0 - d + 1) - (1 + \theta_0\phi)(\theta_0 + \phi)(d_0 - d - 2)\right](d_0 - d)}{(d_0 - d - 1)(d_0 - d - 2)}. \quad (22)$$

The derivatives $\partial K(d, \phi)/\partial \phi$ and $\partial K(d, \phi)/\partial d$ can be readily determined from (22) and hence the pseudo-true values $d_1$ and $\phi_1$ evaluated.

It is clear from (22) that for given values of $|\theta_0| < 1$ we can treat $K(d, \phi)$ as a function of $\widetilde{d} = (d_0 - d)$ and $\phi$, and hence treat $Q(d, \phi) = Q(\boldsymbol{\eta})$ similarly. Figure 2 depicts the contours of $Q(d, \phi)$ graphed as a function of $\widetilde{d}$ and $\phi$ for the values of $\theta_0 = \{-0.7, -0.637014, -0.3\}$ when $\sigma_\varepsilon^2 = \sigma_0^2$. Pre-empting the discussion to come in the following section, the values of $\theta_0$ are deliberately chosen to coincide with $d^* = d_0 - d_1$ being respectively greater than, equal to and less than 0.25. The three graphs in Figure 2 clearly demonstrate the divergence in the asymptotic criterion function that occurs as $\widetilde{d} = (d_0 - d)$ approaches 0.5 and they illustrate that although the location of $(d_1, \phi_1)$ may be unambiguous, the sensitivity of $Q(d, \phi)$ to perturbations in $(d, \phi)$ can be very different depending on the value of $d^* = d_0 - d_1$.[3] In Figure 2a the contours indicate that when $d^* > 0.25$ the limiting criterion function has hyperbolic profiles in a small neighbourhood of the pseudo-true parameter point $(d_1, \phi_1)$, with similar but more locally quadratic behaviour exhibited in Figure 2b when $d^* = 0.25$. The contours of $Q(d, \phi)$ in Figure 2c, corresponding to $d^* < 0.25$, are more elliptical and suggest that in this case the limiting criterion function is far closer to being globally quadratic around $(d_1, \phi_1)$. It turns out that these three different forms of $Q(d, \phi)$, reflecting the most, intermediate, and the least mis-specified cases, correspond to the three different forms of asymptotic distribution presented in the following section.

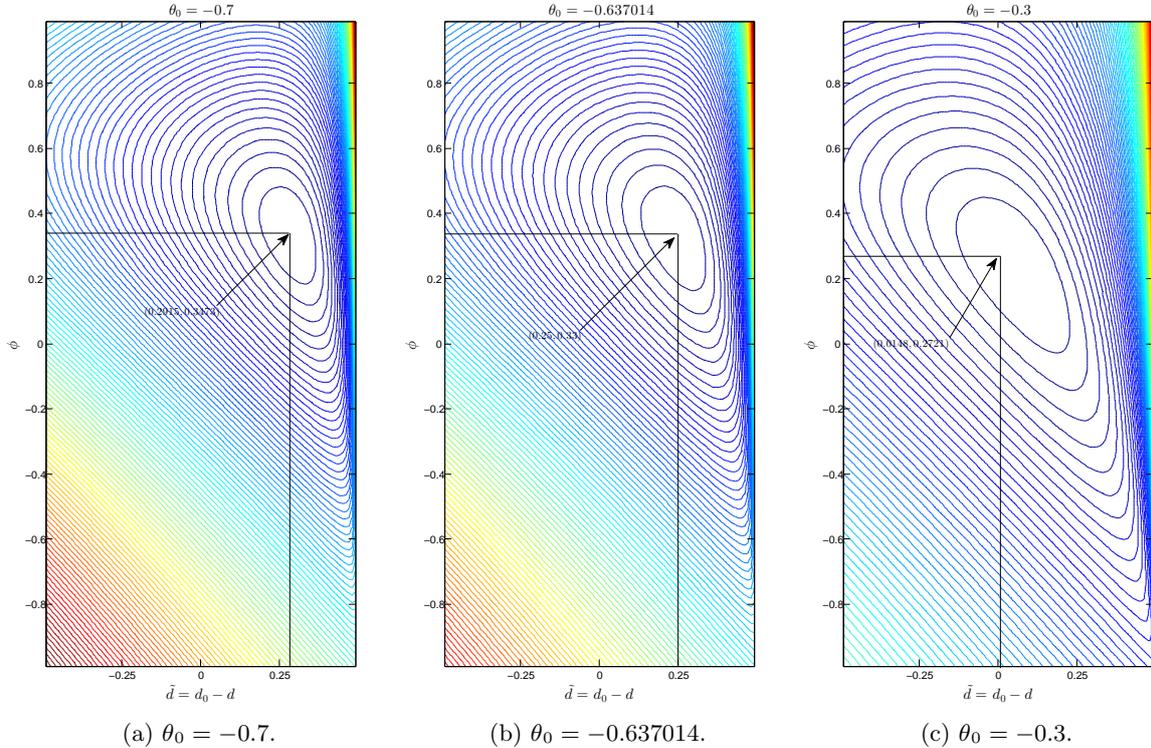

(a) $\theta_0 = -0.7$.  (b) $\theta_0 = -0.637014$.  (c) $\theta_0 = -0.3$.

Figure 2: Contour plot of $Q(d, \phi)$ against $\widetilde{d} = d_0 - d$ and $\phi$ for the mis-specification of an ARFIMA$(0, d_0, 1)$ TDGP by an ARFIMA$(1, d, 0)$ MM; $\widetilde{d} \in (-0.5, 0.5)$, $\phi \in (-1, 1)$. Pseudo-true coordinates $(d_0 - d_1, \phi_1)$ are (a) $(0.2915, 0.3473)$, (b) $(0.25, 0.33)$ and (c) $(0.0148, 0.2721)$.

---

[3] All the numerical results presented in this paper have been produced using MATLAB *2011b*, version *7.13.0.564 (R2011b)*.

## 4 Asymptotic Distributions

In this section we show that the asymptotic distribution of the FML estimator derived in Chen and Deo (2006) in the context of long range dependence is also applicable to the Whittle, TML and CSS estimators, and that all four estimators are, hence, asymptotically equivalent under mis-specification. As was highlighted by Chen and Deo, the rate of convergence and the nature of the asymptotic distribution of the FML estimator is determined by the deviation of the pseudo-true value $d_1$ from the true value $d_0$.[4] Theorem 3 shows that in the event that any one of the FML, Whittle, TML or CSS estimators possesses one of the asymptotic distributions as described in the theorem, then all four estimators will share the same asymptotic distribution, and this will hold for all three forms of memory in the TDGP and the mis-specified model. We comment further on this matter below.

For each of the estimators the asymptotic distributions are obtained via the usual Taylor series expansion of the score function, having first established convergence, and consequently stronger smoothness conditions are required to establish the asymptotic distribution theory and to ensure that the asymptotic variance-covariance matrix of the estimators is well defined. We will therefore suppose:

(A.5) The function $g_1(\boldsymbol{\beta}, \lambda)$ of the MisM is thrice differentiable with continuous third derivatives. Furthermore, the derivatives satisfy;

A.5.1 $\sup_\lambda \sup_{\boldsymbol{\beta}} \left| \frac{\partial g_1(\boldsymbol{\beta}, \lambda)}{\partial \beta_i} \right| < \infty$, $1 \leqslant i \leqslant l$,

A.5.2 $\sup_\lambda \sup_{\boldsymbol{\beta}} \left| \frac{\partial^2 g_1(\boldsymbol{\beta}, \lambda)}{\partial \beta_i \partial \beta_j} \right| < \infty$, $\sup_\lambda \sup_{\boldsymbol{\beta}} \left| \frac{\partial^2 g_1(\boldsymbol{\beta}, \lambda)}{\partial \beta_i \partial \lambda} \right| < \infty$, $1 \leqslant i, j \leqslant l$, and

A.5.3 $\sup_\lambda \sup_{\boldsymbol{\beta}} \left| \frac{\partial^3 g_1(\boldsymbol{\beta}, \lambda)}{\partial \beta_i \partial \beta_j \partial \beta_k} \right| < \infty$, $1 \leqslant i, j, k \leqslant l$.

Assumptions A.2-A.5 are similar to the assumptions adopted by Fox and Taqqu (1986) and Dahlhaus (1989) in the context of correct specification, and they are in essence equivalent to the conditions used in the work of Chen and Deo (2006) on the mis-specified case. In order to derive the asymptotic distribution we will assume that $\{\varepsilon_t\}$ is a strictly stationary, regular process that satisfies the following weak dependence and moment conditions.

(A.1′) The innovation $\{\varepsilon_t\}$ satisfies Assumption (A.1). Moreover, $E_0[|\varepsilon_t|^{4+p}] < \infty$ for some $p \in (0, \infty)$ and there exist finite constants $\mu_3$ and $\mu_4$ such that $E_0[\varepsilon_t^3 | \mathbb{F}_{t-1}] = \mu_3$ and $E_0[\varepsilon_t^4 | \mathbb{F}_{t-1}] = \mu_4$ a.s. for all $t \in \mathbb{Z}$.

Assumption (A.1′) states that $\{\varepsilon_t\}$ is a martingale difference sequence, a not unreasonable assumption that has almost become standard in the asymptotic analysis of time series, and it is typical to assume finite bounds on the first four moments of the innovation process (see Cavaliere, Nielsen and Taylor, 2017, for a detailed explanation of the importance of such bounds). Assumption (A.1′) implies that $\{\varepsilon_t\}$ is completely regular, and that $\{\varepsilon_t^p\}$ is a uniformly integrable sequence for any $p \leq 4$.[5]

---

[4] As already noted, the results in Chen and Deo presupposed that the parameter space of the estimated model coincided with the long memory region assumed for the TDGP. Since $d_1$ is only defined for $(d_0 - d_1) < 0.5$ it follows that the distributional results they presented for the FML estimator were only valid for this region, something that was not explicitly mentioned in their original derivation.

[5] This assumption is closely related to Assumption (A.1) of Lahiri (2003), which specifies a set of weak dependence and moment conditions on $\{\varepsilon_t\}$ based on $\alpha$−mixing.

**Theorem 3** *Suppose that the TDGP of $\{y_t\}$ is as prescribed in equations (1) and (2), and that the MisM is specified as in (3), and assume that Assumptions A.1' and A.2 − A.5 hold. Let*

$$\mathbf{B} = -\frac{\sigma_0^2}{\pi}\int_{-\pi}^{\pi}\frac{f_0(\lambda)}{f_1^3(\boldsymbol{\eta}_1,\lambda)}\frac{\partial f_1(\boldsymbol{\eta}_1,\lambda)}{\partial \boldsymbol{\eta}}\frac{\partial f_1(\boldsymbol{\eta}_1,\lambda)}{\partial \boldsymbol{\eta}^T}d\lambda + \frac{\sigma_0^2}{2\pi}\int_{-\pi}^{\pi}\frac{f_0(\lambda)}{f_1^2(\boldsymbol{\eta}_1,\lambda)}\frac{\partial^2 f_1(\boldsymbol{\eta}_1,\lambda)}{\partial \boldsymbol{\eta}\partial \boldsymbol{\eta}^T}d\lambda, \quad (23)$$

*and set $\boldsymbol{\mu}_n = \mathbf{B}^{-1}E_0\left(\frac{\partial Q_n(\boldsymbol{\eta}_1)}{\partial \boldsymbol{\eta}}\right)$ where $Q_n(\cdot)$ denotes the objective function that defines $\widehat{\boldsymbol{\eta}}_1$.*[6] *Let $\widehat{\boldsymbol{\eta}}_1$ denote the estimator obtained by minimising $Q_n(\boldsymbol{\eta})$ over the compact set $\mathbb{E}_\delta$ where $\boldsymbol{\eta}_1 \in \mathbb{E}_\delta$ and assume that $\boldsymbol{\eta}_1 \ni \partial \mathbb{E}_\delta$ where $\partial \mathbb{E}_\delta$ is the boundary of the set $\mathbb{E}_\delta$. Then the FML, Whittle, TML or CSS estimators are asymptotically equivalent with a common limiting distribution as delineated in Cases 1, 2 and 3:*

*Case 1:* When $d^* = d_0 - d_1 > 0.25$,

$$\frac{n^{1-2d^*}}{\log n}(\widehat{\boldsymbol{\eta}}_1 - \boldsymbol{\eta}_1 - \boldsymbol{\mu}_n) \xrightarrow{D} \mathbf{B}^{-1}\left(\sum_{j=1}^{\infty}W_j, 0, ...0\right)^T, \quad (24)$$

*where $\sum_{j=1}^{\infty}W_j$ is defined as the mean-square limit of the random sequence $\sum_{j=1}^{s}W_j$ as $s \to \infty$, wherein*

$$W_j = \frac{(2\pi)^{1-2d^*}g_0(0)}{j^{2d^*}g_1(\boldsymbol{\beta},0)}\left[U_j^2 + V_j^2 - E_0(U_j^2 + V_j^2)\right],$$

*and $\{U_j\}$ and $\{V_k\}$ denote sequences of Gaussian random variables with zero mean and covariances $Cov_0(U_j, U_k) = Cov_0(U_j, V_k) = Cov_0(V_j, V_k)$ with*

$$Cov_0(U_j, V_k) = \iint_{[0,1]^2}\{\sin(2\pi jx)\cos(2\pi ky) + \sin(2\pi kx)\cos(2\pi jy)\}|x-y|^{2d_0-1}dxdy.$$

*Case 2:* When $d^* = d_0 - d_1 = 0.25$,

$$n^{1/2}\left[\overline{\Lambda}_{dd}\right]^{-1/2}(\widehat{\boldsymbol{\eta}}_1 - \boldsymbol{\eta}_1) \xrightarrow{D} \mathbf{B}^{-1}(Z, 0, ..., 0)^T, \quad (25)$$

*where*

$$\overline{\Lambda}_{dd} = \frac{1}{n}\sum_{j=1}^{n/2}\left(\frac{\sigma_0^2 f_0(\lambda_j)}{2\pi f_1(\boldsymbol{\eta}_1,\lambda_j)}\frac{\partial \log f_1(\boldsymbol{\eta}_1,\lambda_j)}{\partial d}\right)^2,$$

*and $Z$ is a standard normal random variable.*

*Case 3:* When $d^* = d_0 - d_1 < 0.25$,

$$\sqrt{n}(\widehat{\boldsymbol{\eta}}_1 - \boldsymbol{\eta}_1) \xrightarrow{D} N(0, \boldsymbol{\Xi}), \quad (26)$$

*where $\boldsymbol{\Xi} = \mathbf{B}^{-1}\boldsymbol{\Lambda}\mathbf{B}^{-1}$,*

$$\boldsymbol{\Lambda} = \frac{\sigma_0^4}{2\pi}\int_0^{\pi}\left(\frac{f_0(\lambda)}{f_1(\boldsymbol{\eta}_1,\lambda)}\right)^2\left(\frac{\partial \log f_1(\boldsymbol{\eta}_1,\lambda)}{\partial \boldsymbol{\eta}}\right)\left(\frac{\partial \log f_1(\boldsymbol{\eta}_1,\lambda)}{\partial \boldsymbol{\eta}}\right)^T d\lambda.$$

---

[6] Heuristically, $\boldsymbol{\mu}_n$ measures the bias associated with the estimator $\widehat{\boldsymbol{\eta}}_1$. That is, $\boldsymbol{\mu}_n \approx E_0(\widehat{\boldsymbol{\eta}}_1) - \boldsymbol{\eta}_1$. Note that the expression for $\boldsymbol{\mu}_n$ given in Chen and Deo (2006, p 263) contains a typographical error; the proofs in that paper use the correct expression. The derivation of $\boldsymbol{\mu}_n$ for all four estimation methods considered in this paper is provided in Appendix B.

A key point to note from the three cases delineated in Theorem 3 is that when the deviation between the true and pseudo-true values of $d$ is sufficiently large ($d^* \geq 0.25$) – something that is related directly to the degree of mis-specification of $g_0(\lambda)$ by $g_1(\boldsymbol{\beta}, \lambda)$ – the $\sqrt{n}$ rate of convergence is lost, with the rate being arbitrarily close to zero depending on the value of $d^*$. For $d^*$ strictly greater than 0.25, asymptotic Gaussianity is also lost, with the limiting distribution being a function of an infinite sum of non-Gaussian variables. For the $d^* \geq 0.25$ case, the limiting distribution – whether Gaussian or otherwise – is degenerate in the sense that the limiting distribution for each element of $\widehat{\boldsymbol{\eta}}_1$ is a different multiple of the same random variable ($\sum_{j=1}^{\infty} W_j$ in the case of $d^* > 0.25$ and $Z$ in the case of $d^* = 0.25$).

For the form of limiting distribution that obtains in Cases 1, 2 and 3 we refer to Chen and Deo (2006, Theorems 1, 3 and 2), wherein these distributions were produced specifically for the FML estimator in the context of long range dependence. Their proofs depend on the Fourier sine and cosine transformations of the observed series being normally distributed with a given covariance structure. In Chen and Deo (2006) the latter properties are derived by assuming that $\{x(t)\}$ is a Gaussian process. Here we achieve the same outcome by employing Assumption $A.1'$ and appealing to results of Lahiri (2003) which imply that the Fourier sine and cosine transformations are asymptotically normal and hence that lemmas of Moulines and Soulier (1999) used by Chen and Deo can be applied in a more general setting.

To prove that these same limiting distributions hold for the Whittle, TML and CSS estimators we establish that $R_n(\widehat{\boldsymbol{\eta}}_1^{(i)} - \widehat{\boldsymbol{\eta}}_1^{(1)}) \to^D 0$ for $i = 2, 3$ and 4, where $R_n$ denotes the convergence rate applicable in the three different cases outlined in the theorem. We use a first-order Taylor expansion of $\partial Q_n^{(\cdot)}(\boldsymbol{\eta}_1)/\partial \boldsymbol{\eta}$ about $\partial Q_n^{(\cdot)}(\widehat{\boldsymbol{\eta}}_1^{(\cdot)})/\partial \boldsymbol{\eta} = \boldsymbol{0}$. This gives

$$\frac{\partial Q_n^{(\cdot)}(\boldsymbol{\eta}_1)}{\partial \boldsymbol{\eta}} = \frac{\partial^2 Q_n^{(\cdot)}(\dot{\boldsymbol{\eta}}_1^{(\cdot)})}{\partial \boldsymbol{\eta} \partial \boldsymbol{\eta}'} \left(\boldsymbol{\eta}_1 - \widehat{\boldsymbol{\eta}}_1^{(\cdot)}\right)$$

and

$$R_n(\widehat{\boldsymbol{\eta}}_1^{(i)} - \widehat{\boldsymbol{\eta}}_1^{(j)}) = \left[\frac{\partial^2 Q_n^{(j)}(\dot{\boldsymbol{\eta}}_1^{(j)})}{\partial \boldsymbol{\eta} \partial \boldsymbol{\eta}'}\right]^{-1} R_n \frac{\partial Q_n^{(j)}(\boldsymbol{\eta}_1)}{\partial \boldsymbol{\eta}} - \left[\frac{\partial^2 Q_n^{(i)}(\dot{\boldsymbol{\eta}}_1^{(i)})}{\partial \boldsymbol{\eta} \partial \boldsymbol{\eta}'}\right]^{-1} R_n \frac{\partial Q_n^{(i)}(\boldsymbol{\eta}_1)}{\partial \boldsymbol{\eta}},$$

where $\|\boldsymbol{\eta}_1 - \dot{\boldsymbol{\eta}}_1^{(\cdot)}\| \leq \|\boldsymbol{\eta}_1 - \widehat{\boldsymbol{\eta}}_1^{(\cdot)}\|$. Since $\operatorname{plim} \widehat{\boldsymbol{\eta}}_1^{(\cdot)} = \boldsymbol{\eta}_1$ it is therefore sufficient to show that there exists a scalar, possibly constant, function $\mathcal{C}_n(\boldsymbol{\eta})$ such that

$$\left\|\frac{\partial^2 \{\mathcal{C}_n(\boldsymbol{\eta}_1) \cdot Q_n^{(i)}(\boldsymbol{\eta}_1) - Q_n^{(j)}(\boldsymbol{\eta}_1)\}}{\partial \boldsymbol{\eta} \partial \boldsymbol{\eta}'}\right\| = o_p(1) \tag{27}$$

and

$$\operatorname{plim}_{n \to \infty} \quad R_n \left\|\mathcal{C}_n(\boldsymbol{\eta}_1) \cdot \frac{\partial Q_n^{(i)}(\boldsymbol{\eta}_1)}{\partial \boldsymbol{\eta}} - \frac{\partial Q_n^{(j)}(\boldsymbol{\eta}_1)}{\partial \boldsymbol{\eta}}\right\| = 0. \tag{28}$$

The condition in (27) is established by showing that $\partial^2 \{Q_n^{(1)}(\boldsymbol{\eta}_1)\}/\partial \boldsymbol{\eta} \partial \boldsymbol{\eta}'$ converges in probability to $\mathbf{B}$, as defined in (23), and that for each $i = 2, 3$ and 4 the corresponding Hessian is proportional to $\partial^2 \{Q_n^{(1)}(\boldsymbol{\eta}_1)\}/\partial \boldsymbol{\eta} \partial \boldsymbol{\eta}'$ with probability approaching one. The proof of (27) is fairly conventional, whereas the proof of (28) – which implicitly invokes the Cramér-Wold device since the moments (cumulants) of the asymptotically normal gradient vector are convergence determining for the limiting distributions in Theorem 3 – is more involved because of the presence of the scaling factor $R_n$. In Appendix A.5 we present the steps necessary to prove (27) and (28) for each estimator, and for TDGPs with fractional indices in the range $-0.5 < d_0 < 0.5$.

Finally, we highlight the fact that the FML and Whittle estimators are mean invariant by virtue of being defined on the non-zero fundamental Fourier frequencies. As a consequence, all convergence results presented in both this and the previous section for the FML and Whittle estimators also hold for a process that has an arbitrary (non-zero) mean, which may be unknown, thereby broadening the applicability of the theoretical results as they pertain to these particular estimators. The same is not true, however, either for the two time domain based methods or for the exact (integral-based) Whittle estimator, as will be demonstrated in Section 5.4 below.

## 5 Finite Sample Performance of the Mis-Specified Parametric Estimators of the Pseudo-True Parameter

### 5.1 Experimental design

In this section we explore the finite sample performance of the alternative methods, as it pertains to estimation of the pseudo-true value of the long memory parameter, $d_1$, under specific types of mis-specification. We give particular focus to the four estimators: $\widehat{d}_1^{(1)}$ (FML), $\widehat{d}_1^{(2)}$ (Whittle), $\widehat{d}_1^{(3)}$ (TML) and $\widehat{d}_1^{(4)}$ (CSS), for which the preceeding theoretical results have been produced. We first document the form of the finite sample distributions for each estimator by plotting the distribution of the standardized versions of the estimators, for which the asymptotic distributions are given in Cases 1, 2 and 3 respectively in Theorem 3. As part of this exercise we develop a method for obtaining the limiting distribution for $d^* > 0.25$, as the distribution does not have a closed form in this case, as well as a method for estimating the bias-adjustment term, $\boldsymbol{\mu}_n$, which is relevant for this distribution. In the figures that follow the 'Limit' curve depicts the limiting distribution of the relevant statistic. Supplementing these graphical results, we then tabulate the bias and MSE of the four different techniques, as estimators of the pseudo-true parameter $d_1$, again under specific types of mis-specification and, hence, for different values of $d^*$. These results are supplemented by bias and MSE values for the exact Whittle estimator, which we refer to as $\widehat{d}_1^{(5)}$.

Data are simulated from a zero-mean Gaussian $ARFIMA(p_0, d_0, q_0)$ process, with the method of Sowell (1992), as modified by Doornik and Ooms (2001), used to compute the exact autocovariance function for the TDGP for any given values of $p_0$, $d_0$ and $q_0$. We have produced results for $n = 100, 200, 500$ and $1000$ and for two versions of mis-specification nested in the general case for which the analytical results are derived in Section 3. However, we report selected results (only) from the full set due to space constraints. The bias and MSE results, plus certain computations needed for the numerical specification for the limiting distribution in the $d^* > 0.25$ case, are produced from $R = 1000$ replications of samples of size $n$ from the relevant TDGP. The two forms of mis-specification considered are:

**Example 1** : *An $ARFIMA(0, d_0, 1)$ TDGP, with parameter values $d_0 = \{-0.2, 0.2, 0.4\}$ and $\theta_0 = \{-0.7, -0.444978, -0.3\}$; and an $ARFIMA(0, d, 0)$ MisM. The value $\theta_0 = -0.7$ corresponds to the case where $d^* > 0.25$ and $\widehat{d}_1^{(i)}$, $i = 1, 2, 3, 4$, have the slowest rate of convergence, $n^{1-2d^*}/\log n$, and to a non-Gaussian distribution. The value $\theta_0 = -0.444978$ corresponds to the case where $d^* = 0.25$, in which case asymptotic Gaussianity is preserved but the rate of convergence is of order $(n/\log^3 n)^{1/2}$. The value $\theta_0 = -0.3$ corresponds to the case where $d^* < 0.25$, with $\sqrt{n}$-convergence to Gaussianity obtaining.*

**Example 2** : *An $ARFIMA(0, d_0, 1)$ TDGP, with parameter values $d_0 = \{-0.2, 0.2, 0.4\}$ and $\theta_0 = \{-0.7, -0.637014, -0.3\}$; and an $ARFIMA(1, d, 0)$ MisM. In this example the value $\theta_0 = -0.7$ corresponds to the case where $d^* > 0.25$, the value $\theta_0 = -0.637014$ corresponds to the case where $d^* = 0.25$, and the value $\theta_0 = -0.3$ corresponds to the case where $d^* < 0.25$.*

In Section 5.2 we document graphically the form of the finite sampling distributions of all four estimators of $d$ under the first type of mis-specification described above, and for $d_0 = 0.2$ only. The corresponding graphs under the different values of $d_0$ (and for all three cases) are qualitatively equivalent to those reported here and, hence, are not included.[7] In Section 5.3 we then report the bias and MSE of all four estimators (in terms of estimating the pseudo-true value $d_1$) under both forms of mis-specification and for all three values of $d_0$. To supplement these results, all of which are based on the assumption that the mean is known, in Section 5.4 we reproduce corresponding bias and MSE results for the two time domain estimators, plus the exact Whittle estimator, using data in which the unknown mean is estimated along with the other parameters.

## 5.2 Finite sample distributions

In this section we consider in turn the three cases listed under Theorem 3. For notational ease and clarity we use $\widehat{d}_1$ to denote the (generic) estimator obtained under mis-specification, remembering that this estimator may be produced by any one of the four estimation methods: $\widehat{d}_1^{(1)}$ to $\widehat{d}_1^{(4)}$. Similarly, we use $Q_n(\cdot)$ to denote the criterion associated with a generic estimator. Only when contrasting the (finite sample) performances of the alternative estimators do we re-introduce the superscript notation.

### 5.2.1 Case 1: $d^* > 0.25$

The limiting distribution for $\widehat{d}_1$ in this case is

$$\frac{n^{1-2d^*}}{\log n} \left( \widehat{d}_1 - d_1 - \mu_n \right) \xrightarrow{D} b^{-1} \sum_{j=1}^{\infty} W_j \,, \tag{29}$$

where $\mu_n = b^{-1} E_0 \left( \frac{\partial Q_n(\boldsymbol{\eta}_1)}{\partial d} \right)$,

$$\begin{aligned}
b &= -2 \frac{\sigma_0^2}{2\pi} \int_{-\pi}^{\pi} \frac{f_0(\lambda)}{f_1^3(\boldsymbol{\eta}_1,\lambda)} \left( \frac{\partial f_1(\boldsymbol{\eta}_1,\lambda)}{\partial d} \right)^2 d\lambda + \frac{\sigma_0^2}{2\pi} \int_{-\pi}^{\pi} \frac{f_0(\lambda)}{f_1^2(\boldsymbol{\eta}_1,\lambda)} \frac{\partial^2 f_1(\boldsymbol{\eta}_1,\lambda)}{\partial d^2} d\lambda \\
&= -2 \int_0^{\pi} (1 + \theta_0^2 + 2\theta_0 \cos(\lambda))(2\sin(\lambda/2))^{-2d^*} (2\log(2\sin(\lambda/2)))^2 d\lambda \,, \tag{30}
\end{aligned}$$

and $W_j = \frac{(2\pi)^{1-2d^*}(1+\theta_0^2)}{j^{2d^*}} \left[ U_j^2 + V_j^2 - E_0(U_j^2 + V_j^2) \right]$, with $\{U_j\}$ and $\{V_k\}$ as defined in Theorem 3. (With reference to Theorem 3, both $\mathbf{B}$ and $\boldsymbol{\mu}_n$ in (24) are here scalars since in Example 1 there is only one parameter to estimate under the MisM, namely $d$. Hence the obvious changes made to notation. All other notation is as defined in the theorem.)

Given that the distribution in (29) is non-standard and does not have a closed form representation, consideration must be given to its numerical evaluation. In finite samples the bias-adjustment term $\mu_n$ (which approaches zero as $n \to \infty$) also needs to be calculated. We tackle each of these issues in turn, beginning with the computation of $\mu_n$.

(1) From Theorem 3 it is apparent that in general the formula for $\mathbf{B}$ is independent of the estimation method, but the calculation of $\boldsymbol{\mu}_n$ requires separate evaluation of $E_0(\partial Q_n(\boldsymbol{\eta}_1)/\partial \boldsymbol{\eta})$ for each estimator. In Appendix B we provide expressions for $E_0(\partial Q_n(\boldsymbol{\eta}_1)/\partial \boldsymbol{\eta})$ for each of the four estimation methods. These formulae are used to evaluate the scalar $\mu_n$ here. Each value is then used in the specification of the standardized estimator $\frac{n^{1-2d^*}}{\log n} \left( \widehat{d}_1 - d_1 - \mu_n \right)$ in the simulation experiments.

---

[7]These additional graphical results are available from the authors on request.

(2) Quantification of the distribution of $\sum_{j=1}^{\infty} W_j$ requires the approximation of the infinite sum of the $W_j$, plus the use of simulation to represent the (appropriately truncated) sum. We truncate the series $\sum_{j=1}^{\infty} W_j$ after $s$ terms where the truncation point $s$ is chosen such that $1 \leqslant s < \lfloor n/2 \rfloor$ with $s \to \infty$ as $n \to \infty$ (*cf.* Lemma 6 of Chen and Deo (2006)). The value of $s$ is determined using the following criterion function. Let

$$S_n = \widehat{Var}_0 \left[ \frac{n^{1-2d^*}}{\log n} \left( \widehat{d}_1 - d_1 - \mu_n \right) \right] \tag{31}$$

denote the empirical finite sample variation observed across the $R$ replications and for each $m$, $1 \leqslant m < \lfloor n/2 \rfloor$, let

$$T_m = S_n - b^{-2} \Omega_m,$$

where $\Omega_m = Var_0 \left( \sum_{j=1}^{m} W_j \right)$. Now set

$$s = \arg \min_{1 \leqslant m < \lfloor n/2 \rfloor} T_m. \tag{32}$$

Given $s$, we generate random draws of $\sum_{j=1}^{s} W_j$ via the underlying Gaussian random variables from which the $W_j$ are constructed, and produce an estimate of the limiting distribution using kernel methods.

To determine $s$ we need to evaluate

$$Var_0 \left( \sum_{j=1}^{m} W_j \right) = \sum_{j=1}^{m} Var_0 (W_j) + 2 \sum_{j=1}^{m} \sum_{\substack{k=1 \\ j \neq k}}^{m} Cov_0 (W_j, W_k). \tag{33}$$

The variance of $W_j$ in this case is

$$Var_0 \left\{ \frac{(2\pi)^{1-2d^*} (1 + \theta_0^2)}{j^{2d^*}} \left[ U_j^2 + V_j^2 - E_0 (U_j^2 + V_j^2) \right] \right\}$$

$$= \frac{(2\pi)^{2-4d^*} (1 + \theta_0^2)^2}{j^{4d^*}} \left\{ E_0 \left( U_j^2 + V_j^2 \right)^2 - \left[ E_0 \left( U_j^2 + V_j^2 \right) \right]^2 \right\}.$$

As $\{U_j\}$ and $\{V_k\}$ are normal random variables with a covariance structure as specified in Theorem 3, standard formulae for the moments of Gaussian random variables yield the result that

$$\begin{aligned} E_0 \left( U_j^2 + V_j^2 \right)^2 &= E_0 \left( U_j^4 \right) + 2 E_0 \left( U_j^2 V_j^2 \right) + E_0 \left( V_j^4 \right) \\ &= 3 \left[ Var_0 (U_j) \right]^2 + 2 \left[ Var_0 (U_j) Var_0 (V_j) + 2 Cov_0 (U_j, V_j) \right] \\ &\quad + 3 \left[ Var_0 (V_j) \right]^2 \\ &= 12 \left[ Var_0 (U_j) \right]^2 \end{aligned}$$

and

$$\begin{aligned} \left[ E_0 \left( U_j^2 + V_j^2 \right) \right]^2 &= \left[ E_0 (U_j^2) + E_0 (V_j^2) \right]^2 \\ &= \left[ Var_0 (U_j) + Var_0 (V_j) \right]^2 \\ &= 4 \left[ Var_0 (U_j) \right]^2. \end{aligned}$$

Thus,

$$Var_0(W_j) = \frac{8 (2\pi)^{2-4d^*} (1 + \theta_0^2)^2}{j^{4d^*}} \left[ Var_0 (U_j) \right]^2.$$

Similarly, the covariance between $W_j$ and $W_k$ when $j \neq k$ can be shown to be equal to

$$\frac{(2\pi)^{2-4d^*}(1+\theta_0^2)^2}{(jk)^{2d^*}} Cov_0\left(U_j^2 + V_j^2, U_k^2 + V_k^2\right)$$

$$= \frac{4(2\pi)^{2-4d^*}(1+\theta_0^2)^2}{(jk)^{2d^*}}\left[Var_0(U_j)Var_0(V_k) + 2Cov_0(U_j, V_k)\right].$$

The expression in (33) can therefore be evaluated numerically using the formula for $Cov_0(U_j, V_k)$ to calculate the necessary moments required to determine $s$ from (32).

The idea behind the use of $T_m$ is simply to minimize the difference between the second-order sample and population moments. The value of $S_n$ in (31) will vary with the estimation method of course; however, we choose $s$ based on $S_n$ calculated from the FML estimates and maintain this choice of $s$ for all other methods. The terms in (33) are also dependent on the form of both the TDGP and the MisM and hence $T_m$ needs to be determined for any specific case. The values of $s$ for the sample sizes used in the particular simulation experiment underlying the results in this section are provided in Table 1.

**Table 1** Truncation values $s$ corresponding to the case $d^* = 0.3723$ for the Ex 1: ARFIMA $(0, d_0, 1)$ TDGP with $d_0 = 0.2$ and $\theta_0 = -0.7$ vis-à-vis MisM: ARFIMA $(0, d, 0)$.

| $n$ | 100 | 200 | 500 | 1000 |
|---|---|---|---|---|
| $s$ | 36 | 75 | 162 | 230 |

Each panel in Figure 3 provides the kernel density estimate of $\frac{n^{1-2d^*}}{\log n}(\widehat{d}_1 - d_1 - \mu_n)$ under the four estimation methods, for a specific $n$ as labeled above each plot, plus the limiting distribution for the given $s$. The particular parameter values employed in the specification of the TGDP are $d_0 = 0.2$ and $\theta_0 = -0.7$, with $d^* = 0.3723$ in this case, and the values of $s$ used are those given in Table 1. From Figure 3 we see that $\frac{n^{1-2d^*}}{\log n}(\widehat{d}_1 - d_1 - \mu_n)$ is centered away from zero for all sample sizes, for all estimation methods. However, as the sample size increases the point of central location of $\frac{n^{1-2d^*}}{\log n}(\widehat{d}_1 - d_1 - \mu_n)$ approaches zero and all distributions of the standardized statistics go close to matching the asymptotic ('limit') distributions. The salient feature to be noted is the clustering that occurs, in particular for $n \leqslant 500$; that is, TML and CSS form one cluster and FML and Whittle form the other, with the time domain estimators being closer to the asymptotic distribution for all three (smaller) sample sizes.

### 5.2.2 Case 2: $d^* = 0.25$

The limiting distribution for $\widehat{d}_1$ in the case of $d^* = 0.25$ is

$$n^{1/2}[\overline{\Lambda}_{dd}]^{-1/2}\left(\widehat{d}_1 - d_1\right) \xrightarrow{D} N(0, b^{-2}), \tag{34}$$

where

$$\overline{\Lambda}_{dd} = \frac{1}{n}\sum_{j=1}^{n/2}(1 + \theta_0^2 + 2\theta_0\cos(\lambda_j))^2(2\sin(\lambda_j/2))^{-1}(2\log(2\sin(\lambda_j/2)))^2 \tag{35}$$

and $b$ is as in (30). In both (35) and (30) $\theta_0 = -0.444978$, as $d^* = 0.25$ occurs at this particular value. Once again, $d_0 = 0.2$ in the TDGP.

Each panel of Figure 4 provides the densities of $n^{1/2}[\overline{\Lambda}_{dd}]^{-1/2}\left(\widehat{d}_1 - d_1\right)$ under the four estimation methods, for a specific $n$ as labeled above each plot, plus the limiting distribution given in (34). Once again we observe a disparity between the time domain and frequency domain kernel estimates, with the pair of time domain methods yielding finite sample distributions that are closer to the limiting distribution, for all sample sizes considered. The discrepancy between the two types of methods declines as the sample size increases, with the distributions of all methods being reasonably close both to one another, and to the limiting distribution, when $n = 1000$.

Figure 3: Kernel density of $\frac{n^{1-2d^*}}{\log n}\left(\widehat{d}_1-d_1-\mu_n\right)$ for an $ARFIMA(0,d_0,1)$ TDGP with $d_0= 0.2$ and $\theta_0= -0.7$, and an $ARFIMA(0,d,0)$ MM; $d^*> 0.25$.

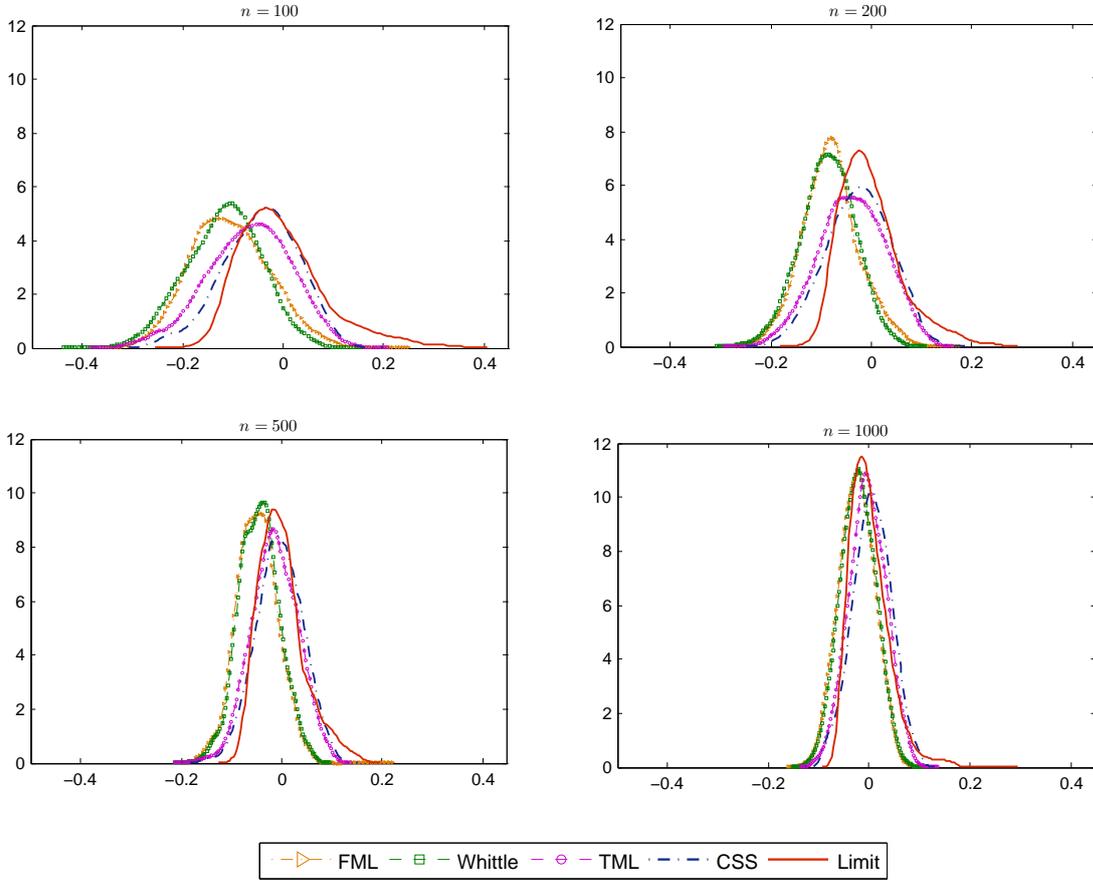

### 5.2.3 Case 3: $d^* < 0.25$

In this case we have

$$\sqrt{n}\left(\widehat{d}_1 - d_1\right) \xrightarrow{D} N(0, v^2), \tag{36}$$

where

$$v^2 = \Lambda_{11}/b^{-2}, \tag{37}$$

with

$$\begin{aligned}\Lambda_{11} &= 2\pi \frac{\sigma_0^2}{2\pi} \int_0^\pi \left(\frac{f_0(\lambda)}{f_1(d_1,\lambda)}\right)^2 \left(\frac{\partial \log f_1(d_1,\lambda)}{\partial d}\right)^2 d\lambda \\ &= 2\pi \int_0^\pi (1+\theta_0^2 + 2\theta_0 \cos(\lambda))^2 (2\sin(\lambda/2))^{-4d^*}(2\log(2\sin(\lambda/2)))^2 d\lambda,\end{aligned}$$

Figure 4: Kernel density of $n^{1/2}[\overline{\Lambda}_{dd}]^{-1/2}\left(\widehat{d}_1 - d_1\right)$ for an $ARFIMA(0, d_0, 1)$ TDGP with $d_0 = 0.2$ and $\theta_0 = -0.444978$, and an $ARFIMA(0, d, 0)$ MM, $d^* = 0.25$.

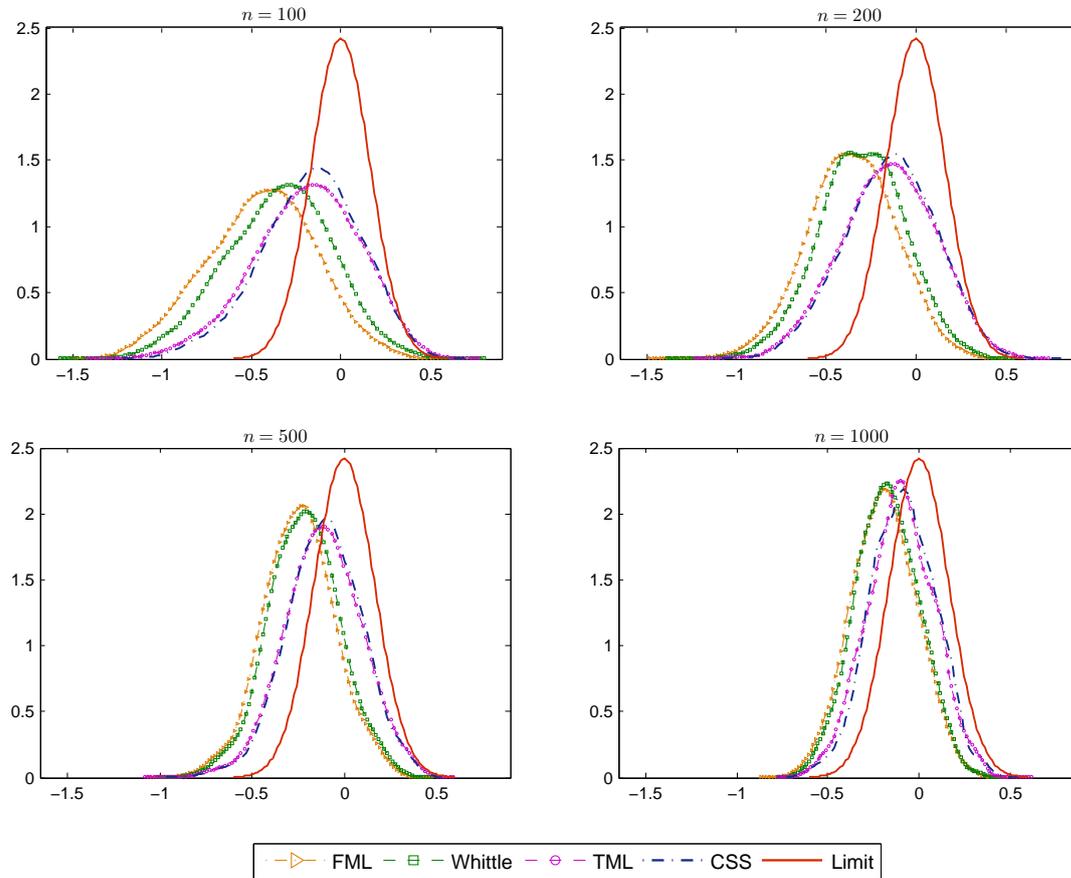

and $b$ as given in (30) evaluated at $\theta_0 = -0.3$ and $d^* = 0.1736$. Each panel in Figure 5 provides the kernel density estimate of the standardized statistic $\sqrt{n}(\widehat{d}_1 - d_1)$, under the four estimation methods, for a specific $n$ as labeled above each plot, plus the limiting distribution given in (36). In this case there is no clear visual differentiation between the time domain and frequency domain methods, for any sample size, and perhaps not surprisingly given the faster convergence rate in this case, all the methods produce finite sample distributions that match the limiting distribution reasonably well by the time $n = 1000$.

## 5.3 Finite sample bias and MSE of estimators of the pseudo-true parameter $d_1$ : known mean case

We supplement the graphical results in the previous section by documenting the finite sample bias and MSE of the four alternative estimators discussed in the previous section, in addition to the exact Whittle estimator, as estimators of the pseudo-true parameter $d_1$. The following

Figure 5: Kernel density of $\sqrt{n}\left(\widehat{d}_1 - d_1\right)$ for an $ARFIMA(0, d_0, 1)$ TDGP with $d_0 = 0.2$ and $\theta_0 = -0.3$, and an $ARFIMA(0, d, 0)$ MM, $d^* < 0.25$.

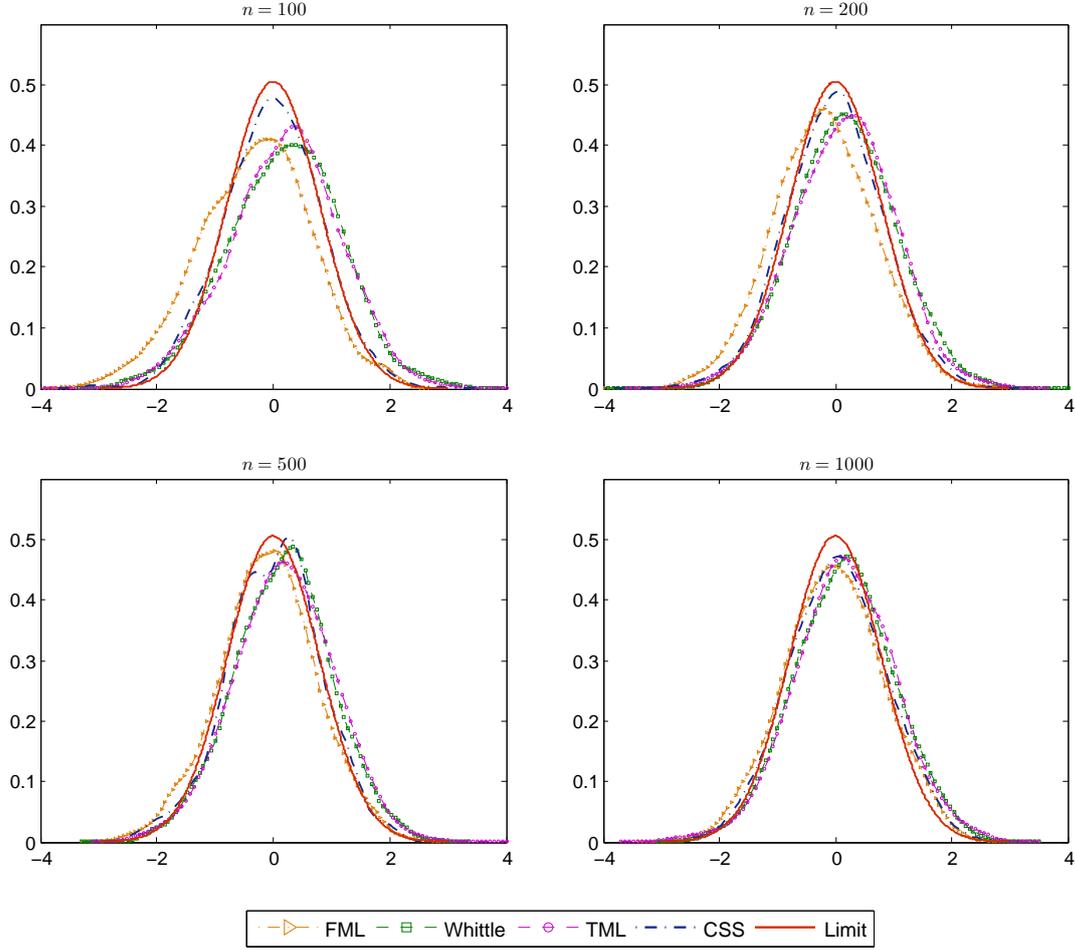

standard formulae,

$$\widehat{\text{Bias}}_0\left(\widehat{d}_1^{(i)}\right) = \frac{1}{R}\sum_{r=1}^{R}\widehat{d}_r^{(i)} - d_1 \tag{38}$$

$$\widehat{Var}_0\left(\widehat{d}_1^{(i)}\right) = \frac{1}{R}\sum_{r=1}^{R}\left(\widehat{d}_{1,r}^{(i)}\right)^2 - \left(\frac{1}{R}\sum_{r=1}^{R}\widehat{d}_{1,r}^{(i)}\right)^2 \tag{39}$$

$$\widehat{\text{MSE}}_0\left(\widehat{d}_1^{(i)}\right) = \widehat{\text{Bias}}_0^2 + \widehat{Var}_0\left(\widehat{d}_1^{(i)}\right) \tag{40}$$

$$\widehat{r.eff}_0\left(\widehat{d}_1^{(i)}, \widehat{d}_1^{(j)}\right) = \frac{\widehat{\text{MSE}}_0\left(\widehat{d}_1^{(i)}\right)}{\widehat{\text{MSE}}_0\left(\widehat{d}_1^{(j)}\right)}, \tag{41}$$

are applied to all five estimators $i, j = 1, ..., 5$. Since all empirical expectations and variances are evaluated under the TDGP, we make this explicit with appropriate subscript notation. Results for known mean are produced for Example 1 and Example 2 in Table 2 and 3 respectively, with selected additional results relevant to both examples recorded in Table 4. Values of $d^* = d_0 - d_1$ are documented across the key ranges, $d^* \lesseqgtr 0.25$, along with associated values for the MA coefficient in the TDGP, $\theta_0$. The minimum values of bias and MSE for each parameter setting are highlighted in bold face in all tables for each sample size, $n$.[8] Corresponding results for the case in which the mean is estimated are recorded in Section 5.4.

Consider first the bias and MSE results for Example 1 with $d_0 = 0.2$, as displayed in the middle panel of Table 2. As is consistent with the theoretical results (and the graphical illustration in the previous section) the bias and MSE of the four parametric estimators FML, Whittle, TML and CSS, show a clear tendency to decline as the sample size increases, for a fixed value of $\theta_0$. In addition, as $\theta_0$ declines in magnitude, and the MisM becomes closer to the TDGP, there is a tendency for the MSE values and the absolute values of the bias to decline. Importantly, the bias is *negative* for all four estimators, with the (absolute) bias of the two frequency domain estimators (FML and Whittle) being larger than that of the two time domain estimators. These results are consistent with the tendency of the standardized sampling distributions illustrated above to cluster, and for the frequency domain estimators to sit further to the left of zero than those of the time domain estimators, at least for the $d^* \geq 0.25$ cases. Again, as is consistent with the theoretical results, the rate of decline in the (absolute) bias and MSE of all estimators, as $n$ increases, is slower for $d^* \geq 0.25$ than for $d^* < 0.25$. The performance of the exact Whittle estimator (as we term it) is (almost) uniformly better than that of the Whittle estimator, but remains inferior to that of the two time domain estimators, with the CSS being clearly the superior estimator overall. The exact Whittle procedure mimics the other four methods in terms of the decline in both (absolute) bias and MSE as $n$ increases, providing numerical evidence that this frequency domain estimator is also consistent for $d_1$.[9]

As indicated by the results in the bottom panel of Table 2 for $d_0 = 0.4$, the impact of an increase in $d_0$ (for any given value of $d^*$ and $n$) is to often (but not uniformly) increase the bias and MSE of all (five) estimator of $d_1$. That is, the ability of the estimators to accurately estimate the pseudo-true parameter tends to decline (overall) as the long memory in the TDGP increases. In contrast, and with reference to the results in the top panel of the table for the antipersistent case, estimation accuracy tends to increase (as a general rule) as the memory in the TDGP declines. Nevertheless, the results in Table 2 show that the *relativities* between the estimators remain essentially the same for the different values of $d_0$, with the CSS estimator remaining preferable overall to all other estimators under mis-specification, and the FML estimator performing the worst of all.[10]

The results recorded in Table 3 for Example 2 illustrate that the presence of an AR term in the MisM means that more severe mis-specification can be tolerated. More specifically, in all (comparable) cases and for all estimators, the finite sample bias and MSE recorded in Table 3 tend to be smaller in (absolute) value than the corresponding values in Table 2. Results not presented here suggest, however, that when the value of $\theta_0$ is near zero, estimation under the MisM with an extraneous AR parameter causes an increase in (absolute) bias and MSE, relative to the case where the MisM is fractional noise (see also the following remark). With due consideration taken of the limited nature of the experimental design, these results suggest that the inclusion of some form of short memory dynamics in the estimated model – even if those dynamics are not of the correct form – acts as an insurance against more extreme mis-specification, but at the possible cost of a decline in performance when the consequences of mis-specification are not severe.

REMARK: When the parameter $\theta_0$ of the $ARFIMA(0, d_0, 1)$ TDGP equals zero the TDGP coincides with the $ARFIMA(0, d, 0)$ model and is nested within the $ARFIMA(1, d, 0)$

---

[8] Only that number which is smallest at the precision of 8 decimal places is bolded. Values highlighted with a '∗' are equally small to 4 decimal places.

[9] The tendency for the exact Whittle estimator to have (in particular) smaller finite sample bias than its inexact counterpart would seem to confirm the speculation in Chen and Deo (2006, Section 2) that the bias term ($\boldsymbol{\mu}_n$ in Theorem 3) may converge to zero more quickly for the first estimator than for the second.

[10] A slight caveat to this statement is that the superiority of the TML estimator over the exact Whittle is slightly less uniform when either $d = -0.2$ or $d = 0.4$. The overall similarity of the performance of these two estimators across the whole parameter space is, however, not surprising.

Table 2: Estimates of the bias and MSE of $\hat{d}_1$ for the FML, Whittle, Exact Whittle, TML and CSS estimators corresponding to Example 1 - TDGP: ARFIMA $(0, d_0, 1)$ vis-à-vis Mis-M: ARFIMA $(0, d, 0)$. Process mean, $\mu = 0$, is known.

|  |  |  |  | FML | | Whittle | | Exact Whittle | | TML | | CSS | |
| --- | --- | --- | --- | --- | --- | --- | --- | --- | --- | --- | --- | --- | --- |
| $d^*$ | $d_0$ | $\theta_0$ | $n$ | Bias | MSE | Bias | MSE | Bias | MSE | Bias | MSE | Bias | MSE |
| 0.3723 | -0.2 | -0.7 | 100 | -0.2416 | 0.0698 | -0.2085 | 0.0578 | -0.1317 | 0.0388 | -0.2055 | 0.0542 | -0.0941 | 0.0216 |
|  |  |  | 500 | -0.1332 | 0.0209 | -0.1243 | 0.0185 | -0.1107 | 0.0128 | -0.1175 | 0.0168 | -0.0631 | 0.0081 |
|  |  |  | 1000 | -0.1053 | 0.0129 | -0.1018 | 0.0121 | -0.1053 | 0.0115 | -0.0984 | 0.0113 | -0.0465 | 0.0050 |
| 0.2500 |  | -0.44 | 100 | -0.1425 | 0.0308 | -0.1068 | 0.0227 | -0.0823 | 0.0170 | -0.0837 | 0.0169 | -0.0401 | 0.0109 |
|  |  |  | 500 | -0.0539 | 0.0048 | -0.0502 | 0.0045 | -0.0460 | 0.0039 | -0.0414 | 0.0037 | -0.0192 | 0.0026 |
|  |  |  | 1000 | -0.0378 | 0.0024 | -0.0371 | 0.0024 | -0.0308 | 0.0022 | -0.0325 | 0.0021 | -0.0092 | 0.0014 |
| 0.1736 |  | -0.3 | 100 | -0.1071 | 0.0219 | -0.0645 | 0.0147 | -0.0468 | 0.0156 | -0.0464 | 0.0111 | -0.0215 | 0.0087 |
|  |  |  | 500 | -0.0345 | 0.0028 | -0.0272 | 0.0025 | -0.0154 | 0.0020 | -0.0212 | 0.0021 | -0.0087 | 0.0018 |
|  |  |  | 1000 | -0.0218 | 0.0013 | -0.0191 | 0.0012 | -0.0167 | 0.0011 | -0.0156 | 0.0010 | -0.0023 | 0.0009 |
| 0.3723 | 0.2 | -0.7 | 100 | -0.1781 | 0.0915 | -0.2466 | 0.0691 | -0.1968 | 0.0550 | -0.1748 | 0.0481 | -0.1427 | 0.0315 |
|  |  |  | 500 | -0.1354 | 0.0211 | -0.1308 | 0.0178 | -0.1017 | 0.0141 | -0.0916 | 0.0138 | -0.0798 | 0.0097 |
|  |  |  | 1000 | -0.1019 | 0.0141 | -0.0996 | 0.0127 | -0.0789 | 0.0097 | -0.0776 | 0.0103 | -0.0670 | 0.0065 |
| 0.2500 |  | -0.44 | 100 | -0.1515 | 0.0393 | -0.1184 | 0.0298 | -0.0795 | 0.0192 | -0.0650 | 0.0170 | -0.0577 | 0.0119 |
|  |  |  | 500 | -0.0544 | 0.0048 | -0.0487 | 0.0042 | -0.0307 | 0.0035 | -0.0257 | 0.0027 | -0.0241 | 0.0021 |
|  |  |  | 1000 | -0.0351 | 0.0023 | -0.0323 | 0.0021 | -0.0233 | 0.0018 | -0.0188 | 0.0015 | -0.0162 | 0.0012 |
| 0.1736 |  | -0.3 | 100 | -0.1082 | 0.0217 | -0.0712 | 0.0146 | -0.0428 | 0.0109 | -0.0340 | 0.0095 | -0.0330 | 0.0087 |
|  |  |  | 500 | -0.0318 | 0.0026 | -0.0251 | 0.0022 | -0.0184 | 0.0018 | **-0.0106** | 0.0017* | -0.0188 | 0.0017 |
|  |  |  | 1000 | -0.0184 | 0.0011 | -0.0149 | 0.0010 | -0.0077 | 0.0010 | **-0.0065** | 0.0009* | -0.0180 | 0.0009 |
| 0.3723 | 0.4 | -0.7 | 100 | -0.2786 | 0.0995 | -0.2456 | 0.0724 | -0.2165 | 0.0586 | -0.2210 | 0.0515 | -0.1957 | 0.0489 |
|  |  |  | 500 | -0.1598 | 0.0213 | -0.1287 | 0.0181 | -0.1142 | 0.0140 | -0.1347 | 0.0137 | -0.0871 | 0.0118 |
|  |  |  | 1000 | -0.1123 | 0.0157 | -0.0939 | 0.0143 | -0.0821 | 0.0113 | -0.0812 | 0.0121 | -0.0648 | 0.0117 |
| 0.2500 |  | -0.44 | 100 | -0.1903 | 0.0475 | -0.1659 | 0.0383 | -0.1086 | 0.0214 | -0.0911 | 0.0201 | -0.0550 | 0.0138 |
|  |  |  | 500 | -0.0796 | 0.0095 | -0.0550 | 0.0082 | -0.0303 | 0.0062 | -0.0249 | 0.0059 | -0.0224 | 0.0038 |
|  |  |  | 1000 | -0.0360 | 0.0048 | -0.0295 | 0.0042 | -0.0202 | 0.0038 | -0.0180 | 0.0035 | -0.0175 | 0.0025 |
| 0.1736 |  | -0.3 | 100 | -0.0990 | 0.0228 | -0.0843 | 0.0152 | -0.0509 | 0.0122 | -0.0422 | 0.0102 | -0.0321 | 0.0092 |
|  |  |  | 500 | -0.0407 | 0.0031 | -0.0276 | 0.0025 | -0.0188 | 0.0025 | -0.0129 | 0.0022 | **-0.0087** | 0.0019 |
|  |  |  | 1000 | -0.0172 | 0.0011 | -0.0163 | 0.0010 | -0.0089 | 0.0009 | -0.0077 | 0.0009 | **-0.0052** | 0.0008 |

Table 3: Estimates of the bias and MSE of $\widehat{d}_1$ for the FML, Whittle, Exact Whittle, TML and CSS estimators corresponding to Example 2 - TDGP: ARFIMA $(0, d_0, 1)$ vis-à-vis Mis-M: ARFIMA $(1, d, 0)$. Process mean, $\mu = 0$, is known.

| $d^*$ | $d_0$ | $\theta_0$ | $n$ | FML Bias | FML MSE | Whittle Bias | Whittle MSE | Exact Whittle Bias | Exact Whittle MSE | TML Bias | TML MSE | CSS Bias | CSS MSE |
|---|---|---|---|---|---|---|---|---|---|---|---|---|---|
| 0.2915 | -0.2 | -0.7 | 100 | -0.1174 | 0.0503 | -0.0809 | 0.0216 | -0.0756 | 0.0197 | -0.0632 | 0.0163 | -0.0526 | 0.0138 |
| | | | 500 | -0.0210 | 0.0115 | -0.0179 | 0.0109 | -0.0153 | 0.0098 | -0.0134 | 0.0085 | **-0.0122** | **0.0082** |
| | | | 1000 | -0.0105 | 0.0054 | -0.0088 | 0.0049 | -0.0079 | 0.0043 | -0.0070 | 0.0044 | **-0.0053** | **0.0042** |
| | 0.25 | -0.64 | 100 | -0.0868 | 0.0198 | -0.0770 | 0.0153 | -0.0641 | 0.0129 | -0.0507 | 0.0106 | -0.0420 | 0.0096 |
| | | | 500 | -0.0274 | 0.0039 | -0.0167 | 0.0026 | -0.0135 | 0.0031 | **-0.0113** | **0.0021** | -0.0121 | 0.0023 |
| | | | 1000 | -0.0166 | 0.0012 | -0.0128 | 0.0011 | -0.0101 | 0.0010 | -0.0094 | 0.0010 | **-0.0085** | **0.0010** |
| | 0.0148 | -0.3 | 100 | -0.0320 | 0.0123 | -0.0242 | 0.0099 | -0.0192 | 0.0103 | -0.0156 | 0.0069 | **-0.0095** | **0.0065** |
| | | | 500 | -0.0161 | 0.0025 | -0.0137 | 0.0031 | -0.0089 | 0.0029 | -0.0064 | 0.0013 | **-0.0030** | **0.0013** |
| | | | 1000 | -0.0094 | 0.0009 | -0.0075 | 0.0007 | -0.0055 | 0.0007 | -0.0037 | 0.0007* | **-0.0011** | **0.0007** |
| 0.2915 | 0.2 | -0.7 | 100 | -0.1612 | 0.0541 | -0.1169 | 0.0342 | -0.0982 | 0.0304 | -0.0950 | 0.0295 | **-0.0671** | **0.0236** |
| | | | 500 | -0.0679 | 0.0165 | -0.0604 | 0.0125 | -0.0557 | 0.0113 | -0.0454 | 0.0110 | **-0.0369** | **0.0089** |
| | | | 1000 | -0.0469 | 0.0089 | -0.0432 | 0.0071 | -0.0339 | 0.0070 | **-0.0250** | 0.0067 | -0.0303 | **0.0059** |
| | 0.25 | -0.64 | 100 | -0.1339 | 0.0279 | -0.0899 | 0.0175 | -0.0764 | 0.0145 | -0.0655 | 0.0138 | **-0.0457** | **0.0110** |
| | | | 500 | -0.0490 | 0.0041 | -0.0415 | 0.0030 | -0.0373 | 0.0032 | -0.0323 | 0.0026 | **-0.0230** | **0.0022** |
| | | | 1000 | -0.0316 | 0.0019 | -0.0281 | 0.0015 | -0.0228 | 0.0016 | -0.0181 | 0.0013 | **-0.0176** | **0.0011** |
| | 0.0148 | -0.3 | 100 | -0.0508 | 0.0139 | -0.0256 | 0.0086 | -0.0206 | 0.0075 | 0.0190 | 0.0067 | **-0.0082** | **0.0054** |
| | | | 500 | -0.0093 | 0.0027 | -0.0080 | 0.0019 | -0.0078 | 0.0018 | 0.0073 | 0.0016 | **-0.0004** | **0.0014** |
| | | | 1000 | -0.0036 | 0.0010 | -0.0023 | 0.0008 | -0.0020 | 0.0007 | 0.0067 | 0.0006* | **0.0003** | **0.0006** |
| 0.2915 | 0.4 | -0.7 | 100 | -0.2299 | 0.0639 | -0.1805 | 0.0419 | -0.1483 | 0.0389 | -0.1279 | 0.0372 | **-0.0699** | **0.0140** |
| | | | 500 | -0.1294 | 0.0197 | -0.1039 | 0.0150 | -0.0948 | 0.0134 | -0.0816 | 0.0126 | **-0.0294** | **0.0101** |
| | | | 1000 | -0.1089 | 0.0125 | -0.0632 | 0.0099 | -0.0521 | 0.0084 | -0.0462 | 0.0081 | **-0.0109** | **0.0069** |
| | 0.25 | -0.64 | 100 | -0.1396 | 0.0257 | -0.0979 | 0.0155 | -0.0804 | 0.0148 | -0.0692 | 0.0145 | **-0.0508** | **0.0103** |
| | | | 500 | -0.0455 | 0.0065 | -0.0342 | 0.0046 | -0.0315 | 0.0052 | -0.0294 | 0.0041 | **-0.0216** | **0.0033** |
| | | | 1000 | -0.0316 | 0.0027 | -0.0192 | 0.0021 | -0.0156 | 0.0020 | -0.0177 | 0.0018 | **-0.0122** | **0.0014** |
| | 0.0148 | -0.3 | 100 | -0.0650 | 0.0162 | -0.0422 | 0.0115 | -0.0368 | 0.0105 | 0.0246 | 0.0082 | **-0.0132** | **0.0067** |
| | | | 500 | -0.0205 | 0.0042 | -0.0133 | 0.0034 | -0.0099 | 0.0031 | 0.0079 | 0.0026 | **-0.0035** | **0.0023** |
| | | | 1000 | -0.0136 | 0.0021 | -0.0088 | 0.0018 | -0.0065 | 0.0015 | 0.0053 | 0.0014 | **-0.0017** | **0.0013** |

model. Thus the value $\theta_0 = 0$ is associated with a match between the TDGP and the model, at which point $d^* = 0$ and there is no mis-specification. That is, neither the $ARFIMA(0, d, 0)$ model estimated in Example 1, nor the $ARFIMA(1, d, 0)$ model estimated in Example 2, is mis-specified (according to our definition) when applied to an ARFIMA $(0, d_0, 0)$ TDGP, although the $ARFIMA(1, d, 0)$ model is incorrect in the sense of being over-parameterized. Table 4 presents the bias and MSE observed when there is such a lack of mis-specification. Under the correct specification of the $ARFIMA(0, d, 0)$ model the TML estimator is now superior, in terms of both bias and MSE. The relative accuracy of the TML estimator seen here is consistent with certain results recorded in Sowell (1992) and Cheung and Diebold (1994), in which the performance of the TML method (under a known mean, as is the case considered here) is assessed against that of various comparators under correct model specification. For the over-parameterized $ARFIMA(1, d, 0)$ model, however, the CSS estimator dominates once more. □

**Table 4** Estimates of the bias and MSE of $\widehat{d}_1$ for the FML, Whittle, Exact Whittle, TML and CSS estimators corresponding to TDGP: ARFIMA $(0, d_0, 0)$ $d_0 = 0.2$, $d^* = 0.0$ with the known process mean, $\mu = 0$.

| | FML | | Whittle | | Exact Whittle | | TML | | CSS | |
|---|---|---|---|---|---|---|---|---|---|---|
| $n$ | Bias | MSE | Bias | MSE | Bias | MSE | Bias | MSE | Bias | MSE |
| | | | | Correct $ARFIMA(0, d, 0)$ model | | | | | | |
| 100 | -0.0502 | 0.0113 | -0.0173 | 0.0102 | -0.0089 | 0.0099 | **0.0066** | **0.0087** | 0.0094 | 0.0096 |
| 500 | -0.0089 | 0.0015 | -0.0062 | 0.0014 | -0.0033 | 0.0014 | **0.0026** | **0.0013** | 0.0031 | 0.0014 |
| 1000 | -0.0045 | 0.0006* | -0.0037 | 0.0006* | -0.0028 | 0.0009 | **0.0016** | **0.0006** | 0.0025 | 0.0006* |
| | | | | Over-Parameterized $ARFIMA(1, d, 0)$ model | | | | | | |
| 100 | -0.0455 | 0.0177 | 0.0371 | 0.0121 | -0.0364 | 0.0122 | 0.0255 | 0.0107 | **0.0158** | **0.0087** |
| 500 | -0.0120 | 0.0065 | 0.0091 | 0.0049 | 0.0083 | 0.0039 | 0.0078 | 0.0043 | **0.0055** | **0.0037** |
| 1000 | -0.0074 | 0.0027 | 0.0055 | 0.0021 | 0.0040 | 0.0020 | 0.0034 | 0.0019 | **0.0028** | **0.0016** |

The results in Tables 2, 3 and 4 highlight that, in all but one case, the CSS estimator has the smallest MSE of all five estimators under mis-specification, and when there is no mis-specification but the model is over-parameterized, and that this result holds for all sample sizes considered. Indeed, the MSE results indicate that the CSS estimator is between about two and three times as efficient as the FML estimator (in particular) in the region of the parameter space ($d^* \geq 0.25$) in which both (absolute) bias and MSE are at their highest for all estimators. The absolute value of the bias of CSS is also the smallest in the vast majority of such cases, for all values of $d^*$. This almost universal superiority of the CSS method presumably reflects a certain in-built robustness of least squares methods.

### 5.4 Finite sample bias and MSE of estimators of the pseudo-true parameter $d_1$ : unknown mean case

Since the FML and Whittle estimators are both mean invariant the results recorded in the previous two sections for these two estimators are applicable to the unknown (zero) mean case without change. What will potentially alter, however, will be the performance of these two frequency domain estimators relative to that of the exact Whittle and time domain estimators when the unknown mean is also estimated, and it is that possibility that we explore in this section.

Table 5: Estimates of the bias and MSE of $\hat{d}_1$ for the Exact Whittle, TML and CSS estimators when the mean is also estimated using the sample mean. The values of d* for each design are given in Table 2.

| | | | Example 1: TDGP: $ARFIMA(0,d_0,1)$; MisM: $ARFIMA(0,d,0)$ | | | | | | | Example 2: TDGP: $ARFIMA(0,d_0,1)$; MisM: $ARFIMA(1,d,0)$ | | | | | | |
|---|---|---|---|---|---|---|---|---|---|---|---|---|---|---|---|---|
| | | | Exact Whittle | | TML | | CSS | | | Exact Whittle | | TML | | CSS | |
| n | $d_0$ | $\theta_0$ | Bias | MSE | Bias | MSE | Bias | MSE | $\theta_0$ | Bias | MSE | Bias | MSE | Bias | MSE |
| 100 | -0.2 | -0.7 | -0.1824 | 0.0468 | -0.2139 | 0.0581 | -0.1543 | 0.0318 | -0.7 | -0.1078 | 0.0356 | -0.0951 | 0.0364 | -0.0704 | 0.0269 |
| 500 | | | -0.1223 | 0.0254 | -0.1286 | 0.0193 | -0.0965 | 0.0122 | | -0.0487 | 0.0113 | -0.0332 | 0.0101 | -0.0265 | 0.0096 |
| 1000 | | | -0.0823 | 0.0148 | -0.1086 | 0.0133 | -0.0814 | 0.0084 | | -0.0224 | 0.0055 | -0.0168 | 0.0048 | -0.0099 | 0.0055 |
| 100 | | -0.44 | -0.1115 | 0.0202 | -0.1223 | 0.0247 | -0.0977 | 0.0186 | -0.64 | -0.0763 | 0.0212 | -0.0642 | 0.0184 | -0.0523 | 0.0142 |
| 500 | | | -0.0526 | 0.0051 | -0.0542 | 0.0049 | -0.0444 | 0.0039 | | -0.0285 | 0.0066 | -0.0173 | 0.0055 | -0.0119 | 0.0059 |
| 1000 | | | -0.0419 | 0.0025 | -0.0399 | 0.0026 | -0.0325 | 0.0021 | | -0.0155 | 0.0029 | -0.0116 | 0.0022 | -0.0094 | 0.0018 |
| 100 | | -0.3 | -0.0461 | 0.0168 | -0.0835 | 0.0165 | -0.0734 | 0.0143 | -0.3 | -0.0365 | 0.0128 | -0.0250 | 0.0110 | -0.0186 | 0.0099 |
| 500 | | | -0.0233 | 0.0030 | -0.0314 | 0.0027 | -0.0275 | 0.0024 | | -0.0125 | 0.0031 | -0.0127 | 0.0023 | -0.0092 | 0.0020 |
| 1000 | | | -0.0170 | 0.0016 | -0.0209 | 0.0013 | -0.0184 | 0.0011 | | -0.0092 | 0.0013 | -0.0072 | 0.0011 | -0.0045 | 0.0009 |
| 100 | 0.2 | -0.7 | -0.2850 | 0.0970 | -0.2765 | 0.0877 | -0.2527 | 0.0739 | -0.7 | -0.1567 | 0.0547 | -0.1478 | 0.0501 | -0.1216 | 0.0564 |
| 500 | | | -0.1509 | 0.0398 | -0.1425 | 0.0235 | -0.1363 | 0.0222 | | -0.0794 | 0.0222 | -0.0650 | 0.0155 | -0.0631 | 0.0146 |
| 1000 | | | -0.1326 | 0.0222 | -0.1210 | 0.0156 | -0.1121 | 0.0145 | | -0.0637 | 0.0096 | -0.0562 | 0.0082 | -0.0516 | 0.0080 |
| 100 | | -0.44 | -0.1728 | 0.0528 | -0.1677 | 0.0462 | -0.1527 | 0.0449 | -0.64 | -0.1245 | 0.0430 | -0.1171 | 0.0381 | -0.1037 | 0.0313 |
| 500 | | | -0.0748 | 0.0150 | -0.0605 | 0.0051 | -0.0582 | 0.0049 | | -0.0510 | 0.0087 | -0.0463 | 0.0072 | -0.0474 | 0.0068 |
| 1000 | | | -0.0433 | 0.0055 | -0.0425 | 0.0025 | -0.0389 | 0.0024 | | -0.0319 | 0.0035 | -0.0291 | 0.0015 | -0.0289 | 0.0015 |
| 100 | | -0.3 | -0.1242 | 0.0188 | -0.1168 | 0.0370 | -0.1018 | 0.0191 | -0.3 | -0.0585 | 0.0218 | -0.0473 | 0.0161 | -0.0382 | 0.0153 |
| 500 | | | -0.0462 | 0.0037 | -0.0357 | 0.0026 | -0.0300 | 0.0026 | | -0.0149 | 0.0029 | 0.0085 | 0.0019 | -0.0082 | 0.0019 |
| 1000 | | | -0.0244 | 0.0018 | -0.0166 | 0.0011 | -0.0153 | 0.0011 | | -0.0085 | 0.0012 | 0.0058 | 0.0008 | -0.0030 | 0.0008 |
| 100 | 0.4 | -0.7 | -0.2729 | 0.0927 | -0.2667 | 0.0838 | -0.2597 | 0.0781 | -0.7 | -0.2078 | 0.0819 | -0.1992 | 0.0732 | -0.2067 | 0.0643 |
| 500 | | | -0.1462 | 0.0258 | -0.1382 | 0.0244 | -0.1312 | 0.0199 | | -0.1321 | 0.0363 | -0.1246 | 0.0283 | -0.1106 | 0.0278 |
| 1000 | | | -0.1166 | 0.0194 | -0.0945 | 0.0162 | -0.1047 | 0.0155 | | -0.0832 | 0.0297 | -0.0712 | 0.0214 | -0.0664 | 0.0201 |
| 100 | | -0.44 | -0.1828 | 0.0578 | -0.1750 | 0.0479 | -0.1734 | 0.0393 | -0.64 | -0.1141 | 0.0453 | -0.1064 | 0.0362 | -0.1261 | 0.0324 |
| 500 | | | -0.0791 | 0.0148 | -0.0692 | 0.0108 | -0.0585 | 0.0099 | | -0.0491 | 0.0096 | -0.0349 | 0.0071 | -0.0388 | 0.0052 |
| 1000 | | | -0.0567 | 0.0057 | -0.0358 | 0.0056 | -0.0285 | 0.0049 | | -0.0282 | 0.0039 | -0.0163 | 0.0016 | -0.0202 | 0.0022 |
| 100 | | -0.3 | -0.0903 | 0.0192 | -0.0892 | 0.0182 | -0.0922 | 0.0183 | -0.3 | -0.0542 | 0.0262 | -0.0513 | 0.0172 | -0.0583 | 0.0167 |
| 500 | | | -0.0385 | 0.0028 | -0.0325 | 0.0027 | -0.0333 | 0.0029 | | -0.0257 | 0.0054 | 0.0172 | 0.0038 | -0.0152 | 0.0036 |
| 1000 | | | -0.0232 | 0.0012 | -0.0129 | 0.0009 | -0.0177 | 0.0010 | | -0.0178 | 0.0026 | 0.0101 | 0.0018 | -0.0090 | 0.0018 |

In Table 5 we record the bias and MSE obtained for the exact Whittle, TML and CSS estimators when the true mean of the process, $\mu$, is estimated using the sample mean; for both mis-specified examples, all three values of $d_0$, and all three sample sizes. Properties observed in the previous section, such as the decline in bias and MSE with an increase in sample size, for a given $\theta_0$, and the overall decline in MSE and (absolute) bias as the estimated model becomes less mis-specified, continue to obtain in Table 5. However, the magnitudes of the bias and MSE figures for all three estimators are now virtually always higher than the corresponding figures in Tables 2 and 3. As a consequence of this, the time domain estimators lose their relative superiority and no longer uniformly dominate the frequency domain techniques. Instead, the Whittle estimator outperforms all four of the other estimators overall (including its exact counterpart), and almost uniformly in the (true) long memory cases ($d_0 = 0.2, 0.4$). Table 6 records the outcomes obtained for the exact Whittle, TML and CSS estimators under the correct and over-parameterized specifications when the mean is estimated. Comparing Table 6 with Table 4 we find (once again) that the Whittle estimator now dominates all other estimators. As the sample size increases the differences between all comparable results for the known and estimated mean cases become less marked, in accordance with the consistency of the estimated mean for the true (zero) mean.[11]

**Table 6** Estimates of the bias and MSE of $\widehat{d_1}$ for the FML, Whittle, WI, TML and CSS estimators corresponding to TDGP: ARFIMA $(0, d_0, 0)$ $d_0 = 0.2$, $d^* = 0.0$. The unknown mean is estimated using the sample mean.

|      | Exact Whittle | | TML | | CSS | |
|------|------|------|------|------|------|------|
| $n$  | Bias | MSE  | Bias | MSE  | Bias | MSE  |
| *Correct $ARFIMA(0, d, 0)$ model* | | | | | | |
| 100  | -0.0639 | 0.0184 | -0.0542 | 0.0149 | -0.0585 | 0.0158 |
| 500  | -0.0163 | 0.0036 | -0.0111 | 0.0027 | -0.0156 | 0.0032 |
| 1000 | -0.0052 | 0.0018 | -0.0047 | 0.0006 | -0.0048 | 0.0006 |
| *Over-Parameterized $ARFIMA(1, d, 0)$ model* | | | | | | |
| 100  | -0.0825 | 0.0282 | -0.0758 | 0.0224 | -0.0701 | 0.0022 |
| 500  | -0.0267 | 0.0081 | -0.0195 | 0.0077 | -0.0188 | 0.0075 |
| 1000 | -0.0102 | 0.0031 | -0.0087 | 0.0029 | -0.0085 | 0.0029 |

# 6 Discussion

This paper presents theoretical and simulation-based results relating to the estimation of mis-specified models for fractionally integrated processes. We show that under mis-specification four classical parametric estimation methods, frequency domain maximum likelihood (FML), Whittle, time domain maximum likelihood (TML) and conditional sum of squares (CSS) converge to the same pseudo-true parameter value. Consistency of the four estimators for the pseudo-true value is proved for fractional exponents of both the true and estimated models in the long memory, short memory and antipersistent ranges. A general closed-form solution for the limiting criterion function for the four alternative estimators is derived in the case of ARFIMA models. This enables us to demonstrate the link between any form of mis-specification of the short memory dynamics and the difference between the true and pseudo-true values of the fractional index, $d$, and, hence, to the resulting (asymptotic) distributional properties of the estimators, having proved that the estimators are asymptotically equivalent.

---

[11] The results recorded here regarding the performance of the different estimators in the unknown mean case parallel the qualitative conclusions drawn by Nielsen and Frederiksen (2005) for correctly specified models. Note also that very similar results are obtained if the sample mean is replaced by a feasible (plug in) version of the (asymptotically) best linear unbiased estimator (BLUE) of $\mu$.

The finite sample performance of all four estimators is then documented. The extent to which the finite sample distributions mimic the (numerically specified) asymptotic distributions is displayed. In the case of more extreme mis-specification, and conditional on the mean of the process being known, the pairs of time domain and frequency domain estimators tend to cluster together for smaller sample sizes, with the former pair mimicking the asymptotic distributions more closely. Further, bias and mean squared error (MSE) calculations demonstrate the superiority overall of the CSS estimator, under mis-specification, and the distinct inferiority of the FML estimator – as estimators of the pseudo-true parameter for which they are both consistent. Numerical results for the time domain estimators in the case where the unknown mean is estimated tell a different story, however, with the Whittle estimator being the superior finite sample performer overall. Numerical results presented for an exact version of the Whittle estimator show a slight superiority over the (approximate) version of the Whittle procedure, in the case where the mean is known; however, the overall ranking of the two methods is reversed when the mean is estimated, with the exact Whittle method not sharing the mean-invariance property of its inexact counterpart.

There are several interesting issues that arise from the results that we have established, including the following: First, although the known (zero) mean assumption is inconsequential for the FML and Whittle estimators, this is not the case for the exact Whittle and time domain estimators, as our bias and mean squared error experimental results obtained using demeaned data show. The deterioration in the overall performance of the exact Whittle and time domain estimators once the estimation of $\mu$ plays a role in their computation might have been anticipated since the rate of convergence of the sample average to $\mu$ is $n^{1/2-d_0}$ (Hosking, 1996, Theorem 8), and thus slower the larger the value of $d_0$. Similarly, estimation of $\mu$ will impact on the limiting distribution of the time domain estimators – because the rate of convergence of the estimators when the true mean is known is $n^{1-2d^*}/\log n$ when $d^* = d_0 - d_1 > 0.25$, $(n/\log\log\log n)^{1/2}$ when $d^* = 0.25$, and $\sqrt{n}$ otherwise – something that we have not pursued for the current paper, but is the subject of other ongoing research. Second, the extension of our results to non-stationary cases will facilitate the consideration of a broader range of circumstances. To some extent non-stationary values of $d$ might be covered by means of appropriate pre-filtering, for example, the use of first-differencing when $d_0 \in [1/2, 3/2)$, but this would require prior knowledge of the structure of the process and opens up the possibility of a different type of mis-specification from the one we have considered here. Explicit consideration of the interval $d \in [0, 3/2)$, say, allowing for both stationary and non-stationary cases perhaps offers a better approach as prior knowledge of the characteristics of the process would then be unnecessary. The latter also seems particularly relevant given that estimates near the boundaries $d = 0.5$ and $d = 1$ are not uncommon in practice. Previous developments in the analysis of non-stationary fractional processes (see, inter alios, Beran, 1995; Tanaka, 1999; Velasco, 1999; Velasco and Robinson, 2000) might offer a sensible starting point for such an investigation. Third, our limiting distribution results can be used in practice to conduct inference on the long memory and other parameters after constructing obvious smoothed periodogram consistent estimates of $\mathbf{B}$, $\boldsymbol{\mu}_n$, $\overline{\Lambda}_{dd}$ and $\boldsymbol{\Lambda}$. But which situation should be assumed in any particular instance, $d^* > 0.25$, $d^* = 0.25$ or $d^* < 0.25$, may be a moot point. Fourth, the relationships between the bias and MSE of the parametric estimators of $d_1$ (denoted respectively below by Bias_$d_1$ and MSE_$d_1$), and the bias and MSE as estimators of the *true* value $d_0$, (Bias_$d_0$ and MSE_$d_0$ respectively) can be expressed simply as follows:

$$\begin{aligned} \text{Bias\_}d_0 &= E_0(\widehat{d}_1) - d_0 \\ &= \left[E_0(\widehat{d}_1) - d_1\right] + (d_1 - d_0) \\ &= Bias\_d_1 - d^*, \end{aligned}$$

where we recall, $d^* = d_0 - d_1$, and

$$\begin{aligned} \text{MSE\_}d_0 &= E_0\left(\widehat{d}_1 - d_0\right)^2 \\ &= E_0\left(\widehat{d}_1 - E_0(\widehat{d}_1)\right)^2 + \left[E_0(\widehat{d}_1) - d_0\right]^2 \\ &= E_0\left(\widehat{d}_1 - E_0(\widehat{d}_1)\right)^2 + \left[[E_0(\widehat{d}_1) - d_1] - d^*\right]^2 \\ &= E_0\left(\widehat{d}_1 - E_0(\widehat{d}_1)\right)^2 + \left[E_0(\widehat{d}_1) - d_1\right]^2 + d^{*2} - 2d^*\left[E_0(\widehat{d}_1) - d_1\right] \\ &= MSE\_d_1 + d^{*2} - 2d^* Bias\_d_1. \end{aligned}$$

Hence, if Bias_$d_1$ is the same sign as $d^*$ at any particular point in the parameter space, then the bias of a mis-specified parametric estimator *as an estimator of* $d_0$, may be less (in absolute value) than its bias as an estimator of $d_1$, depending on the magnitude of the two quantities. Similarly, MSE_$d_0$ may be less than MSE_$d_1$ if Bias_$d_1$ and $d^*$ have the same sign, with the final result again depending on the magnitude of the two quantities. These results imply that it is possible for the ranking of mis-specified parametric estimators to be altered, once the reference point changes from $d_1$ to $d_0$. This raises the following questions: Does the dominance of the CSS estimator (within the parametric set of estimators) – and in the known mean case – still obtain when the true value of $d$ is the reference value? And more critically from a practical perspective; Are there circumstances where a mis-specified parametric estimator out-performs semi-parametric alternatives in finite samples, the lack of consistency (for $d_0$) of the former notwithstanding? Such topics remain the focus of current and ongoing research.

# References

Beran, J. (1994). *Statistics for Long-Memory Processes*, vol. 61 of *Monographs on Statistics and Applied Probability*. Chapman and Hall, New York.

Beran, J. (1995). Maximum likelihood estimation of the differencing parameter for invertible short and long memory autoregressive integrated moving average models. *Journal of the Royal Statistical Society* **B 57** 654–672.

Boes, D. C., Davis, R. A. and Gupta, S. N. (2006). Parameter estimation in low order fractionally differenced ARMA processes. *Stochastic Hydrology and Hydraulics* **3** 97–110.

Brockwell, P. J. and Davis, R. A. (1991). *Time Series: Theory and Methods*. Springer Series in Statistics. Springer-Verlag, New York, 2nd ed.

Cavaliere, G., Nielsen, M. O. and Taylor, A. M. R. (2017). Quasi-Maximum Likelihood Estimation and Bootstrap Inference in Fractional Time Series Models with Heteroskedasticity of Unknown Form. *Journal of Econometrics* **198** 165 – 188.

Chen, W. W. and Deo, R. S. (2006). Estimation of mis-specified long memory models. *Journal of Econometrics* **53** 257 – 281.

CHEUNG, Y. W. and DIEBOLD, F. X. (1994). On maximum likelihood estmation of the differencing parameter of fractionally integrated noise with unknown mean. *Journal of Econometrics* **62** 301–316.

CRATO, N. and TAYLOR, H. M. (1996). Stationary persistent time series misspecified as nonstationary ARIMA. *Statistical Papers* **37** 215–266.

DAHLHAUS, R. (1989). Efficient parameter estimation for self-similar processes. *Annals of Statistics* **17** 1749–1766.

DOORNIK, J. A. and OOMS, M. (2001). Computational aspects of maximum likelihood estimation of autoregressive fractionally integrated moving average models. *Computational Statistics & Data Analysis* **42** 333–348. Also a 2001 Nuffield discussion paper.

FOX, R. and TAQQU, M. S. (1986). Large sample properties of parameter estimates for strongly dependent stationary gaussian time series. *Annals of Statistics* **14** 517–532.

GIRAITIS, L. and SURGAILIS, D. (1990). A central limit theorem for quadratic forms in strongly dependent linear variables and its application to asymptotical normality of Whittle's estimate. *Probability Theory and Related Fields* **86** 87–104.

GRADSHTEYN, I. S. and RYZHIK, I. M. (2007). *Tables of Integrals, Series and Products*. Academic Press, Sydney.

GRENANDER, U. and SZEGO, G. (1958). *Toeplitz Forms and Their Application*. University of California Press, Berkeley.

HANNAN, E. J. (1973). The asymptotic theory of linear time series models. *Advances in Applied Probability* **10** 130–145.

HASSLER, U. (1994). (Mis)Specification of long memory in seasonal time series. *Journal of Time Series Analysis* **14** 19–30.

HUALDE, J. and ROBINSON, P.M. (2011). Gaussian pseudo-maximum likelihood estimation of fractional time series models. Annals of Statistics. **39**, 3152–3181.

HURVICH, C. M. and BELTRAO, K. I. (1993). Asymptotics for the low-frequency ordinates of the periodogram of a long memory time series. *Journal of Time Series Analysis* **14** 455–472.

HOSKING, J. R. M. (1996). Asymptotic distributions of the sample mean, autocovariances, and autocorrelations of long memory time series. *Journal of Econometrics* **73** 261–284.

LAHIRI, S. N. (2003). A necessary and sufficient condition for asymptotic independence of discrete fourier transforms under short- and long-range dependence. *Annals of Statistics* **31** 613-641

LIEBERMAN, O. and PHILLIPS, P. C. B. (2004). Error bounds and asymptotic expansions for toeplitz product functionals of unbounded spectra. *Journal of Time Series Analysis* **25** 732–753.

MOULINES, E. and SOULIER, P. (1999). Broadband log-periodogram regression of times series with long-range dependence. *Annals of Statistics* **27** 1415-1439.

NIELSEN, M. O. and FREDERIKSEN, P. H. (2005). Finite sample comparison of parametric, semiparametric, and wavelet estimators of fractional integration. *Econometric Reviews* **24** 405–443.

ROBINSON, P. M. (1995). Gaussian semiparametric estimation of long-range dependence *Annals of Statistics*, **23** 1630-1661.

ROBINSON, P. M. (2006). *Time Series and Related Topics: In Memory of Ching-Zong Wei*, vol. 52 of *IMS Lecture Notes–Monograph Series*, chap. Conditional-Sum-of-Squares Estimation of Models for Stationary Time Series with Long Memory. Institue of Mathematical Statistics, Beachwood, 130–137.

ROBINSON, P. M. (2014). The estimation of misspecified long memory models *Journal of Econometrics*, **178** 225–230.

SOWELL, F. (1992). Maximum likelihood estmation of stationary univariate fractionally integrated time series models. *Journal of Econometrics* **53** 165–188.

TANAKA, K. (1999). The nonstationary fractional unit root. *Econometric Theory* **15** 549–582.

VELASCO, C. (1999). Gaussian semiparametric estimation of non-stationary time series. *Journal of Time Series Analysis* **20** 87–127.

VELASCO, C. and ROBINSON, P. M. (2000). Whittle pseudo-maximum likelihood estimation for nonstationary time series. *Journal of the American Statistical Association* **95:452** 1229–1243.

WHITTLE, P. (1953). Estimation and information in stationary time series. *Ark. Mat.* **2** 423–434.

WHITTLE, P. (1984). *Prediction and Regulation By Linear Least Squares Methods*. Basil Blackwell.

YAJIMA, Y. (1992). Asymptotic properties of estimates in incorrect ARMA models for long-memory time series. In D. Brillinger, P. Caines, J. Geweke, E. Parzen, M. Rosenblatt and M. Taqqu, eds., *New Directions in Time Series Analysis, Part II*, vol. 46 of *IMA Volumes in Mathematics and Its Applications*. Springer-Verlag, New York, 375–382.

# A  Proofs

## A.1  Proof of Lemmas 1, 2, 3 and 4

### A.1.1  Proof of Lemma 1:

The proof of the lemma uses a method that parallels that employed by Fox and Taqqu in the proof of their Lemma 1 (see Fox and Taqqu, 1986, pages 523-524), which in turn employs an argument first developed by Hannan in the proof of his Lemma 1 (See Hannan, 1973, pages 133-134.). To describe the approach, set

$$c_n(\tau) = c_n(-\tau) = \frac{1}{n}\sum_{t=1}^{n-\tau} y_t y_{t+\tau}, \quad \tau \geq 0,$$

and let

$$k_M(\boldsymbol{\eta}, \lambda) = \sum_{r=-M}^{M} \kappa(r)\left(1 - \frac{|r|}{M}\right)\exp(i\lambda r)$$

denote the Cesaro sum of the first $M$ terms of the Fourier series of $(f_1(\boldsymbol{\eta},\lambda)+\nu_f)^{-1}$ where $M$ is chosen such that $|(f_1(\boldsymbol{\eta},\lambda) + \nu_f)^{-1} - k_M(\boldsymbol{\eta},\lambda))| < \epsilon$ uniformly in $\boldsymbol{\eta} \in \mathbb{E}_\delta^0$. Then following the same steps as in the derivation presented in Hannan (1973, pages 133-134) we have

$$\left|\frac{4\pi}{n}\sum_{j=1}^{[n/2]} I(\lambda_j)\left\{(f_1(\boldsymbol{\eta},\lambda_j) + \nu_f)^{-1} - k_M(\boldsymbol{\eta},\lambda_j)\right\}\right| < \epsilon c_n(0),$$

and

$$\lim_{n\to\infty}\left|\frac{4\pi}{n}\sum_{j=1}^{[n/2]} I(\lambda_j)k_M(\boldsymbol{\eta},\lambda_j) - \sum_{r=-M}^{M} \kappa(r)\left(1-\frac{|r|}{M}\right)\gamma_0(r)\right| = 0$$

almost surely, the latter result since $I(\lambda) = (2\pi)^{-1}\sum_{r=1-n}^{n-1} c_n(r)\exp(-i\lambda r)$ and $c_n(r)$ converges to $\gamma_0(r)$ almost surely by ergodicity. Moreover,

$$\sum_{r=-M}^{M} \kappa(r)\left(1-\frac{|r|}{M}\right)\gamma_0(r) = \frac{\sigma_0^2}{2\pi}\int_{-\pi}^{\pi} f_0(\lambda)k_M(\boldsymbol{\eta},\lambda)d\lambda$$

differs from the required limiting value by a quantity bounded by $\epsilon\gamma_0(0)$, from which the desired result follows because $\epsilon$ is arbitrary.

An alternative proof of this lemma can be obtained by extending the arguments adopted by (Brockwell and Davis, 1991, Section 10.8.2, pages 378-379), in the proof of their Proposition 10.8.2, to the stationary fractional case, as suggested in Brockwell and Davis (1991, page 528).

### A.1.2 Proof of Lemma 2:

The proof parallels the proof of Lemma 1, only now we use the Cesaro sum of $M$ terms of the Fourier series of $h_1(\boldsymbol{\eta},\lambda)^{-1}$. Denote this sum by $c_M(\boldsymbol{\eta},\lambda) > 0$. Since by construction $h_1(\boldsymbol{\eta},\lambda) > 0$, $M$ can be chosen so that $|h_1(\boldsymbol{\eta},\lambda)^{-1} - c_M(\boldsymbol{\eta},\lambda)| < \varepsilon$ uniformly on $\mathbb{E}_\delta^0$ since the Cesaro sum converges uniformly in $(\boldsymbol{\eta},\lambda)$ for $\boldsymbol{\eta} \in \mathbb{E}_\delta^0$. Once again the detailed steps follow Hannan (1973, pages 133-134), as above, or (Brockwell and Davis, 1991, Section 10.8.2, pages 378-379)

### A.1.3 Proof of Lemma 3:

Observe that $f_1(\boldsymbol{\eta},\lambda) > 0$ when $d \geq 0$ and hence for $\delta$ sufficiently small we have $h_1(\boldsymbol{\eta},\lambda) = f_1(\boldsymbol{\eta},\lambda)$ for all $\lambda \in [-\pi,\pi]$. It follows immediately from Lemma 2 that $\lim_{n\to\infty}|Q_n^{(1)}(\boldsymbol{\eta}) - Q(\boldsymbol{\eta})| = 0$ almost surely and uniformly in $\boldsymbol{\eta}$ on $\mathbb{E}_\delta^0$ when $d \geq 0$. We have thus established Lemma 3 in the case where $d \geq 0$, (*cf.* Chen and Deo, 2006, Lemma 2). To establish that Lemma 3 also holds on $\mathbb{E}_\delta^0$ when $d < 0$, observe that Lemma 1 implies that $Q(\boldsymbol{\eta})$ provides a limit inferior for $Q_n^{(1)}(\boldsymbol{\eta})$ and it therefore only remains for us to establish that $Q(\boldsymbol{\eta})$ also provides a limit superior for $Q_n^{(1)}(\boldsymbol{\eta})$ on $\boldsymbol{\eta} \in \mathbb{E}_\delta^0$ when $d < 0$.

In the latter case $f_1(\boldsymbol{\eta},\lambda) = |\lambda|^{2|d|}L(\lambda)$ where $L(\lambda)$ is slowly varying and bounded as $\lambda \to 0$ and there exists an $\epsilon \in (0, 2|d|)$ and a $K > 0$, that may depend on $\epsilon$, such that $f_1(\boldsymbol{\eta},\lambda) = |\lambda|^{2|d|}K|\lambda|^{-\epsilon}$. We therefore have that $f_1(\boldsymbol{\eta},\lambda) > K|\lambda|^{2|d|}$ when $|\lambda| < 1$ and $h_1(\boldsymbol{\eta},\lambda) \neq f_1(\boldsymbol{\eta},\lambda)$ whenever $\lambda < (K^{-1}\delta)^{1/(2|d|-\epsilon)}$, from which it follows that

$$Q_n^{(1)}(\boldsymbol{\eta}) \leq \frac{2\pi}{n}\sum_{j=k_\delta+1}^{\lfloor n/2 \rfloor} \frac{I(\lambda_j)}{h_1(\boldsymbol{\eta},\lambda_j)} + \frac{1}{K}\left(\frac{2\pi}{n}\right)^{1-2|d|}\sum_{j=1}^{k_\delta} I(\lambda_j) \qquad (A.1)$$

where $k_\delta = \lfloor (K^{-1}\delta)^{1/(2|d|-\epsilon)}(2\pi/n)\rfloor + 1$. The inequality in (A.1) follows because for all $\lambda_j < (K^{-1}\delta)^{1/(2|d|-\epsilon)} < 2\pi k_\delta/n$ we have

$$\left(\frac{h_1(\boldsymbol{\eta},\lambda_j)}{f_1(\boldsymbol{\eta},\lambda_j)} - 1\right) \leq \left(\frac{\delta}{K}\left(\frac{n}{2\pi}\right)^{2|d|} - 1\right)$$

and

$$Q_n^{(1)}(\boldsymbol{\eta}) = \frac{2\pi}{n} \sum_{j=1}^{\lfloor n/2 \rfloor} \frac{I(\lambda_j)}{h_1(\boldsymbol{\eta}, \lambda_j)} + \frac{2\pi}{n} \sum_{j=1}^{k_\delta} I(\lambda_j) \left( \frac{1}{f_1(\boldsymbol{\eta}, \lambda_j)} - \frac{1}{h_1(\boldsymbol{\eta}, \lambda_j)} \right)$$

$$\leq \frac{2\pi}{n} \sum_{j=1}^{\lfloor n/2 \rfloor} \frac{I(\lambda_j)}{h_1(\boldsymbol{\eta}_1, \lambda_j)} + \frac{2\pi}{n} \sum_{j=1}^{k_\delta} \frac{I(\lambda_j)}{h_1(\boldsymbol{\eta}, \lambda_j)} \left( \frac{\delta}{K} \left( \frac{n}{2\pi} \right)^{2|d|} - 1 \right)$$

$$= \frac{2\pi}{n} \sum_{j=k_\delta+1}^{\lfloor n/2 \rfloor} \frac{I(\lambda_j)}{h_1(\boldsymbol{\eta}, \lambda_j)} + \frac{1}{K} \left( \frac{2\pi}{n} \right)^{1-2|d|} \sum_{j=1}^{k_\delta} I(\lambda_j).$$

Applying Lemma 2 to the first term on the right hand side in (A.1) gives a limit of

$$\frac{\sigma_0^2}{2\pi} \int_{(K^{-1}\delta)^{1/(2|d|-\epsilon)}}^{\pi} \frac{f_0(\lambda)}{f_1(\boldsymbol{\eta}_1, \lambda)} d\lambda.$$

Similarly

$$\lim_{n\to\infty} \frac{2\pi}{n} \sum_{j=1}^{k_\delta} I(\lambda_j) = \frac{\sigma_0^2}{2\pi} \int_0^{(K^{-1}\delta)^{1/(2|d|-\epsilon)}} f_0(\lambda) d\lambda = \frac{\sigma_0^2}{2\pi} f_0(\lambda')(K^{-1}\delta)^{1/(2|d|-\epsilon)}$$

for some $\lambda' \in [0, (K^{-1}\delta)^{1/(2|d|-\epsilon)}]$ by the first mean value theorem for integrals. Setting $\delta = (2\pi)^{2|d|-\epsilon}/n^p$ where $p > 2|d| - \epsilon$, we find that

$$\frac{1}{K} \left( \frac{2\pi}{n} \right)^{1-2|d|} \sum_{j=1}^{k_\delta} I(\lambda_j) \sim \frac{1}{K} \left( \frac{n}{2\pi} \right)^{2|d|} \frac{\sigma_0^2}{2\pi} f_0(\lambda') \frac{2\pi k_\delta}{n}$$

$$\sim \frac{1}{K} \left( \frac{2\pi}{n} \right)^{1-2|d|} \frac{\sigma_0^2}{2\pi} f_0(\lambda') \left( \frac{1}{n} \right)^{\frac{p-2|d|+\epsilon}{(2|d|-\epsilon)}}$$

and hence we can conclude that

$$\limsup_{n\to\infty} Q_n^{(1)}(\boldsymbol{\eta}) \leq Q(\boldsymbol{\eta})$$

uniformly in $\boldsymbol{\eta} \in \mathbb{E}_\delta^0$, as required.

### A.1.4 Proof of Lemma 4:

Let $L_1(\boldsymbol{\eta}, \lambda) = \lambda^{2d} f_1(\boldsymbol{\eta}, \lambda)$ and suppose that $\boldsymbol{\eta} \in \overline{\mathbb{E}}_{\delta 1}^0 \cup \overline{\mathbb{E}}_{\delta 2}^0 \neq \emptyset$. Then

$$\liminf_{n\to\infty} Q_n^{(1)}(\boldsymbol{\eta}) = \liminf_{n\to\infty} \frac{2\pi}{n} \sum_{j=1}^{\lfloor n/2 \rfloor} \frac{I(\lambda_j)}{f_1(\boldsymbol{\eta}, \lambda_j)}$$

$$= \liminf_{n\to\infty} \frac{2\pi}{n} \sum_{j=1}^{\lfloor n/2 \rfloor} \frac{I(\lambda_j)\lambda_j^{2d}}{L_1(\boldsymbol{\eta}, \lambda_j)}$$

$$\geq \liminf_{n\to\infty} \frac{(2\pi)^{-2\delta}}{n} \sum_{j=1}^{\lfloor n/2 \rfloor} \frac{I(\lambda_j)}{L_1(\boldsymbol{\eta}, \lambda_j)\lambda_j^{1-2(d_0+\delta)}}, \quad (A.2)$$

where the inequality in (A.2) arises because for all $\boldsymbol{\eta} \in \overline{\mathbb{E}}_{\delta 1}^0 \cup \overline{\mathbb{E}}_{\delta 2}^0$ we have $(d_0 - d) > 0.5 - \delta$ and it follows that $\lambda_j^{-2(d_0-d)} \geq (2\pi)^{-(2\delta+1)} \lambda_j^{2\delta-1}$ for all $\lambda_j = 2\pi j/n$, $j = 1, \ldots, \lfloor n/2 \rfloor$.

Applying Lemma 1 and Lemma 2 to (A.2) by replacing $f_1(\boldsymbol{\eta}, \lambda_j)$ by $L_1(\boldsymbol{\eta}, \lambda)\lambda^{1-2(d_0+\delta)}$, and then letting the constant $\nu_f > 0$ in the lemmas approach zero, it follows from Fatou's theorem that

$$\lim_{n\to\infty}\left|\frac{(2\pi)^{-2\delta}}{n}\sum_{j=1}^{\lfloor n/2 \rfloor}\frac{I(\lambda_j)}{L_1(\boldsymbol{\eta},\lambda_j)\lambda_j^{1-2(d_0+\delta)}} - \frac{1}{(2\pi)^{2\delta+1}}\int_0^\pi \frac{(\sigma_0^2/2\pi)f_0(\lambda)\lambda^{2d_0}}{L_1(\boldsymbol{\eta},\lambda)\lambda^{1-2\delta}}d\lambda\right| = 0\,,$$

wherein we recognise that $0 \leq 1-2(d_0+\delta) \leq 2(1-2\delta)$ and that $L_1(\boldsymbol{\eta}, \lambda) = (\sigma_1^2/2\pi)g_1(\boldsymbol{\beta},\lambda)\operatorname{sinc}(\lambda/2)^{-2d}$ and $(\sigma_0^2/2\pi)f_0(\lambda)\lambda^{2d_0} = (\sigma_0^2/2\pi)g_0(\lambda)\operatorname{sinc}(\lambda/2)^{-2d_0}$ where $\operatorname{sinc}(x) = \sin(x)/x$, the cardinal sine function. Since $2/\pi \leq \operatorname{sinc}(\lambda/2) \leq 1$ for $0 \leq \lambda \leq \pi$, it follows from Assumption A.3 and *Conditions A* that there exists a finite positive constant $R$ such that

$$\begin{aligned}\frac{1}{(2\pi)^{2\delta+1}}\int_0^\pi \frac{(\sigma_0^2/2\pi)f_0(\lambda)\lambda^{2d_0}}{L_1(\boldsymbol{\eta},\lambda)\lambda^{1-2\delta}}d\lambda &\geq \frac{R}{(2\pi)^{2\delta+1}}\cdot\int_0^\pi \lambda^{2\delta-1}d\lambda \\ &= \frac{R}{(2\pi)^{2\delta+1}}\cdot\frac{\pi^{2\delta}}{2\delta} \\ &\geq \frac{R}{8\pi}\cdot\frac{1}{\delta}\,.\end{aligned} \quad (A.3)$$

The statements in Lemma 4 now follow from (A.3), directly in the case of $\boldsymbol{\eta} \in \overline{\mathbb{E}}^0_{\delta 1}$, and for $\boldsymbol{\eta} \in \overline{\mathbb{E}}^0_{\delta 2}$ on setting $\delta < R/(8\pi C)$ and letting $\delta \to 0$ as $C \to \infty$. ∎

## A.2 Proof of Proposition 1

Let $\boldsymbol{\eta}_n$ denote a sequence in $\mathbb{E}^0_\delta$ that converges to $\boldsymbol{\eta}$. For any $\nu_f > 0$ we have

$$\left|\frac{1}{f_1(\boldsymbol{\eta}_n,\lambda)+\nu_f} - \frac{1}{f_1(\boldsymbol{\eta},\lambda)+\nu_f}\right| = \left|\frac{|f_1(\boldsymbol{\eta}_n,\lambda)-f_1(\boldsymbol{\eta},\lambda)|}{(f_1(\boldsymbol{\eta}_n,\lambda)+\nu_f)(f_1(\boldsymbol{\eta},\lambda)+\nu_f)}\right|$$
$$\leq \frac{|f_1(\boldsymbol{\eta}_n,\lambda)-f_1(\boldsymbol{\eta},\lambda)|}{\nu_f^2}\,.$$

Moreover, by assumption $f_1(\boldsymbol{\eta}, \lambda)$ is continuous for all $\lambda \neq 0$ and hence uniformly continuous for $\lambda$ in any closed interval of the form $[\varepsilon, \pi]$, $\varepsilon > 0$. Consequently we can determine a value $n'$ such that for $n \geq n'$ there exists an $\varepsilon$ sufficiently small that $|f_1(\boldsymbol{\eta}_n, \lambda) - f_1(\boldsymbol{\eta}, \lambda)| < \nu_f^3$ and

$$\left|\frac{2\pi}{n}\sum_{j=1}^{\lfloor n/2 \rfloor}\frac{I(\lambda_j)}{f_1(\boldsymbol{\eta}_n,\lambda_j)+\nu_f} - \frac{2\pi}{n}\sum_{j=1}^{\lfloor n/2 \rfloor}\frac{I(\lambda_j)}{f_1(\boldsymbol{\eta},\lambda_j)+\nu_f}\right| \leq \frac{2\nu_f\pi}{n}\sum_{j=1}^{\lfloor n/2 \rfloor}I(\lambda_j)\,. \quad (A.4)$$

Using Lemma 1 in conjunction with (A.4), it follows that

$$\begin{aligned}\liminf_{n\to\infty}Q_n^{(1)}(\boldsymbol{\eta}_n) &\geq \liminf_{n\to\infty}\frac{2\pi}{n}\sum_{j=1}^{\lfloor n/2 \rfloor}\frac{I(\lambda_j)}{f_1(\boldsymbol{\eta}_n,\lambda_j)+\nu_f} \\ &\geq \lim_{n\to\infty}\left\{\frac{2\pi}{n}\sum_{j=1}^{\lfloor n/2 \rfloor}\frac{I(\lambda_j)}{f_1(\boldsymbol{\eta},\lambda_j)+\nu_f} - \frac{2\nu_f\pi}{n}\sum_{j=1}^{\lfloor n/2 \rfloor}I(\lambda_j)\right\} \\ &= \frac{\sigma_0^2}{2\pi}\int_0^\pi \frac{f_0(\lambda)}{f_1(\boldsymbol{\eta},\lambda)+\nu_f}d\lambda - \nu_f\pi\gamma_0(0)\,,\end{aligned}$$

where $\gamma_0(0)$ is the variance of the TDGP. Letting $\nu_f \to 0$ and applying Lebegue's monotone convergence theorem gives

$$\liminf_{n\to\infty}Q_n^{(1)}(\boldsymbol{\eta}_n) \geq \frac{\sigma_0^2}{2\pi}\int_0^\pi \frac{f_0(\lambda)}{f_1(\boldsymbol{\eta},\lambda)}d\lambda = Q(\boldsymbol{\eta})\,.$$

Since by definition $\boldsymbol{\eta}_1$ minimises $Q(\boldsymbol{\eta})$ it follows that $Q(\boldsymbol{\eta}_1)$ provides a lower bound to the limit inferior of $Q_n^{(1)}(\boldsymbol{\eta}_n)$ for any sequence in $\mathbb{E}_\delta^0$.

Now let $\boldsymbol{\eta}_n$ denote a sequence in $\overline{\mathbb{E}}_{\delta 1}^0 \cup \overline{\mathbb{E}}_{\delta 2}^0$ that converges to $\boldsymbol{\eta}$. Setting

$$\delta \ll \min\left\{\frac{\sigma_0^2}{4(2\pi)^2}\frac{C_l}{C_u}\frac{1}{(Q(\boldsymbol{\eta}_1)+q)}, 0.25 - 0.5(d_0 - d_1)\right\} \quad \text{where} \quad q \gg 0,$$

and applying Lemma 4 in conjunction with (A.4) implies that

$$\liminf_{n\to\infty} Q_n^{(1)}(\boldsymbol{\eta}_n) \gg Q(\boldsymbol{\eta}_1) + q.$$

Hence we can conclude that for any sequence $\boldsymbol{\eta}_n \in \overline{\mathbb{E}}_{\delta 1}^0 \cup \overline{\mathbb{E}}_{\delta 2}^0$ the criterion value $Q_n^{(1)}(\boldsymbol{\eta}_n)$ will, for all $n$ sufficiently large, exceed $Q(\boldsymbol{\eta}_1)$, which equals $\lim_{n\to\infty} Q_n^{(1)}(\boldsymbol{\eta}_1)$ by Lemma 3.

By definition of $\widehat{\boldsymbol{\eta}}_1^{(1)}$, however, $Q_n^{(1)}(\widehat{\boldsymbol{\eta}}_1^{(1)}) \leq Q_n^{(1)}(\boldsymbol{\eta}_1)$ and it follows from Lemma 3 that

$$\limsup_{n\to\infty} Q_n^{(1)}(\widehat{\boldsymbol{\eta}}_1^{(1)}) \leq \limsup_{n\to\infty} Q_n^{(1)}(\boldsymbol{\eta}_1) = Q(\boldsymbol{\eta}_1).$$

We can therefore conclude that $|Q_n^{(1)}(\widehat{\boldsymbol{\eta}}_1^{(1)}) - Q(\boldsymbol{\eta}_1)| \to 0$ almost surely and an argument by contradiction then shows that $\widehat{\boldsymbol{\eta}}_1^{(1)} \to \boldsymbol{\eta}_1$ with probability one. ∎

## A.3 Proof of Theorem 1:

In what follows we assume that the mean is known, and without loss of generality set $\mu = 0$ and suppose that the data is mean corrected.

### A.3.1 The Whittle estimator

Concentrating $Q_n^{(2)}(\sigma^2, \boldsymbol{\eta})$ with respect to $\sigma^2$ and setting $n \cdot \lfloor n/2 \rfloor = 0.5$ yields the profile (negative) log-likelihood

$$Q_n^{(2)}(\boldsymbol{\eta}) = \frac{2\pi}{2}\log\left(\frac{\widehat{\sigma}^2(\boldsymbol{\eta})}{2\pi}\right) + \frac{2\pi}{n}\sum_{j=1}^{\lfloor n/2 \rfloor} \log f_1(\boldsymbol{\eta}, \lambda_j) + \pi$$

where $\widehat{\sigma}^2(\boldsymbol{\eta}) = 2Q_n^{(1)}(\boldsymbol{\eta})$ and $Q_n^{(1)}(\boldsymbol{\eta})$ is as given in (7)). Now, following Beran (1994, p. 116), we have

$$\frac{2\pi}{n}\sum_{j=1}^{\lfloor n/2 \rfloor} \log f_1(\boldsymbol{\eta}, \lambda_j) = \frac{1}{2}\sum_{r=-\infty}^{\infty} \rho_1(\boldsymbol{\eta}, rn) \to \frac{1}{2}\int_{-\pi}^{\pi} \log f_1(\boldsymbol{\eta}, \lambda) d\lambda,$$

where the Fourier coefficients $\rho_1(\boldsymbol{\eta}, r) = \int_{-\pi}^{\pi} \log f_1(\boldsymbol{\eta}, \lambda) \exp(\imath\lambda r) d\lambda$ form a convergent series and

$$\int_{-\pi}^{\pi} \log f_1(\boldsymbol{\eta}, \lambda) d\lambda = \int_{-\pi}^{\pi} \log\left(g_1(\boldsymbol{\beta}, \lambda)|2\sin(\lambda/2)|^{-2d}\right) d\lambda$$

$$= \int_{-\pi}^{\pi} \log g_1(\boldsymbol{\beta}, \lambda) d\lambda - 2d\int_{-\pi}^{\pi} \log|2\sin(\lambda/2)| d\lambda.$$

By Assumption A.2 $\int_{-\pi}^{\pi} \log g_1(\boldsymbol{\beta},\lambda)d\lambda = 0$, and from standard results for trigonometric integrals (Gradshteyn and Ryzhik, 2007, p.583)

$$\int_{-\pi}^{\pi} \log|2\sin(\lambda/2)|d\lambda = 2\int_{0}^{\pi} \log|2\sin(\lambda/2)|d\lambda = 0\,.$$

Furthermore, since $\log f_1(\boldsymbol{\eta}, \lambda)$ is integrable, and continuously differentiable for all $\lambda \neq 0$ by Assumption A.3, $\rho_1(\boldsymbol{\eta}, n) = o(1/n)$, which implies that

$$\frac{2\pi}{n}\sum_{j=1}^{\lfloor n/2 \rfloor} \log f_1(\boldsymbol{\eta}, \lambda_j) = \sum_{r=1}^{\infty} \rho_1(\boldsymbol{\eta}, rn) = O(n^{-1}\log n)\,.$$

Hence it follows that

$$\left|Q_n^{(2)}(\boldsymbol{\eta}) - \pi \log Q_n^{(1)}(\boldsymbol{\eta}) - \pi(\log \pi + 1)\right| = O(n^{-1}\log n)\,, \tag{A.5}$$

almost surely and uniformly in $\boldsymbol{\eta}$. From this we can deduce that

$$\lim_{n\to\infty}\left|Q_n^{(2)}(\widehat{\boldsymbol{\eta}}_1^{(2)}) - \pi \log Q_n^{(1)}(\widehat{\boldsymbol{\eta}}_1^{(1)}) - \pi(\log \pi + 1)\right| = 0 \quad a.s.,$$

where $\widehat{\boldsymbol{\eta}}_1^{(1)}$ is the value of $\boldsymbol{\eta}$ that minimises the profile log-likelihood, *having first deleted the term* $2\pi\sum_{j=1}^{\lfloor n/2 \rfloor} \log f_1(\boldsymbol{\eta}, \lambda_j)/n$, namely $\widehat{\boldsymbol{\eta}}_1^{(1)} = \arg\min_{\boldsymbol{\eta}} Q_n^{(1)}(\boldsymbol{\eta})$. We are thereby lead directly to the conclusion that $\widehat{\boldsymbol{\eta}}_1^{(2)}$ and $\widehat{\boldsymbol{\eta}}_1^{(1)}$ converge, *i.e.* $\lim_{n\to\infty}\|\widehat{\boldsymbol{\eta}}_1^{(2)} - \widehat{\boldsymbol{\eta}}_1^{(1)}\| = 0$.

### A.3.2  The TML estimator

Using the argument employed by Hannan (1973, pages 134-135) in the proof of his Lemma 4, following the detailed steps given by (Brockwell and Davis, 1991, Section 10.8.2, pages 380-382) in their proof of their Proposition 10.8.3, shows that

$$\lim_{n\to\infty}\left|\frac{1}{n}\mathbf{Y}^T\boldsymbol{\Sigma}_{\boldsymbol{\eta}}^{-1}\mathbf{Y} - \frac{4\pi}{n}\sum_{j=1}^{\lfloor n/2 \rfloor}\frac{I(\lambda_j)}{f_1(\boldsymbol{\eta},\lambda_j)}\right| = 0 \quad a.s., \tag{A.6}$$

and the convergence is uniform in $\boldsymbol{\eta}$ on $\mathbb{E}_{\delta}^0$. From a theorem due to Grenander and Szego (1958, Chapter 5) we know that

$$\frac{1}{n}\log|\boldsymbol{\Sigma}_{\boldsymbol{\eta}}| = \frac{1}{2\pi}\int_{-\pi}^{\pi}\log f_1(\boldsymbol{\eta},\lambda)d\lambda + O(n^{-1}), \tag{A.7}$$

for the second term in (11). That the convergence in (A.7) is uniform in $\boldsymbol{\eta}$ is not stated in Grenander and Szegö, although it follows from the uniformity of the order relations used in their proof. Their proof depends on approximating $f_1(\boldsymbol{\eta}, \lambda)$ by trigonometric polynomials, and since $f_1(\boldsymbol{\eta}, \lambda)$ is a continuous function of $\boldsymbol{\eta}$ and $\lambda$ for all $\lambda \neq 0$ by Assumption A.3 the Stone-Weierstrass Theorem implies that $f_1(\boldsymbol{\eta}, \lambda)$ can be so approximated uniformly. It follows that

$$\lim_{n\to\to\infty}\left|Q_n^{(3)}(\sigma^2, \boldsymbol{\eta}) - \log\sigma^2 - \frac{2Q_n^{(1)}(\boldsymbol{\eta})}{\sigma^2}\right| = 0$$

almost surely, and the convergence is uniform in $\boldsymbol{\eta}$ on $\mathbb{E}_\delta^0$.

The almost sure limit of the criterion function $Q_n^{(3)}(\sigma^2, \boldsymbol{\eta})$ is therefore

$$\mathcal{Q}^{(3)}(\sigma^2, Q(\boldsymbol{\eta})) = \log \sigma^2 + \frac{2Q(\boldsymbol{\eta})}{\sigma^2},$$

uniformly in $\boldsymbol{\eta}$ on $\mathbb{E}_\delta^0$ by Lemma 3, whereas $Q_n^{(3)}(\sigma^2, \boldsymbol{\eta})$ is either arbitrarily large for $\delta$ sufficiently small or divergent on $\overline{\mathbb{E}}_{\delta 1}^0 \cup \overline{\mathbb{E}}_{\delta 2}^0$ by Lemma 4. Concentrating $\mathcal{Q}^{(3)}(\sigma^2, Q(\boldsymbol{\eta}))$ with respect to $\sigma^2$ we find that the minimum of the asymptotic criterion function is given by $\log(2Q(\boldsymbol{\eta}_1)) + 1$. Once again $\boldsymbol{\eta}_1 = \arg\min_{\boldsymbol{\eta}} Q(\boldsymbol{\eta})$ is the pseudo-true parameter for the estimator under mis-specification and we can conclude that $\lim_{n\to\infty} \widehat{\boldsymbol{\eta}}_1^{(3)} = \boldsymbol{\eta}_1$ and $\lim_{n\to\infty} \|\widehat{\boldsymbol{\eta}}_1^{(3)} - \widehat{\boldsymbol{\eta}}_1^{(1)}\| = 0$.

### A.3.3 The CSS estimator

Let $\mathbf{T}_\eta$ and $\mathbf{H}_\eta$ denote the $n \times n$ upper triangular Toeplitz matrix with non-zero elements $\tau_{|i-j|}(\boldsymbol{\eta})$, $i,j = 1,\ldots,n$, and the $n \times \infty$ reverse Hankel matrix with typical element $\tau_{n-i+j}(\boldsymbol{\eta})$, $i=1,\ldots,n$, $j=1,\ldots,\infty$, respectively. Let $\mathbf{A}_\eta = [a_{s-r}(\eta)]$ where

$$a_{s-r}(\eta) = \int_{-\pi}^{\pi} \frac{1}{f_1(\boldsymbol{\eta},\lambda)} \exp(i(s-l)\lambda) d\lambda, \quad r,s = 1,\ldots,n. \tag{A.8}$$

Then from (A.8) we can deduce that $\mathbf{A}_\eta = \mathbf{T}_\eta \mathbf{T}_\eta^T + \mathbf{H}_\eta \mathbf{H}_\eta^T$ and from (13) and (14) it follows that $Q_n^{(4)}(\boldsymbol{\eta}) = \frac{1}{n} \mathbf{Y}^T \mathbf{T}_\eta \mathbf{T}_\eta^T \mathbf{Y}$. Replacing $\boldsymbol{\Sigma}_\eta^{-1}$ by $\mathbf{A}_\eta$ in (A.6) and adapting the argument used to establish (A.6) accordingly, in a manner similar to the proof of Lemma 1, shows that

$$\lim_{n\to\infty} \left| \frac{1}{n} \mathbf{Y}^T \mathbf{A}_\eta \mathbf{Y} - \frac{4\pi}{n} \sum_{j=1}^{\lfloor n/2 \rfloor} \frac{I(\lambda_j)}{f_1(\boldsymbol{\eta},\lambda_j)} \right| = 0 \quad a.s. \tag{A.9}$$

and that the convergence is uniform in $\boldsymbol{\eta}$ on $\mathbb{E}_\delta^0$. It is also shown below that $\frac{1}{n} \mathbf{Y}^T \mathbf{H}_\eta \mathbf{H}_\eta^T \mathbf{Y} = o(1)$ when $|d| < 0.5$, $|d_0| < 0.5$ and $d_0 - d < 0.5$.

We can therefore conclude that $\left| Q_n^{(4)}(\boldsymbol{\eta}) - 2Q_n^{(1)}(\boldsymbol{\eta}) \right|$ converges to zero almost surely when $\boldsymbol{\eta} \in \mathbb{E}_\delta$, and hence that the limiting value of the criterion function $Q_n^{(4)}(\boldsymbol{\eta})$ is $2Q(\boldsymbol{\eta})$ by Lemma 3. When $\boldsymbol{\eta} \in \overline{\mathbb{E}}_{\delta 1}^0 \cup \overline{\mathbb{E}}_{\delta 2}^0$, expression (A.9) and Lemma 4, together with the equality $Q_n^{(4)}(\boldsymbol{\eta}) = \frac{1}{n} \mathbf{Y}^T \mathbf{A}_\eta \mathbf{Y} - \frac{1}{n} \mathbf{Y}^T \mathbf{H}_\eta \mathbf{H}_\eta^T \mathbf{Y}$, imply that $\liminf_{n\to\infty} Q_n^{(4)}(\boldsymbol{\eta}) \geq \lim_{n\to\infty} \frac{1}{n} \mathbf{Y}^T \mathbf{A}_\eta \mathbf{Y}$ and the CSS criterion function is either arbitrarily large for $\delta$ sufficiently small or divergent. That the pseudo-true parameter for the CSS estimator under mis-specification is $\boldsymbol{\eta}_1 = \arg\min_{\boldsymbol{\eta}} Q(\boldsymbol{\eta})$ and $\lim \widehat{\boldsymbol{\eta}}_1^{(4)} = \boldsymbol{\eta}_1$ and $\lim_{n\to\infty} \|\widehat{\boldsymbol{\eta}}_1^{(4)} - \widehat{\boldsymbol{\eta}}_1^{(1)}\| = 0$ follows directly.

It remains for us to establish that $\frac{1}{n} \mathbf{Y}^T \mathbf{H}_\eta \mathbf{H}_\eta^T \mathbf{Y} = o(1)$ in regions of the parameter space where $d_0 - d < 0.5$. Suppressing the dependence on the parameter $\boldsymbol{\eta}$ for notational simplicity, set $\mathbf{M} = \mathbf{H}\mathbf{H}^T$. Then $\mathbf{M} = [m_{ij}]_{i,j=1,\ldots,n}$ where $m_{ij} = \sum_{u=0}^{\infty} \tau_{u+n-i} \tau_{u+n-j}$, and

$$E_0[\mathbf{Y}^T \mathbf{M} \mathbf{Y}] = \text{tr}(\mathbf{M}\boldsymbol{\Sigma}_0) = \sum_{i=1}^{n} \sum_{j=1}^{n} m_{ij} \gamma_0(j-i)$$

where $\gamma_0(\tau)$, $\tau = 0, \pm 1, \pm 2, \ldots$, denotes the autocovariance function of the TDGP. Since $|\tau_k| \sim k^{-(1+d)} \mathcal{C}_\tau$, $\mathcal{C}_\tau < \infty$, the series $\sum_{k=0}^{\infty} |\tau_k|^2 \sim \mathcal{C}_\tau^2 \zeta(2(d+1))$ for all $d > -1/2$, where $\zeta(\cdot)$ denotes the Riemann zeta function, from which we can deduce that $|m_{ij}| \sim \{(n-i+1)(n-j+1)\}^{-(1+d)} \mathcal{C}'_m$ for some $\mathcal{C}'_m < \infty$. Hence on setting $r = n-i+1$ and $s = n-j+1$ we have that

$$0 \leq \sum_{i=1}^{n} \sum_{j=1}^{n} m_{ij} \gamma_0(j-i) \sim \mathcal{C}_m n^{-2(d+1)} \sum_{r=1}^{n} \sum_{s=1}^{n} |\gamma_0(r-s)|, \tag{A.10}$$

where $\mathcal{C}_m < \infty$. But $|\gamma_0(\tau)| \leq \mathcal{C}_\varrho \gamma_0(0)|\tau|^{2d_0-1}$, $\mathcal{C}_\varrho < \infty$, for all $\tau \neq 0$, and

$$\begin{aligned}
n^{-2(d+1)} \sum_{r=1}^{n} \sum_{s=1}^{n} |\gamma_0(r-s)| &\leq n^{-2(d+1)} \gamma_0(0)(n + 2\mathcal{C}_\varrho \sum_{k=1}^{n-1}(n-k)k^{2d_0-1}) \\
&\leq n^{-(2d+1)} \gamma_0(0)(1 + 2\mathcal{C}_\varrho \sum_{k=1}^{n-1} k^{2d_0-1}) \\
&\sim \frac{\gamma_0(0)}{n^{(2d+1)}} \times \begin{cases} 1 + 2\mathcal{C}_\varrho \zeta(1-2d_0), & d_0 < 0; \\ 1 + 2\mathcal{C}_\varrho \log n, & d_0 = 0; \\ 1 + 2\mathcal{C}_\varrho n^{2d_0}/2d_0, & d_0 > 0. \end{cases}
\end{aligned}$$

It follows that for all $d$ where $|d| < 0.5$

$$E_0[\mathbf{Y}^T \mathbf{M} \mathbf{Y}] \leq \frac{\mathcal{C}_m \gamma_0(0)}{n^{1-2(d_0-d)}} \times \begin{cases} 1 + 2\mathcal{C}_\varrho \zeta(1-2d_0)/n^{2d_0}, & d_0 < 0; \\ 1 + 2\mathcal{C}_\varrho \log n, & d_0 = 0; \\ 1 + \mathcal{C}_\varrho/d_0, & d_0 > 0; \end{cases}$$

We can therefore conclude that

$$Pr\left(n^{-1}\mathbf{Y}^T \mathbf{M} \mathbf{Y} > \epsilon\right) = \begin{cases} O(n^{-2(d+1)}), & 0.5 < d_0 < 0; \\ O(\log n / n^{2(d+1)}), & d_0 = 0; \\ O(n^{2(d_0-d)-2}), & 0 < d_0 < 0.5; \end{cases} \tag{A.11}$$

for all $\epsilon > 0$ by Markov's inequality. Since $\epsilon$ is arbitrary it follows that when $|d| < 0.5$ and $|d_0| < 0.5$ the almost sure limit of $n^{-1}\mathbf{Y}^T\mathbf{M}\mathbf{Y}$ is zero whenever $d_0 - d < 0.5$, by the Borell-Cantelli lemma, giving the desired result. ∎

## A.4 Proof of Theorem 2:

First note that

$$Q_N(\boldsymbol{\eta}) = \left\{\frac{\sigma_0^2 \Gamma(1-2(d_0-d))}{2\Gamma^2(1-(d_0-d))}\right\} K_N(\boldsymbol{\eta}) \tag{A.12}$$

by the same argument that gives (19). Now let $\Delta C_N(z) = \sum_{j=N+1}^{\infty} c_j z^j = C(z) - C_N(z)$. Then

$$|C(e^{i\lambda})|^2 = |C_N(e^{i\lambda})|^2 + C_N(e^{i\lambda})\Delta C_N(e^{-i\lambda}) \\ + \Delta C_N(e^{i\lambda})C_N(e^{-i\lambda}) + |\Delta C_N(e^{i\lambda})|^2$$

and the remainder term can be decomposed as $R_N = R_{1N} + R_{2N}$ where

$$R_{1N} = \left(\frac{\sigma_{\varepsilon 0}^2}{2\pi}\right) \int_0^{\pi} |\Delta C_N(e^{i\lambda})|^2 |2\sin(\lambda/2)|^{-2(d_0-d)} d\lambda \tag{A.13}$$

and

$$R_{2N} = \left(\frac{\sigma_{\varepsilon 0}^2}{2\pi}\right) \int_{-\pi}^{\pi} \Delta C_N(e^{i\lambda}) C_N(e^{-i\lambda}) |2\sin(\lambda/2)|^{-2(d_0-d)} d\lambda. \tag{A.14}$$

The first integral in (A.13) equals

$$\left\{\frac{\sigma_0^2 \Gamma(1-2(d_0-d))}{2\Gamma^2(1-(d_0-d))}\right\} \left(\sum_{j=N+1}^{\infty} c_j^2 + 2\sum_{k=N+1}^{\infty} \sum_{j=k+1}^{\infty} c_j c_k \rho(j-k)\right).$$

Because $B(z) \neq 0$, $|z| \leq 1$, it follows that $|c_j| < \mathcal{C}\zeta^j$, $j = 1, 2, \ldots$, for some $\mathcal{C} < \infty$ and $\zeta \in (0,1)$, and hence that

$$\sum_{j=N+1}^{\infty} c_j^2 < \zeta^{2(N+1)} \frac{\mathcal{C}^2}{(1-\zeta^2)}.$$

Furthermore, since $|d_0 - d| < 0.5$ Sterling's approximation can be used to show that $|\rho(h)| < \mathcal{C}'^{2(d_0-d)-1}$, $h = 1, 2, \ldots$, for some $\mathcal{C}' < \infty$. This implies that

$$\left| \sum_{k=N+1}^{\infty} \sum_{j=k+1}^{\infty} c_j c_k \rho(j-k) \right| < \sum_{r=0}^{\infty} \sum_{s=r+1}^{\infty} \mathcal{C}^2 \mathcal{C}'^{2(N+1)} \zeta^r \zeta^s (s-r)^{2(d_0-d)-1}$$

$$< \zeta^{2(N+1)} \frac{\mathcal{C}^2 \mathcal{C}'}{(1-\zeta)^2}.$$

Thus we can conclude that $R_{1N} < const. \zeta^{2(N+1)}$ where $0 < \zeta < 1$. Applying the Cauchy-Schwarz inequality to the second integral in (A.14) enables us to bound $|R_{2N}|$ by $2(\sigma_{\varepsilon 0}/\sigma)\sqrt{I_N \cdot R_{1N}}$. It therefore follows from the preceding analysis that $|R_{2N}| < const. \zeta^{(N+1)}$. Since $|R_N| \leq R_{1N} + |R_{2N}|$ and $(N+1)/\exp(-(N+1)\log\zeta) \to 0$ as $N \to \infty$ it follows that $R_N = o(N^{-1})$, as stated. The gradient vector of $Q(\boldsymbol{\eta})$ with respect to $\boldsymbol{\eta}$ is

$$\frac{\partial Q(\boldsymbol{\eta})}{\partial \boldsymbol{\eta}} = \left(\frac{\sigma_0^2}{2\pi}\right) \int_{-\pi}^{\pi} \frac{C(e^{i\lambda})}{|2\sin(\lambda/2)|^{(d_0-d)}} \frac{\partial}{\partial \boldsymbol{\eta}} \{\frac{C(e^{-i\lambda})}{|2\sin(\lambda/2)|^{(d_0-d)}}\} d\lambda$$

and substituting $C(z) = C_N(z) + \Delta C_N(z)$ gives $\partial Q(\boldsymbol{\psi})/\partial \eta_j = \partial Q_N(\boldsymbol{\eta})/\partial \eta_j + R_{3N} + R_{4N}$ for the typical $j$'th element where

$$R_{3N} = \left(\frac{\sigma_0^2}{2\pi}\right) \int_{-\pi}^{\pi} \frac{C_N(e^{i\lambda})}{|2\sin(\lambda/2)|^{(d_0-d)}} \frac{\partial}{\partial \eta_j} \{\frac{\Delta C_N(e^{-i\lambda})}{|2\sin(\lambda/2)|^{(d_0-d)}}\} d\lambda$$

and

$$R_{4N} = \left(\frac{\sigma_0^2}{2\pi}\right) \int_{-\pi}^{\pi} \frac{\Delta C_N(e^{i\lambda})}{|2\sin(\lambda/2)|^{(d_0-d)}} \frac{\partial}{\partial \eta_j} \{\frac{C(e^{-i\lambda})}{|2\sin(\lambda/2)|^{(d_0-d)}}\} d\lambda.$$

The Cauchy-Schwarz inequality now yields the inequalities

$$|R_{3N}|^2 \leq R_{1N} \left(\frac{\sigma_0^2}{2\pi}\right) \int_{-\pi}^{\pi} \frac{|C_N(e^{i\lambda})|^2}{|2\sin(\lambda/2)|^{2(d_0-d)}} \left| \frac{\partial}{\partial \eta_j} \{\log \frac{\Delta C_N(e^{-i\lambda})}{|2\sin(\lambda/2)|^{(d_0-d)}}\} \right|^2 d\lambda$$

and

$$|R_{4N}|^2 \leq R_{1N} \left(\frac{\sigma_0^2}{2\pi}\right) \int_{-\pi}^{\pi} \left| \frac{\partial}{\partial \eta_j} \{\frac{C(e^{-i\lambda})}{|2\sin(\lambda/2)|^{(d_0-d)}}\} \right|^2 d\lambda,$$

from which we can infer that $|R_{3N} + R_{4N}| \leq const. \zeta^{(N+1)} = o(N^{-1})$, thus completing the proof. ∎

## A.5   Proof of Theorem 3:

The distributions exhibited in the three cases presented in Theorem 3 correspond to those given in Theorems 1, 3 and 2 of Chen and Deo (2006), and in the following lemmas we state the properties necessary to generalise the applicability of these distributions and establish their validity under the current scenario and assumptions. Although the distributions are non-standard, the proof proceeds standardly via the use of the mean value theorem and convergence in probability of a Hessian in a neighbourhood of $\boldsymbol{\eta}_1$, plus the application to the criterion differential function of an appropriate central limit theorem.

**Lemma A.1** *Let*

$$\frac{1}{\sqrt{2\pi n}} \sum_{t=1}^{n} y_t \exp(-i\lambda t) = \xi(\lambda) = \xi_c(\lambda) - \imath \xi_s(\lambda)$$

*and set* $\mathbf{X}^T = (\xi_c(\lambda_1), \xi_s(\lambda_1), \ldots, \xi_c(\lambda_{\lfloor n/2 \rfloor}), \xi_s(\lambda_{\lfloor n/2 \rfloor})) \mathbf{F}_0^{-1/2}$ *where*

$$\mathbf{F}_0 = diag(f_0(\lambda_1), f_0(\lambda_1), \ldots, f_0(\lambda_{\lfloor n/2 \rfloor}), f_0(\lambda_{\lfloor n/2 \rfloor})).$$

*Assume that Conditions A hold. Then under Assumption A.1′ the vector* $\mathbf{X}^T$ *converges in distribution to a Gaussian random variable with zero mean and variance-covariance matrix* $\mathbf{\Omega} = \frac{1}{2}(\mathbf{I}+\mathbf{\Delta})$, $\mathbf{X}^T \xrightarrow{D} \boldsymbol{\xi} \sim N(\mathbf{0}, \mathbf{\Omega})$, *where* $\mathbf{\Delta} = [\triangle_{rc}]$, $\triangle_{rc} = O(j^{-d_0} k^{d_0-1} \log k)$ *for* $r = 2j-1$ *or* $r = 2j$, *and* $c = 2k-1$ *or* $2k$, $1 \leq j \leq k \leq \lfloor n/2 \rfloor$.

**Proof of Lemma A.1:** Assumption $(A.1')$ implies that Assumption $(A.1)$ of Lahiri (2003) holds. Since *Conditions A* imply that Assumption A.3 of Lahiri (2003) also holds, the asymptotic normality of $\mathbf{X}^T$ follows from Theorem 2.1 of Lahiri (2003). The stated covariance structure follows from Lemmas 1 and 4 of Moulines and Soulier (1999) in which the moment properties of $\xi_c(\lambda_j)$ and $\xi_s(\lambda_j)$ are derived supposing that *exact* Gaussianity holds for the sine and cosine transforms for all $n$, with bounds that are uniform with respect to $n$ for each $j = 1, \ldots, \lfloor n/2 \rfloor$. See also Corollary 5.2 of Lahiri (2003) and the discussion in Lahiri (2003, page 624). ∎

Since the limiting joint distribution of the sine and cosine transforms is Gaussian, and the sine and cosine transforms are uniformly integrable, the form of the asymptotic distribution and covariance properties of the corresponding periodogram ordinates are determined by the limit law of $\xi_c(\lambda_j)$ and $\xi_s(\lambda_j)$, $j = 1, \ldots, \lfloor n/2 \rfloor$.

**Corollary A.1** *Assume that the conditions of Lemma A.1 hold, and for each* $j = 1, \ldots \lfloor n/2 \rfloor$ *set* $Z_j = I(\lambda_j)/f_0(\lambda_j) = |\xi(\lambda_j)|^2/f_0(\lambda_j)$ *and let* $\rho_j = Cov_0[\xi_c(\lambda_j)\xi_s(\lambda_j)]/f_0(\lambda_j)$. *Then* $Z_j - \rho_j \xi_c(\lambda_j) \xi_s(\lambda_j)/f_0(\lambda_j)$ *converges in distribution to* $\frac{1}{2}\chi^2(2)(1+\triangle_{2j2j})(1-\rho_j^2)$ *where* $\chi^2(2)$ *denotes a Chi-squared random variable with two degrees of freedom. Furthermore,* $E_0[Z_j] = 1 + O(\log j/j)$, $Var_0[Z_j] = 1 + O(\log j/j)$ *and* $Cov_0[Z_j Z_k] = O(j^{-2|d_0|} k^{2|d_0|-2} \log^2 k)$ *for* $1 \leq j < k \leq \lfloor n/2 \rfloor$.

**Proof of Corollary A.1:** For $j = 1, \ldots \lfloor n/2 \rfloor$ set

$$U_j = \frac{\xi_c(\lambda_j) - \xi_s(\lambda_j)]}{\sqrt{f_0(\lambda_j)(1 + \triangle_{2j2j})(1 - \rho_j)}} \quad \text{and} \quad V_j = \frac{\xi_c(\lambda_j) + \xi_s(\lambda_j)]}{\sqrt{f_0(\lambda_j)(1 + \triangle_{2j2j})(1 + \rho_j)}}.$$

Then the Continuous Mapping Theorem implies that

$$\frac{Z_j - \rho_j \xi_c(\lambda_j) \xi_s(\lambda_j)/f_0(\lambda_j)}{(1 + \triangle_{2j2j})(1 - \rho_j^2)} = U_j^2 + V_j^2 \xrightarrow{D} \tfrac{1}{2}\chi^2(2)$$

since by Lemma A.1 $\mathbf{X}^T \xrightarrow{D} \boldsymbol{\xi} \sim N(\mathbf{0}, \mathbf{\Omega})$. Let $\mathbf{A}$ and $\mathbf{B}$ be any $\lfloor n/2 \rfloor \times \lfloor n/2 \rfloor$ symmetric selection matrices. Then $E[\boldsymbol{\xi}^\top \mathbf{A} \boldsymbol{\xi}] = \text{tr}\mathbf{\Omega}\mathbf{A}$ and $E[(\boldsymbol{\xi}^\top \mathbf{A} \boldsymbol{\xi})(\boldsymbol{\xi}^\top \mathbf{B} \boldsymbol{\xi})] = \text{tr}\mathbf{\Omega}\mathbf{A}\text{tr}\mathbf{\Omega}\mathbf{B} + \text{tr}\mathbf{\Omega}\mathbf{A}\mathbf{\Omega}\mathbf{B}$, from which the stated moments can be derived via appropriate choice of $\mathbf{A}$ and $\mathbf{B}$. Note, in particular, that $\rho_j = Cov_0[\xi_c(\lambda_j)\xi_s(\lambda_j)]/f_0(\lambda_j) = \frac{1}{2}\triangle_{(2j-1)2j} = O(\log j/j)$ and $Cov[\xi_j^2 \xi_k^2] = (E[\xi_j \xi_k])^2 = \frac{1}{4}\triangle_{2j2k}^2 = O(j^{-2|d_0|} k^{2|d_0|-2} \log^2 k)$ for $1 \leq j < k \leq \lfloor n/2 \rfloor$. ∎

The remaining steps in the proof of Theorem 3 are based on Taylor expansions of the gradient vector (or score function) of the criterion functions. For the FML estimator we have

$$\mathbf{0} = \frac{\partial Q_n^{(1)}(\boldsymbol{\eta}_1)}{\partial \boldsymbol{\eta}} + \frac{\partial^2 Q_n^{(1)}(\overline{\boldsymbol{\eta}}_1)}{\partial \boldsymbol{\eta} \partial \boldsymbol{\eta}'} (\widehat{\boldsymbol{\eta}}_1 - \boldsymbol{\eta}_1)$$

where

$$\frac{\partial Q_n^{(1)}(\boldsymbol{\eta}_1)}{\partial \boldsymbol{\eta}} = -\frac{2\pi}{n} \sum_{j=1}^{\lfloor n/2 \rfloor} \frac{I(\lambda_j)}{f_1(\boldsymbol{\eta}, \lambda_j)^2} \frac{\partial f_1(\boldsymbol{\eta}, \lambda_j)}{\partial \boldsymbol{\eta}} = -\frac{2\pi}{n} \sum_{j=1}^{\lfloor n/2 \rfloor} \frac{I(\lambda_j)}{f_0(\lambda_j)} \mathbf{w}(\boldsymbol{\eta}, \lambda_j)$$

$$\mathbf{w}(\boldsymbol{\eta}, \lambda_j) = \frac{f_0(\lambda_j)}{f_1(\boldsymbol{\eta}, \lambda_j)} \frac{\partial \log f_1(\boldsymbol{\eta}, \lambda_j)}{\partial \boldsymbol{\eta}},$$

and

$$\frac{\partial^2 Q_n^{(1)}(\boldsymbol{\eta})}{\partial \boldsymbol{\eta} \partial \boldsymbol{\eta}'} = \frac{2\pi}{n} \sum_{j=1}^{\lfloor n/2 \rfloor} \frac{I(\lambda_j)}{f_1(\boldsymbol{\eta}, \lambda_j)} \mathbf{H}(\boldsymbol{\eta}, \lambda_j),$$

$$\mathbf{H}(\boldsymbol{\eta}, \lambda_j) = 2 \frac{\partial \log(f_1(\boldsymbol{\eta}, \lambda_j))}{\partial \boldsymbol{\eta}} \frac{\partial \log(f_1(\boldsymbol{\eta}, \lambda_j))}{\partial \boldsymbol{\eta}'} - \frac{1}{f_1(\boldsymbol{\eta}, \lambda_j)} \frac{\partial^2 f_1(\boldsymbol{\eta}, \lambda_j)}{\partial \boldsymbol{\eta} \partial \boldsymbol{\eta}'},$$

and the components of $\overline{\boldsymbol{\eta}}_1$ lie on the line segment between $\widehat{\boldsymbol{\eta}}_1$ and $\boldsymbol{\eta}_1$. Existence of the Taylor expansion is justified by convexity and Assumptions A.3 and A.5.

**Lemma A.2** *Let* $\mathrm{d}Q_n^{(1)}(\boldsymbol{\eta}; \mathbf{t})$, *where* $\mathbf{t} = (t_1, \ldots, t_{l+1})^\top$, *denote the differential of* $Q_n^{(1)}(\boldsymbol{\eta})$. *Then under the assumptions of Theorem 3*

$$\frac{n}{2\pi} \left( \sum_{j=1}^{\lfloor n/2 \rfloor} \left( \mathbf{t}^\top \mathbf{w}(\boldsymbol{\eta}_1, \lambda_j) \right)^2 \right)^{-\frac{1}{2}} \left( \mathrm{d}Q_n^{(1)}(\boldsymbol{\eta}_1; \mathbf{t}) - E[\mathrm{d}Q_n^{(1)}(\boldsymbol{\eta}_1; \mathbf{t})] \right) \xrightarrow{D} Z \sim N(0, 1)$$

*for all* $\mathbf{t} \in \mathbb{R}^{l+1}$, $0 < \|\mathbf{t}\| < \infty$.

**Proof of Lemma A.2:** By Assumption A.3 the differential of $Q_n^{(1)}(\boldsymbol{\eta})$ exists and is given by $\partial Q_n^{(1)}(\boldsymbol{\eta}_1)/\partial \boldsymbol{\eta}^\top \mathbf{t}$, from which it follows that

$$\mathrm{d}Q_n^{(1)}(\boldsymbol{\eta}_1; \mathbf{t}) - E[\mathrm{d}Q_n^{(1)}(\boldsymbol{\eta}_1; \mathbf{t})] = -\frac{2\pi}{n} \sum_{j=1}^{\lfloor n/2 \rfloor} (Z_j - E[Z_j]) \mathbf{w}(\boldsymbol{\eta}_1, \lambda_j)^\top \mathbf{t}.$$

Theorem 2 of Moulines and Soulier (1999) provides a generalisation of central limit theorems for triangular arrays of martingale differences and weakly dependent sequences to similarly weighted sums of correlated variables. Replacing Moulines and Soulier's $\eta_{nj}$ by $Z_j - E[Z_j]$ and their $b_{n,j}$ by $\mathbf{w}(\boldsymbol{\eta}, \lambda_j)^\top \mathbf{t}$, recognising from Corollary A.1 that $Z_j - E[Z_j]$, $j = 1, \ldots, \lfloor n/2 \rfloor$, share the same moment structure and order of correlation as Moulines and Soulier's $\eta_{nj}$, the proof of the lemma follows Moulines and Soulier's proof of their Theorem 2 presented in Moulines and Soulier (1999, Appendix B). Conditions (i) and (ii) of Theorem 2 of Moulines and Soulier (1999) are satisfied because $C_1 \lambda_j^{-2d^*} \log \lambda_j \leq \|\mathbf{w}(\boldsymbol{\eta}, \lambda_j)\| \leq C_2 \lambda_j^{-2d^*} \log \lambda_j$ for some constants $C_1$ and $C_2$ (see Chen and Deo, 2006, expression (21) page 276) and

$$\lim_{n \to \infty} \sup_{j=1,\ldots,\lfloor n/2 \rfloor} \left( \sum_{j=1}^{\lfloor n/2 \rfloor} \left( \mathbf{t}^\top \mathbf{w}(\boldsymbol{\eta}_1, \lambda_j) \right)^2 \right)^{-1} (\mathbf{w}(\boldsymbol{\eta}_1, \lambda_j)^\top \mathbf{t})^2 = 0. \blacksquare$$

The following lemma parallels Lemma 3 of Chen and Deo (2006) and is derived in a similar fashion. The lemma and its proof are presented here for completeness.

**Lemma A.3** *Let* $\mathbb{E}_c$ *denote a compact convex subset of* $\mathbb{E}_\delta^0$ *and denote the second order differential of the FML criterion function by* $\mathrm{d}^2 Q_n^{(1)}(\boldsymbol{\eta}; \mathbf{t}) = \mathbf{t}^\top \left( \partial^2 Q_n^{(1)}(\boldsymbol{\eta})/\partial \boldsymbol{\eta} \partial \boldsymbol{\eta}' \right) \mathbf{t}$. *Then for all* $\mathbf{t}$, $\|\mathbf{t}\| < \infty$,

$$plim_{n \to \infty} \sup_{\boldsymbol{\eta} \in \mathbb{E}_c} \left| \mathrm{d}^2 Q_n^{(1)}(\boldsymbol{\eta}; \mathbf{t}) - \mathrm{d}^2 Q(\boldsymbol{\eta}; \mathbf{t}) \right| = 0.$$

*under Assumptions* A.1′ *and* A.2 − A.5.

**Proof of Lemma A.3:** By definition of the second order differential we have

$$E_0\left[\mathrm{d}^2 Q_n^{(1)}(\boldsymbol{\eta};\mathbf{t})\right] = E_0\left[\frac{2\pi}{n}\sum_{j=1}^{\lfloor n/2 \rfloor}\frac{I(\lambda_j)}{f_1(\boldsymbol{\eta},\lambda_j)}\mathbf{t}^\top \mathbf{H}(\boldsymbol{\eta},\lambda_j)\mathbf{t}\right]$$

$$=\frac{2\pi}{n}\sum_{j=1}^{\lfloor n/2 \rfloor}\frac{f_0(\lambda_j)}{f_1(\boldsymbol{\eta},\lambda_j)}\mathbf{t}^\top \mathbf{H}(\boldsymbol{\eta},\lambda_j)\mathbf{t} + \left(\frac{E_0[I(\lambda_j)]}{f_0(\lambda_j)}-1\right)\frac{f_0(\lambda_j)}{f_1(\boldsymbol{\eta},\lambda_j)}\mathbf{t}^\top \mathbf{H}(\boldsymbol{\eta},\lambda_j)\mathbf{t},$$

where $E_0[I(\lambda_j)]/f_0(\lambda_j) - 1 = O(\log j/j)$, by Corollary A.1, and $\mathbf{t}^\top \mathbf{H}(\boldsymbol{\eta},\lambda_j)\mathbf{t} = O(\log^2 \lambda_j)$ since $\sup_{\boldsymbol{\eta}} \partial \log f_1(\boldsymbol{\eta},\lambda_j)/\partial \boldsymbol{\eta}$ is of order $O(\log \lambda_j)$ by Assumptions A.2 and A.3 and $\|\mathbf{t}\| < \infty$. Thus we can conclude that

$$\frac{2\pi}{n}\sum_{j=1}^{\lfloor n/2 \rfloor}\left(\frac{E_0[I(\lambda_j)]}{f_0(\lambda_j)}-1\right)\frac{f_0(\lambda_j)}{f_1(\boldsymbol{\eta},\lambda_j)}\mathbf{t}^\top \mathbf{H}(\boldsymbol{\eta},\lambda_j)\mathbf{t} = O\left(\frac{2\pi}{n}\sum_{j=1}^{\lfloor n/2 \rfloor}\frac{\log j}{j}\lambda_j^{-2d^*}\log^2 \lambda_j\right)$$

$$= O\left(n^{2d^*-1}\sum_{j=1}^{\lfloor n/2 \rfloor}\frac{\log j}{j}j^{-2d^*}\log^2(j/n)\right)$$

$$= \begin{cases} O(n^{2d^*-1}\log^2 n), & 0 < d^* < 0.5; \\ O(n^{-1}\log^4 n), & -1.0 < d^* \leq 0, \end{cases}$$

and hence that $E_0\left[\mathrm{d}^2 Q_n^{(1)}(\boldsymbol{\eta};\mathbf{t})\right] \to \mathbf{t}^\top \frac{\partial^2 Q(\boldsymbol{\eta})}{\partial \boldsymbol{\eta}\partial \boldsymbol{\eta}'}\mathbf{t} = \mathrm{d}^2 Q(\boldsymbol{\eta};\mathbf{t})$. Similarly, setting $h(\boldsymbol{\eta};\mathbf{t},\lambda_j,\lambda_k) = \mathbf{t}^\top \mathbf{H}(\boldsymbol{\eta},\lambda_j)\mathbf{t} \cdot \mathbf{t}^\top \mathbf{H}(\boldsymbol{\eta},\lambda_k)\mathbf{t}$ and invoking Corollary A.1 once again we have

$$Var_0\left[\mathrm{d}^2 Q_n^{(1)}(\boldsymbol{\eta};\mathbf{t})\right] = \left(\frac{2\pi}{n}\right)^2 \sum_{j=1}^{\lfloor n/2 \rfloor}\sum_{k=1}^{\lfloor n/2 \rfloor}\frac{f_0(\lambda_j)}{f_1(\boldsymbol{\eta},\lambda_j)}\frac{f_0(\lambda_k)}{f_1(\boldsymbol{\eta},\lambda_k)}h(\boldsymbol{\eta};\mathbf{t},\lambda_j,\lambda_k)Cov_0\left[\frac{I(\lambda_j)}{f_0(\lambda_j)}\frac{I(\lambda_k)}{f_0(\lambda_k)}\right]$$

$$= O\left(\frac{1}{n^2}\sum_{j=1}^{\lfloor n/2 \rfloor}\sum_{k\geq j}^{\lfloor n/2 \rfloor}\lambda_j^{-2d^*}\lambda_k^{-2d^*}\log^2 \lambda_j \log^2 \lambda_k j^{-2|d_0|}k^{2|d_0|-2}\log^2 k\right)$$

$$= O\left(n^{4d^*-2}\sum_{j=1}^{\lfloor n/2 \rfloor}j^{-2(d^*+|d_0|)}\log^2(j/n)\sum_{k=1}^{\lfloor n/2 \rfloor}k^{-2(d^*-|d_0|)-2}\log^2 k \log^2(k/n)\right)$$

$$= \begin{cases} O(n^{4d^*-2}\log^4 n), & d^* + |d_0| > 0.5 \quad 0 < d^* < 0.5; \\ O(n^{-(1+2(|d_0|-d^*))}\log^5 n), & d^* + |d_0| \leq 0.5 \quad 0 < d^* < 0.5; \\ O(n^{-(1+2|d_0|)}\log^5 n), & d^* + |d_0| \leq 0.5 \quad -1 < d^* \leq 0. \end{cases}$$

It therefore follows from Markov's inequality that $\mathrm{d}^2 Q_n^{(1)}(\boldsymbol{\eta};\mathbf{t})$ converges in probability to $\mathrm{d}^2 Q(\boldsymbol{\eta};\mathbf{t})$.

Now, by the Mean Value Theorem, for any $\boldsymbol{\eta}_1$ and $\boldsymbol{\eta}_2$ in $\mathbb{E}_c$

$$\left|\mathrm{d}^2 Q_n^{(1)}(\boldsymbol{\eta}_1;\mathbf{t}) - \mathrm{d}^2 Q_n^{(1)}(\boldsymbol{\eta}_2;\mathbf{t})\right| \leq \left\|\frac{\partial\left\{\mathrm{d}^2 Q_n^{(1)}(\overline{\boldsymbol{\eta}};\mathbf{t})\right\}}{\partial \boldsymbol{\eta}}\right\| \cdot \|\boldsymbol{\eta}_1 - \boldsymbol{\eta}_2\|$$

for some $\overline{\boldsymbol{\eta}}$ between $\boldsymbol{\eta}_1$ and $\boldsymbol{\eta}_2$. Moreover,

$$E_0\left[\frac{\partial\left\{\mathrm{d}^2 Q_n^{(1)}(\boldsymbol{\eta};\mathbf{t})\right\}}{\partial \boldsymbol{\eta}}\right] = \frac{2\pi}{n}\sum_{j=1}^{\lfloor n/2 \rfloor}\left\{\frac{f_0(\lambda_j)}{f_1(\boldsymbol{\eta},\lambda_j)} + \left(\frac{E_0[I(\lambda_j)]}{f_0(\lambda_j)}-1\right)\frac{f_0(\lambda_j)}{f_1(\boldsymbol{\eta},\lambda_j)}\right\}\cdot \mathbf{k}(\boldsymbol{\eta};\mathbf{t},\lambda_j)$$

$$= \frac{2\pi}{n}\sum_{j=1}^{\lfloor n/2 \rfloor}\frac{f_0(\lambda_j)}{f_1(\boldsymbol{\eta},\lambda_j)}\mathbf{k}(\boldsymbol{\eta};\mathbf{t},\lambda_j) + \mathbf{r}_n \qquad (A.15)$$

where
$$\mathbf{k}(\boldsymbol{\eta}; \mathbf{t}, \lambda_j) = \frac{\partial \mathbf{t}^\top \mathbf{H}(\boldsymbol{\eta}, \lambda_j) \mathbf{t}}{\partial \boldsymbol{\eta}} - \mathbf{t}^\top \mathbf{H}(\boldsymbol{\eta}, \lambda_j) \mathbf{t} \frac{\partial \log f_1(\boldsymbol{\eta}, \lambda_j)}{\partial \boldsymbol{\eta}} = O(\log^3 \lambda_j)$$

and the remainder

$$\begin{aligned}
\mathbf{r}_n &= \frac{2\pi}{n} \sum_{j=1}^{\lfloor n/2 \rfloor} \left( \frac{E_0[I(\lambda_j)]}{f_0(\lambda_j)} - 1 \right) \frac{f_0(\lambda_j)}{f_1(\boldsymbol{\eta}, \lambda_j)} \mathbf{k}(\boldsymbol{\eta}; \mathbf{t}, \lambda_j) \\
&= O\left( \frac{2\pi}{n} \sum_{j=1}^{\lfloor n/2 \rfloor} \frac{\log j}{j} \lambda_j^{-2d^*} \log^3 \lambda_j \right) \\
&= \begin{cases} O(n^{2d^*-1} \log^3 n), & 0 < d^* < 0.5\,; \\ O(n^{-1} \log^5 n), & -1 < d^* \leq 0, \end{cases}
\end{aligned}$$

From Assumption A.3 and A.5 it follows that the components of the first term on the right hand side of (A.15) converge to finite constants, and hence that

$$\left| \mathrm{d}^2 Q_n^{(1)}(\boldsymbol{\eta}_1; \mathbf{t}) - \mathrm{d}^2 Q_n^{(1)}(\boldsymbol{\eta}_2; \mathbf{t}) \right| \leq C_n \|\boldsymbol{\eta}_1 - \boldsymbol{\eta}_2\|$$

where

$$C_n = \sup_{\boldsymbol{\eta} \in \mathbb{E}_c} \left\| \frac{\partial \left\{ \mathrm{d}^2 Q_n^{(1)}(\boldsymbol{\eta}; \mathbf{t}) \right\}}{\partial \boldsymbol{\eta}} \right\| = O_p(1)$$

since $\sup_n E_0 \left[ \partial \left\{ \mathrm{d}^2 Q_n^{(1)}(\boldsymbol{\eta}; \mathbf{t}) \right\} / \partial \boldsymbol{\eta} \right] < \infty$ for all $\boldsymbol{\eta} \in \mathbb{E}_c$. We can therefore conclude that $\mathrm{d}^2 Q_n^{(1)}(\boldsymbol{\eta}_1; \mathbf{t})$ is stochastically equicontinuous, and hence that

$$\mathrm{plim}_{n \to \infty} \sup_{\boldsymbol{\eta} \in \mathbb{E}_c} \left| \mathrm{d}^2 Q_n^{(1)}(\boldsymbol{\eta}; \mathbf{t}) - \mathrm{d}^2 Q(\boldsymbol{\eta}; \mathbf{t}) \right| = 0\,,$$

for all $\mathbf{t}$, $\|\mathbf{t}\| < \infty$, as required. ∎

That the FML estimator possesses the asymptotic distributions as specified in Theorem 3 now follows by replacing Lemma 5 of Chen and Deo (2006) by Lemma and Corollary A.1, Lemmas 8 and 9 by Lemma A.2, and Lemma 3 of Chen and Deo (2006) by Lemma A.3. Having made these replacements we then find that the convergence rates and asymptotic approximations given in Chen and Deo's Lemma 4 and for their Cases 1, 2 and 3 in their lemmas 6, 7, 10, 11 and 12 remain valid, thus establishing Theorem 3 for the FML estimator.

For the Whittle estimator we have, via definition of the differential and application of the chain rule, that

$$\begin{aligned}
\left| \mathrm{d} Q_n^{(2)}(\boldsymbol{\eta}; \mathbf{t}) - \frac{\mathrm{d} Q_n^{(1)}(\boldsymbol{\eta}; \mathbf{t})}{Q_n^{(1)}(\boldsymbol{\eta})} \right| &\leq |\nabla Q_n^{(2)}(\boldsymbol{\eta}; \mathbf{t}) - \mathrm{d} Q_n^{(2)}(\boldsymbol{\eta}; \mathbf{t})| + \\
&\quad |\nabla Q_n^{(2)}(\boldsymbol{\eta}; \mathbf{t}) - \nabla \log Q_n^{(1)}(\boldsymbol{\eta}; \mathbf{t})| + \\
&\quad \left| \nabla \log Q_n^{(1)}(\boldsymbol{\eta}; \mathbf{t}) - \frac{\mathrm{d} Q_n^{(1)}(\boldsymbol{\eta}; \mathbf{t})}{Q_n^{(1)}(\boldsymbol{\eta})} \right| \\
&\leq 2\epsilon \|\mathbf{t}\| + |\nabla Q_n^{(2)}(\boldsymbol{\eta}; \mathbf{t}) - \nabla \log Q_n^{(1)}(\boldsymbol{\eta}; \mathbf{t})|
\end{aligned}$$

where

$$\nabla \log Q_n^{(1)}(\boldsymbol{\eta}; \mathbf{t}) = \log Q_n^{(1)}(\boldsymbol{\eta} + \mathbf{t}) - \log Q_n^{(1)}(\boldsymbol{\eta}) \quad \text{and} \quad \nabla Q_n^{(2)}(\boldsymbol{\eta}; \mathbf{t}) = Q_n^{(2)}(\boldsymbol{\eta} + \mathbf{t}) - Q_n^{(2)}(\boldsymbol{\eta})$$

and $\epsilon \to 0$ as $\|\mathbf{t}\| \to 0$. Setting $\|\mathbf{t}\| = O(n^{-1} \log n)$, noting that (A.5) implies that the difference in differences $|\nabla Q_n^{(2)}(\boldsymbol{\eta}; \mathbf{t}) - \nabla \log Q_n^{(1)}(\boldsymbol{\eta}; \mathbf{t})| = O(n^{-1} \log n)$, we find that

$$\left| \mathrm{d} Q_n^{(2)}(\boldsymbol{\eta}; \mathbf{t}) - \frac{\mathrm{d} Q_n^{(1)}(\boldsymbol{\eta}; \mathbf{t})}{Q_n^{(1)}(\boldsymbol{\eta})} \right| \leq O(n^{-1} \log n)\,. \tag{A.16}$$

Equation (A.16) leads, in turn, to the conclusion that

$$\left| d^2 Q_n^{(2)}(\boldsymbol{\eta}; \mathbf{t}) - \frac{d^2 Q_n^{(1)}(\boldsymbol{\eta}; \mathbf{t})}{Q_n^{(1)}(\boldsymbol{\eta})} \right| \leq \left\{ \frac{dQ_n^{(1)}(\boldsymbol{\eta}; \mathbf{t})}{Q_n^{(1)}(\boldsymbol{\eta})} \right\}^2 + O(n^{-1} \log n) \quad (A.17)$$

But by Lemma 4 of Chen and Deo (2006)

$$E\left[dQ_n^{(1)}(\boldsymbol{\eta}_1; \mathbf{t})\right] = -\frac{2\pi}{n} \sum_{j=1}^{\lfloor n/2 \rfloor} E\left[\frac{I(\lambda_j)}{f_0(\lambda_j)}\right] \mathbf{w}(\boldsymbol{\eta}_1, \lambda_j)^\top \mathbf{t}$$

$$= \begin{cases} O(n^{2d^*-1} \log n), & 0 < d^* < 0.5 ; \\ O(n^{-1} \log^3 n), & -1.0 < d^* \leq 0 . \end{cases} \quad (A.18)$$

In addition,

$$Var_0\left[dQ_n^{(1)}(\boldsymbol{\eta}; \mathbf{t})\right] = O\left(\frac{1}{n^2} \sum_{j=1}^{\lfloor n/2 \rfloor} \sum_{k \geq j}^{\lfloor n/2 \rfloor} \lambda_j^{-2d^*} \lambda_k^{-2d^*} \log \lambda_j \log \lambda_k j^{-2|d_0|} k^{2|d_0|-2} \log^2 k\right)$$

$$= O\left(n^{4d^*-2} \sum_{j=1}^{\lfloor n/2 \rfloor} j^{-2(d^*+|d_0|)} \log(j/n) \sum_{k=1}^{\lfloor n/2 \rfloor} k^{-2(d^*-|d_0|)-2} \log^2 k \log(k/n)\right)$$

$$= \begin{cases} O(n^{4d^*-2} \log^2 n), & d^* + |d_0| > 0.5 \quad 0 < d^* < 0.5 ; \\ O(n^{-(1-2(d^*-|d_0|))} \log^3 n), & d^* + |d_0| \leq 0.5 \quad 0 < d^* < 0.5 ; \\ O(n^{-(1+2|d_0|)} \log^3 n), & d^* + |d_0| \leq 0.5 \quad -1.0 < d^* \leq 0 . \end{cases} \quad (A.19)$$

The asymptotic equivalence of the FML and Whittle estimators now follows since: by Lemma 3 $Q_n^{(1)}(\boldsymbol{\eta}_1)$ converges almost surely to $Q(\boldsymbol{\eta}_1) \geq \sigma_0^2 > 0$; equations (A.17), (A.18) and (A.19) imply that $|d^2 Q_n^{(2)}(\boldsymbol{\eta}_1; \mathbf{t}) - d^2 Q_n^{(1)}(\boldsymbol{\eta}_1; \mathbf{t}) / Q_n^{(1)}(\boldsymbol{\eta}_1)| = o_p(1)$; and equation (A.16) implies that

$$\frac{n}{2\pi} \left(\sum_{j=1}^{\lfloor n/2 \rfloor} \left(\mathbf{t}^\top \mathbf{w}(\boldsymbol{\eta}_1, \lambda_j)\right)^2\right)^{-\frac{1}{2}} \left| dQ_n^{(2)}(\boldsymbol{\eta}_1; \mathbf{t}) - \frac{dQ_n^{(1)}(\boldsymbol{\eta}_1; \mathbf{t})}{Q_n^{(1)}(\boldsymbol{\eta}_1)} \right|$$

$$= \begin{cases} O(n^{-2(d_0-d_1)}), & 0.25 < d_0 - d_1 < 0.5 ; \\ O((n \log n)^{-\frac{1}{2}}), & -1.0 < d_0 - d_1 \leq 0.25 . \end{cases}$$

since

$$\sum_{j=1}^{\lfloor n/2 \rfloor} \left\{\mathbf{t}^\top \mathbf{w}(\boldsymbol{\eta}, \lambda_j)\right\}^2 = O\left(\sum_{j=1}^{\lfloor n/2 \rfloor} \lambda_j^{-4d^*} \log^2 \lambda_j\right)$$

$$= \begin{cases} O(n^{4d^*} \log^2 n), & 0.25 < d^* < 0.5 ; \\ O(n \log^3 n), & -1.0 < d^* \leq 0.25 . \end{cases}$$

This establishes that Lemma A.2 also holds with $dQ_n^{(1)}(\boldsymbol{\eta}_1; \mathbf{t})$ replaced by $Q(\boldsymbol{\eta}_1) dQ_n^{(2)}(\boldsymbol{\eta}_1; \mathbf{t})$.

For the TML estimator we begin by noting that

$$\left|\frac{2\pi}{n} \sum_{j=1}^{\lfloor n/2 \rfloor} \log f_1(\boldsymbol{\eta}, \lambda_j) - \frac{\pi}{n} \log |\boldsymbol{\Sigma}_\eta|\right| = O(n^{-1} \log n),$$

and concentrating $Q_n^{(2)}(\sigma^2, \boldsymbol{\eta})$ and $Q_n^{(3)}(\sigma^2, \boldsymbol{\eta})$ with respect to $\sigma^2$ yields the inequality

$$|Q_n^{(2)}(\boldsymbol{\eta}) - \pi Q_n^{(3)}(\boldsymbol{\eta})| \leq O(n^{-1} \log n) + |\log 2Q_n^{(1)}(\boldsymbol{\eta}) - \log(2\pi/n) \mathbf{Y}^T \boldsymbol{\Sigma}_\eta^{-1} \mathbf{Y}|. \quad (A.20)$$

If we let $\mathbf{U}$ denote the $n \times n$ unitary matrix with entries $n^{-\frac{1}{2}} \exp(\imath 2\pi(r-1)(c-1)/n)$ in row $r$ and column $c$, $r, c = 1, \ldots, n$, then the off diagonal entries in $\mathbf{U}\boldsymbol{\Sigma}_\eta \mathbf{U}^*$ are of order $O(n^{-1})$, and the diagonal entries are

$$\sum_{s=-(n-1)}^{n-1} \left(1 - \frac{|s|}{n}\right) \frac{\gamma_1(s)}{\sigma^2} \exp(\imath 2\pi(j-1)s/n) \quad j = 1, \ldots, n.$$

Since $f_1(\boldsymbol{\eta}, \lambda)$ is absolutely integrable on $[-\pi, \pi]$, and by Assumptions 3 and 5 $f_1(\boldsymbol{\eta}, \lambda)$ is continuously differentiable for all $\lambda \neq 0$, from Fejer's Theorem it follows that $\mathbf{U}\boldsymbol{\Sigma}_\eta \mathbf{U}^* - \mathbf{F}_1 = O(n^{-1})$ where $\mathbf{F}_1$ equals

$$\begin{cases} \operatorname{diag}(Csf_1, f_1(\boldsymbol{\eta}, \lambda_1) \ldots, f_1(\boldsymbol{\eta}, \lambda_{\lfloor n/2 \rfloor}), f_1(\boldsymbol{\eta}, \lambda_{\lfloor n/2 \rfloor}), \ldots, f_1(\boldsymbol{\eta}, \lambda_1)), & \text{for } n \text{ odd;} \\ \operatorname{diag}(Csf_1, f_1(\boldsymbol{\eta}, \lambda_1) \ldots, f_1(\boldsymbol{\eta}, \lambda_{(n-2)/2}), f_1(\boldsymbol{\eta}, \lambda_{\lfloor n/2 \rfloor}), f_1(\boldsymbol{\eta}, \lambda_{(n-2)/2}), \ldots, f_1(\boldsymbol{\eta}, \lambda_1)), & \text{for } n \text{ even,} \end{cases}$$

and the Cèsaro sum

$$Csf_1 = \sum_{s=-(n-1)}^{n-1} \left(1 - \frac{|s|}{n}\right) \frac{\gamma_1(s)}{\sigma^2} = \begin{cases} O(n^{2d} \log n), & 0 < d < 0.5 \\ O(1), & -0.5 < d \leq 0. \end{cases}$$

*Conditions A* and Assumption A.3 imply that $\boldsymbol{\Sigma}_\eta$ and $\mathbf{F}_1$ are positive definite and it therefore follows, upon application of the Rayleigh-Ritz theorem, that

$$(2\pi/n)\left|\mathbf{Y}^T \boldsymbol{\Sigma}_\eta^{-1} \mathbf{Y} - \mathbf{Y}^T \mathbf{U}\mathbf{F}_1^{-1}\mathbf{U}^*\mathbf{Y}\right| = \left|(2\pi/n)\mathbf{Y}^T \boldsymbol{\Sigma}_\eta^{-1}\mathbf{Y} - 2Q_n^{(1)}(\boldsymbol{\eta})\right|$$
$$= n^{-1}|\mathbf{Y}^T \mathbf{R}_\eta \mathbf{Y}|$$
$$\leq n^{-1} \max_{i=1,\ldots,n} \{|\mu_i(\mathbf{R}_\eta)|\} \|\mathbf{Y}\|^2$$

where $\mu_i(\mathbf{R}_\eta)$, $i = 1, \ldots, n$, are the eigenvalues of the residual $\mathbf{R}_\eta = \boldsymbol{\Sigma}_\eta^{-1} - \mathbf{U}\mathbf{F}_1^{-1}\mathbf{U}^* = O(n^{-1})$. Evaluating the characteristic polynomial of $\mathbf{R}_\eta$ via the leading principle minors, or using the Faddeev-Leverrier method, then indicates that $|\mu_i(\mathbf{R}_\eta)|^n \leq |\mu_i(\mathbf{R}_\eta)|^{n-1} O(n^{-1})$ and the spectral radius of $\mathbf{R}_\eta$ is $O(n^{-1})$.

We can therefore use the method leading to (A.16) and (A.17) to deduce from the inequality in (A.20) that the first and second differentials satisfy $|\mathrm{d}Q_n^{(2)}(\boldsymbol{\eta}_1; \mathbf{t}) - \pi \mathrm{d}Q_n^{(3)}(\boldsymbol{\eta}_1; \mathbf{t})| = O(n^{-1} \log n)$ and $|\mathrm{d}^2 Q_n^{(2)}(\boldsymbol{\eta}_1; \mathbf{t}) - \pi \mathrm{d}^2 Q_n^{(3)}(\boldsymbol{\eta}_1; \mathbf{t})| = o_p(1)$. It therefore follows that the Whittle estimator and the TML estimator converge in distribution as

$$\frac{n}{2\pi} \left(\sum_{j=1}^{\lfloor n/2 \rfloor} \left(\mathbf{t}^\top \mathbf{w}(\boldsymbol{\eta}_1, \lambda_j)\right)^2\right)^{-\frac{1}{2}} \left|\mathrm{d}Q_n^{(2)}(\boldsymbol{\eta}_1; \mathbf{t}) - \pi \mathrm{d}Q_n^{(3)}(\boldsymbol{\eta}_1; \mathbf{t})\right|$$
$$= \begin{cases} O(n^{-2(d_0-d_1)}), & 0.25 < d_0 - d_1 < 0.5; \\ O((n \log n)^{-\frac{1}{2}}), & -1.0 < d_0 - d_1 \leq 0.25. \end{cases}$$

For the CSS estimator we have $Q_n^{(4)}(\boldsymbol{\eta}_1) = \left\{\mathbf{Y}^T \mathbf{A}_\eta \mathbf{Y} - \mathbf{Y}^T \mathbf{M}_\eta \mathbf{Y}\right\}/n$. Replacing $\boldsymbol{\Sigma}_\eta$ by $\mathbf{A}_\eta$ and adapting the argument used previously shows that $\mathbf{U}\mathbf{A}_\eta \mathbf{U}^* = 2\pi \mathbf{F}_1^{-1} + O(n^{-1})$ and hence, using (A.11), that

$$|Q_n^{(4)}(\boldsymbol{\eta}_1) - 2Q_n^{(1)}(\boldsymbol{\eta}_1)| \leq O(n^{-1}) + o_p(n^{-\frac{1}{2}}).$$

Apart from notational changes, the remaining steps in showing that the CSS and FML estimators converge in distribution are the same as those used in establishing the equivalence of the FML, Whittle and TML estimators, and are therefore omitted.

The preceding derivations imply that Lemma A.2 also holds with $\mathrm{d}Q_n^{(1)}(\boldsymbol{\eta}_1; \mathbf{t})$ replaced by $Q(\boldsymbol{\eta}_1)\mathrm{d}Q_n^{(2)}(\boldsymbol{\eta}_1; \mathbf{t})$, $\pi Q(\boldsymbol{\eta}_1)\mathrm{d}Q_n^{(3)}(\boldsymbol{\eta}_1; \mathbf{t})$ and $(1/2)\mathrm{d}Q_n^{(4)}(\boldsymbol{\eta}_1; \mathbf{t})$. As with the FML estimator, we then find that the convergence rates and asymptotic approximations given in lemmas 4, 6, 7, 10, 11 and 12 of Chen and Deo (2006) remain valid, thus establishing Theorem 3 for the Whittle, TML and CSS estimators, and hence confirming that the four estimators $\widehat{\boldsymbol{\eta}}_1^{(1)}$, $\widehat{\boldsymbol{\eta}}_1^{(2)}$, $\widehat{\boldsymbol{\eta}}_1^{(3)}$ and $\widehat{\boldsymbol{\eta}}_1^{(4)}$ are asymptotically equivalent. ∎

## B  Evaluation of Bias Correction Term

For the FML estimator we have

$$
\begin{aligned}
E_0\left(\frac{\partial Q_n^{(1)}(\boldsymbol{\eta})}{\partial \boldsymbol{\eta}}\right) &= \frac{2\pi}{n} \sum_{j=1}^{\lfloor n/2 \rfloor} E_0(I(\lambda_j)) \frac{\partial f_1(\boldsymbol{\eta}, \lambda_j)^{-1}}{\partial \boldsymbol{\eta}} \\
&= \frac{2\pi}{n} \sum_{j=1}^{\lfloor n/2 \rfloor} \left( \sum_{|k|<n} \left(1 - \frac{|k|}{n}\right) \gamma_0(k) \exp(ik\lambda_j) \right) \frac{\partial f_1(\boldsymbol{\eta}, \lambda_j)^{-1}}{\partial \boldsymbol{\eta}},
\end{aligned}
$$

where $\gamma_0(k)$ denotes the autocovariance at lag $k$ of the TDGP (see, for example, Brockwell and Davis, 1991, Proposition 10.3.1). Similarly, for the Whittle estimator we have

$$
\begin{aligned}
E_0\left(\frac{\partial Q_n^{(2)}(\sigma_\varepsilon^2, \boldsymbol{\eta})}{\partial \boldsymbol{\eta}}\right) &= \frac{4}{n} \sum_{j=1}^{\lfloor n/2 \rfloor} \frac{\partial \log f_1(\boldsymbol{\eta}_1, \lambda_j)}{\partial \boldsymbol{\eta}} \\
&\quad + \frac{8\pi}{\sigma^2 n} \sum_{j=1}^{\lfloor n/2 \rfloor} \left( \sum_{|k|<n} \left(1 - \frac{|k|}{n}\right) \gamma_0(k) \exp(ik\lambda_j) \right) \frac{\partial f_1(\boldsymbol{\eta}, \lambda_j)^{-1}}{\partial \boldsymbol{\eta}}.
\end{aligned}
$$

Differentiating the TML criterion function with respect to $\boldsymbol{\eta}$ gives

$$
\frac{\partial Q_n^{(3)}(\sigma^2, \boldsymbol{\eta})}{\partial \boldsymbol{\eta}} = \frac{1}{n}\text{tr}\boldsymbol{\Sigma}_\eta^{-1} \frac{\partial \boldsymbol{\Sigma}_\eta}{\partial \boldsymbol{\eta}} + \frac{1}{n\sigma^2} \mathbf{Y}^T \frac{\partial \boldsymbol{\Sigma}_\eta^{-1}}{\partial \boldsymbol{\eta}} \mathbf{Y},
$$

which has expectation

$$
E_0\left(\frac{\partial Q_n^{(3)}(\sigma^2, \boldsymbol{\eta})}{\partial \boldsymbol{\eta}}\right) = \frac{1}{n}\text{tr}\boldsymbol{\Sigma}_\eta^{-1} \frac{\partial \boldsymbol{\Sigma}_\eta}{\partial \boldsymbol{\eta}} - \frac{1}{n\sigma^2}\text{tr}\boldsymbol{\Sigma}_\eta^{-1} \frac{\partial \boldsymbol{\Sigma}_\eta}{\partial \boldsymbol{\eta}} \boldsymbol{\Sigma}_\eta^{-1} \boldsymbol{\Sigma}_0,
$$

where $\boldsymbol{\Sigma}_0 = [\gamma_0(|i-j|)]$ and $\sigma^2 \boldsymbol{\Sigma}_\eta = [\gamma_1(|i-j|)]$, $i, j = 1, 2, ..., n$. The criterion function for the CSS estimator in (13) can be re-written as

$$
Q_n^{(4)}(\boldsymbol{\eta}) = \frac{1}{n} \sum_{t=1}^{n} \left( \sum_{i=0}^{t-1} \tau_i y_{t-i} \right)^2 = \frac{1}{n} \sum_{t=1}^{n} \sum_{i=0}^{t-1} \sum_{j=0}^{t-1} \tau_i \tau_j y_{t-i} y_{t-j},
$$

where $\tau_i$ is as defined in (15). The gradient of $Q_n^{(4)}(\boldsymbol{\eta})$, recalling that $\tau_i = \tau_i(\boldsymbol{\eta})$, is thus given by

$$
\frac{\partial Q_n^{(4)}(\boldsymbol{\eta})}{\partial \boldsymbol{\eta}} = \frac{1}{n} \sum_{t=1}^{n} \sum_{i=0}^{t-1} \sum_{j=0}^{t-1} \left( \tau_i \frac{\partial \tau_j}{\partial \boldsymbol{\eta}} + \tau_j \frac{\partial \tau_i}{\partial \boldsymbol{\eta}} \right) y_{t-i} y_{t-j},
$$

and the expected value of the gradient is

$$
E_0\left(\frac{\partial Q_n^{(4)}(\boldsymbol{\eta})}{\partial \boldsymbol{\eta}}\right) = \frac{1}{n} \sum_{t=1}^{n} \sum_{i=0}^{t-1} \sum_{j=0}^{t-1} \left( \tau_i \frac{\partial \tau_j}{\partial \boldsymbol{\eta}} + \tau_j \frac{\partial \tau_i}{\partial \boldsymbol{\eta}} \right) \gamma_0(i-j).
$$